\documentclass[aip,amsmath,amssymb,reprint,dvipsnames]{revtex4-2}

\usepackage{graphicx}
\usepackage{dcolumn}
\usepackage{tabularx}
\usepackage{multirow}
\usepackage{bm}
\usepackage{svg}
\usepackage[utf8]{inputenc}
\usepackage[T1]{fontenc}
\usepackage{mathptmx}
\usepackage{float}

\usepackage{xcolor}
\definecolor{Blue4}{rgb}{0.2,0.2,0.6}

\usepackage[hidelinks]{hyperref} 

\usepackage{xr} 
\begin{document}

\preprint{AIP/123-QED}

\title{Memorisation and forgetting in a learning Hopfield neural network: bifurcation mechanisms, attractors and basins}


\author{Adam~E.~Essex}
\email[E-mail: ]{A.Essex@lboro.ac.uk}
\author{Natalia~B.~Janson}
\email[E-mail: ]{N.B.Janson@lboro.ac.uk}
\author{Rachel~A.~Norris}
\affiliation{Department of Mathematical Sciences, Loughborough University, Loughborough LE11 3TU, UK}
\author{Alexander~G.~Balanov}
\email[E-mail: ]{A.Balanov@lboro.ac.uk}
\affiliation{Department of Physics, Loughborough University, Loughborough LE11 3TU, UK}

\begin{abstract}

Despite explosive expansion  of artificial intelligence based on artificial neural networks (ANNs), these are employed as ``black boxes'',  since it is unclear  how, during learning, they form memories  or  develop unwanted features, 
such as 
spurious memories and catastrophic forgetting. 
Much research is available on isolated aspects of learning ANNs, but due to their high dimensionality and non-linearity, their comprehensive analysis remains a challenge. 
Here we  comprehensively analyse mechanisms of memory formation and destruction in  an $81$-neuron Hopfield network
undergoing Hebbian learning. We propose a method to analyse bifurcations in a learning ANN in terms of the stage of learning. We discover that individual memories, represented by attractors and their basins, are formed and destroyed thanks to bifurcations. We show that the same bifurcations cause both the abrupt memory loss while the ANN is trained on a single task, and the conventional catastrophic forgetting while it switches between training tasks. We point to a possible mechanism of spurious memory formation. We reveal the structure of attractor basins formed at the end of learning. 
Our strategy to analyse high-dimensional learning  ANNs  is 
principally applicable to  recurrent ANNs of any form. The demonstrated mechanisms of memory formation and of catastrophic forgetting shed light on the operation of a wider class of recurrent ANNs and could aid the development of approaches to mitigate their flaws.

\end{abstract}


\maketitle


\begin{quotation}

Artificial neural networks (ANNs)  are now abundant in   a broad range 
of pivotal technologies, including communication, security, health and  finance. 
Being rough imitations of the biological brain, they  provide a powerful  paradigm for pattern recognition, categorisation, data processing and optimisation tasks, which are conventionally attributed to artificial intelligence (AI). However, one of the major drawbacks of such AI  is the lack of transparency in the mechanisms related to memory formation, reasoning and decision-making, known as the explainability problem.  Here, by analysing the dynamics and, more specifically, the phase space of an archetypal high-dimensional Hopfield neural network, we study mechanisms of its memory formation as it learns from external stimuli. Our analysis reveals typical bifurcations, which occur in the course of training and lead to the formation of  attractors and of their basins associated with memorised categories.  We show that bifurcations play the dual role in learning: besides creating memories, they are responsible for the well-known  deficiencies of ANNs, namely, spurious memories and catastrophic forgetting.
 Our results constitute an uncommon instance of analyses and visualisations conducted on genuinely high-dimensional dynamical systems, and of a comprehensive study of the full range of dynamical phenomena involved in learning.
 They  shed light on the dynamical mechanisms of memory formation,  spurious memories and catastrophic forgetting in a wide class of recurrent ANNs, and thus provide a promising ground for further understanding of reasoning and decision-making, as well as for designing ANNs with better performance.

\end{quotation}

\section{Introduction}

Artificial neural networks (ANNs) are embedded in many aspects of modern technologies, including sound and image recognition, optimisation, data analysis, dynamics prediction and medical diagnostics \cite{LeCun:2015aa, Alzubaidi:2021aa}. Inspired by a biological intelligent system  (the brain), such networks are formed by many interconnected functional non-linear elements,  called neurons, in which the inter-neuron coupling strengths (weights) are adjusted in the course of learning. Here we focus exclusively on \emph{recurrent} ANNs, i.e. those 
modelled as generally non-linear
dynamical systems (DSs), and in what follows we often omit the descriptor ``recurrent''. 

{\bf Memory representation in NNs.} On the one hand, in machine learning  (ML) literature it is overwhelmingly assumed that in ANNs the knowledge  is contained in their weights.  On the other hand, it is widely accepted that memories or categories are represented as attractors \cite{Amari1972, Hopfield82, Hopfield_NN}, or more appropriately, as attractor basins \cite{Pineda_Hopfield_NN_bifurcation_JC88,Barack_Two_views_on_cognitive_brain_NatRevNeurosci21}.  
Both paradigms are well justified and also linked. Indeed,  for the given ANN, the weights as control parameters \emph{shape} the vector field of the respective DS and  hence fully determine the configuration of all attractors and basins. However, weights and attractors play different roles.  

Namely, in an operational ANN, individual attractor basins serve as representations of individual memories, or categories. Memory recall, or pattern recognition, is implemented by converting the input data into initial conditions, which fall into one of the attractor basins, and by letting the state of the ANN to automatically converge to the respective attractor. Memory of this type is called associative\cite{Hopfield_NN}. Importantly, for the given ANN, the collection of all weights contains information about \emph{all} memories at once. However, 
from the weights alone one cannot extract information about \emph{individual} memories e.g. by estimating some functions of weights.

 {\bf Preparing a NN for operation.} In order to bring an ANN to an operational form, one needs to appropriately shape its vector field, which can be achieved by setting the weights at appropriate values.

In the context of recurrent ANNs, the main goal of ML can be described as designing an optimal method to identify the configuration of weights, which for the given ANN would  deliver the vector field ensuring the best configuration of attractors and basins that can enable accurate pattern recognition. For different collections of patterns -- or memories -- there would be different best configurations of attractors and hence of weights. 
Due to the non-linearity of ANNs, it is usually impossible to calculate the required weights simply as functions of vectors representing patterns to be recognised (except for binary Hopfield NNs\cite{Hertz_book_91} and approximately for continuous-state Hopfield NNs \cite{Hopfield2_Dong_Connections}). 
 
 Instead,  the suitable weights are usually found by  \emph{training} the ANN. Namely,  the weights are set at some initial values, and the network is presented one by one with some typical patterns (coded by vectors) from various categories of interest. Between the consecutive time moments,  the weights are adjusted according to some pre-defined rule. Such training -- or learning -- algorithms usually incorporate the activity of the network in response to these training vectors. There are many learning algorithms \cite{Mandic_NN_recurrent_learning_algirithms_book01}, and their further development is an ongoing process.
 
{\bf Deficiencies of NNs.} Despite the extensive development of ANNs and the growing area of their application, they exhibit critical flaws. These include spurious memories \cite{Amit_89} when the emerging attractors do not represent any real categories, or  catastrophic forgetting \cite{McCloskey_catastrophic_forgetting_PLM89,French_catastrophic_forgetting_in_NNs_TCS99,Kirkpatrick_overcoming_catastrophic_forgetting_in_NNs_PNAS17,Kemker_catastrophic_forgetting_deep_NNs_AAAI18} when previously learnt categories disappear abruptly as the ANN continues to learn from a different set of inputs (i.e. to learn a different task). The latter presents a major obstacle to the realisation of lifelong learning in ANNs 
 \cite{Yoon_Lifelong_Learning_with_Dynamically_Expandable_Netw_arxive18,Parisi_lifelong_learning_NN_review_NN19}. 

 Because of the non-linearity of ANNs, the configuration of attractors and their basins cannot be determined analytically from the weight values. 
 Modifications of weights can induce bifurcations creating or destroying attractors -- and hence memories -- but these cannot be predicted analytically, either.  For the reasons above, with any learning rule,
 the mechanisms underlying the formation of memories and categories during learning are not transparent, and ANNs are employed as ``black boxes'' \cite{Ali_AI_explainable_review_IF23}. 
For this reason, common  methods for mitigating catastrophic forgetting are empirical and not directly tied to its underlying causes \cite{Ratcliff90,Robins98,Kirkpatrick_overcoming_catastrophic_forgetting_in_NNs_PNAS17,Rusu_progressive_NNs_against_catastrophic_forgetting_arxiv16,Yoon_Lifelong_Learning_with_Dynamically_Expandable_Netw_arxive18,Fayek_progressive_learning_expanding_eliminate_catastr_forgt_NN20}.

Significant efforts are being directed  towards the development of techniques to provide insights into the internal workings of ANNs, and/or to improve  their interpretability  \cite{Nauta_NN_Explainability_Methods_Review,Saleem_DNN_Explainability_Overview}. However, understanding  the mechanisms of memory formation and decision-making remains one of the major challenges of ANN research. 

{\bf Role of bifurcations in NNs.} The influence of bifurcations on learning in ANNs has been mostly seen as 
detrimental \cite{Doya_NN_bifurcations_leanring_IEEE93,Pascanu_NN_learning_bifurcations_conf13}. With this, in 1988 it was theoretically suggested that bifurcations are \emph{required} for memory formation \cite{Pineda_Hopfield_NN_bifurcation_JC88}, which implies that they should be crucial for successful learning. These two viewpoints present a paradox needing resolution. 

Despite the presumed importance of bifurcations in learning ANNs, the problem of their identification at various stages of learning has not been formulated rigorously within the \emph{standard} DS theory.  The available results that 
are described as demonstrations of  
bifurcations in learning NNs\cite{Ribeiro_NN_bifurcations_exploding_gradients_conf20,Haputhanthri_Why_NN_learn_bifurcations_24}, while being useful as indirect evidence of bifurcations, are mathematically non-rigorous, incomplete and somewhat inconclusive, as explained in Sec.~\ref{sec:challenge}. Rigorous bifurcation analysis of learning ANNs in terms of the stage of learning is still lacking.

Ideally, learning algorithms  should take account of bifurcations and either use them, or avoid them as required. However, before such algorithms could be developed, it is necessary to reveal bifurcation mechanisms operating within the existing learning algorithms in standard ANNs. 

{\bf  Goals and results.} In this paper we address the \emph{mechanisms} of memory formation and destruction, including those of spurious memory formation and catastrophic forgetting, in a paradigmatic recurrent ANN -- continuous-time Hopfield NN \cite{Hopfield_NN} --  undergoing foundational Hebbian learning \cite{Hebb49,Gerstner_Hebbian_learning_BC02}. 
 For this NN, we resolve the paradox mentioned above by establishing that memory formation and destruction are both caused by, and are manifestations of, bifurcations occurring in the course of learning. 

In connection with memory loss, we consider two different settings, which in much ML literature are conventionally regarded as conceptually distinct, namely, while the NN learns a \emph{single} task, and when after learning one task it \emph{switches} to learning a new task. In the latter scenario, the abrupt loss of previously acquired memory is conventionally called ``catastrophic forgetting'' \cite{McCloskey_catastrophic_forgetting_PLM89,French_catastrophic_forgetting_in_NNs_TCS99,Kirkpatrick_overcoming_catastrophic_forgetting_in_NNs_PNAS17,Kemker_catastrophic_forgetting_deep_NNs_AAAI18}. 
NNs learning a \emph{single} task were observed to experience memory loss similar to the conventionally understood catastrophic forgetting \cite{Toneva_NN_forgetting_single_task_arXiv19,Eisenmann_NN_discrete_2-6_neuron_bifurcation_NeurIPS23}.
However, these two cases of memory loss have been deemed conceptually disconnected 
because they were linked to different settings where they were observed, while their exact causes -- or mechanisms -- have been unknown. With this, if two phenomena with similar manifestations have identical causes, they would be conceptually connected. Here we discover that in both settings, the Hopfield NN with Hebbian learning forgets previously acquired memories abruptly and as a result of the same bifurcations, and hence establish the conceptual connection between the two cases of memory loss.

In order to provide firm evidence of bifurcations as underlying causes of memorisation and forgetting,
we rigorously formulate the problem of identifying bifurcations in a \emph{learning} NN. 
We also reveal the basins of attraction formed by the end of  training on the first task, and hence the shape of memory imprints developed in the NN. 

{\bf  Structure.} The paper has the following structure. In Sec.~\ref{sec:model} we describe the ANN 
analysed 
in this study. Section~\ref{sec:attr} reports  evidence 
of bifurcations occurring in the course of learning and leading to the formation  or the disappearance of memories. In Sec.~\ref{sec:bif} a systematic bifurcation analysis is performed  in terms of the stage of learning,
 which reveals how memories are formed.  In Sec.~\ref{sec:memfog} we formulate a hypothesis about the bifurcation mechanisms of 
 forgetting. In Sec.~\ref{sec:forget} we verify this hypothesis for the NN learning a single task,  
and in Sec.~\ref{sec_cat} for what is conventionally known as catastrophic forgetting. Section~\ref{sec:bas} discusses the structure of the boundaries of attraction basins and their relation to  the network memory. 
The results are summarised in Sec.~\ref{sec:summary}, and the broader significance of our findings is discussed in Sec.~\ref{sec:disc}.

\section{Hopfield neural Network with Hebbian learning}
\label{sec:model}

Hopfield NNs \cite{Hopfield82, Hopfield_NN} are a type of ANN  capable of pattern recognition, which  feature the so-called associative memory.   Despite the recent proliferation of deep (multi-layer) NNs, a single-layer Hopfield NN is an important paradigm of a NN and remains highly relevant to date \cite{Yu2020,Krotov:2023aa}. It consists of a set of coupled neurons (functional units), which are connected to each other with certain strengths, usually referred to as  weights.

Here  we use the NN model in the same form as in the foundational paper \cite{Hopfield2_Dong_Connections}:
\begin{eqnarray}
\label{eqn:eq2.1}
\frac{dx_{i}}{dt} &=& -x_{i}+g\sum_{j=1}^{N} F\left(\omega_{ij}\right)F\left(x_{j}\right)+AI_{i}(t), \\
\label{eqn:eq2.4}
\frac{d\omega_{ij}}{dt} &=& \frac{1}{B_{ij}}\left( -\omega_{ij} + F\left(x_{i}\right)F\left(x_{j}\right)\right),
\end{eqnarray}
where  $x_{i}(t)$ $\left(i=1,\dots ,N\right)$ is the state of $i^{\textrm{th}}$ neuron at time $t$, $\omega_{ij}(t)$ is the weight of connection between $i^{\textrm{th}}$ and $j^{\textrm{th}}$ neurons, and $I_i(t)$ is external input to $i^{\textrm{th}}$ neuron, which is used as the training stimulus applied in the course of NN learning.  Following  Ref.~\cite{Hopfield2_Dong_Connections}, we consider the case without self-connections $\omega_{ii}$$=$$0$ and with symmetric couplings $\omega_{ij}=\omega_{ji}$ for all $i,j$$=$$1,$$\dots$$,N$. We note, that in Ref.~\cite{Hopfield2_Dong_Connections} the non-linear activation function $F$ was piecewise-linear, thus making the DS (\ref{eqn:eq2.1})--(\ref{eqn:eq2.4}) non-smooth. We are interested in the study of the Hopfield NN as a smooth DS, and for this reason we choose $F$ to take the form of a smooth sigmoid function from Ref.~\cite{Hopfield_NN}:
\begin{equation}
\label{eqn:eq2.2}
F\left(x\right) = \left(\frac{2}{\pi}\right)\arctan\left(\frac{\lambda \pi x}{2}\right).
\end{equation}
The parameter $g$ is the coupling gain parameter, and  $A$ is the strength of the external input. 

The network learns from the input signal $\mathbf{I}(t)$$=$$(I_1(t),\ldots,I_N(t))$ containing a sequence of $N$-dimensional vectors, each coding some pattern, by adjusting its weights $\omega_{ij}$ according to some built-in rule. Often this is the Hebbian learning rule governed by (\ref{eqn:eq2.4}). The rate of learning is regulated by parameter $B_{ij}$. 

Once the patterns are learned, the weights $\omega_{ij}$ are fixed at their latest values, and $\mathbf{I}(t)$ is set to zero. After that the network (\ref{eqn:eq2.1}) can recognise a new pattern as belonging to some stored category, or ``recall a memory''. This is achieved by setting initial conditions of (\ref{eqn:eq2.1}) to a vector representing this new pattern and observing how the state $\mathbf{x}$$=$$(x_1,\ldots,x_N)$ converges to an attractor to whose basin the initial conditions belong. The respective attractor is assumed to represent the most typical element of the given category, or the undistorted memory. 

Hopfield NNs have some important properties.    They are robust to noise and can recall memories from cues, which are distorted or incomplete versions of the memorised patterns.  They also have the property of a content-addressable memory \cite{Koberle_NN_content-addressable_memory_CPC89}, meaning that the network can recall a stored pattern based on its content, rather than its physical address as in random-access memory (RAM) of conventional computers. This is achieved by \emph{identifying} the content with location: any pattern is coded by a state vector $\mathbf{x}$, but  instead of the physical location within a real device, any state vector has its unique location  in the phase space of the NN.  When recalling a memory from a cue, the NN automatically goes through the states in its phase space towards the required attractor. 

Operation of the Hopfield NN is often explained in terms of an energy function \cite{Hopfield2_Dong_Connections}, to whose local minimum the state of the NN converges with time. Namely, each state of the network is associated with an energy level, and the network evolves from the states with higher energy to those with lower energy. The state with the locally lowest energy represents an attractor, which in the Hopfield NN is the stable fixed point. 

Remarkably, if  in (\ref{eqn:eq2.2}) $\lambda \to \infty$, thus making $F(x)$ a Heaviside step function, for $x_{j}$$ \neq$$0$ with $j$$=$$1,$$\dots$$,N$,  Eq.~(\ref{eqn:eq2.1}) can be rewritten as
\begin{equation}
\label{eqn:eq2.6}
\frac{dx_{i}}{dt} = -\frac{\partial V}{\partial x_{i}},
\end{equation}
where 
\begin{equation}
\label{eqn:eq2.5}
V = \sum_{i=1}^{N} \left[\frac{\left(x_{i}-AI_{i}\left(t\right)\right)^{2}}{2} - g{x_{i} \sum_{j=1}^{N}{F\left(\omega_{ij}\right)F\left(x_{j}\right)}}\right].\end{equation}
In this form, Eq.~(\ref{eqn:eq2.1}) has a  mechanical interpretation as a model describing the motion of an overdamped  particle with location  $\mathbf{x}(t)$ in the potential energy landscape $V$$=$$V(\mathbf{x},t)$. Thus, the shape of $V$  changes with time, as determined by the dynamics of the weights $\omega_{ij}$ and stimulus $I_i$. In that case, the model (\ref{eqn:eq2.1})--(\ref{eqn:eq2.4}) can be regarded as a particular case of a dynamical system with plastic self-organising vector field \cite{Janson:2017aa, Janson_NonAutonomous_Attractors}. Note, that even in the form Eq.~(\ref{eqn:eq2.6}) the NN is not a gradient system because its energy landscape explicitly depends on time due to the presence of the  stimuli $I_i(t)$.

\begin{figure}
 \includegraphics[width=0.5\textwidth]{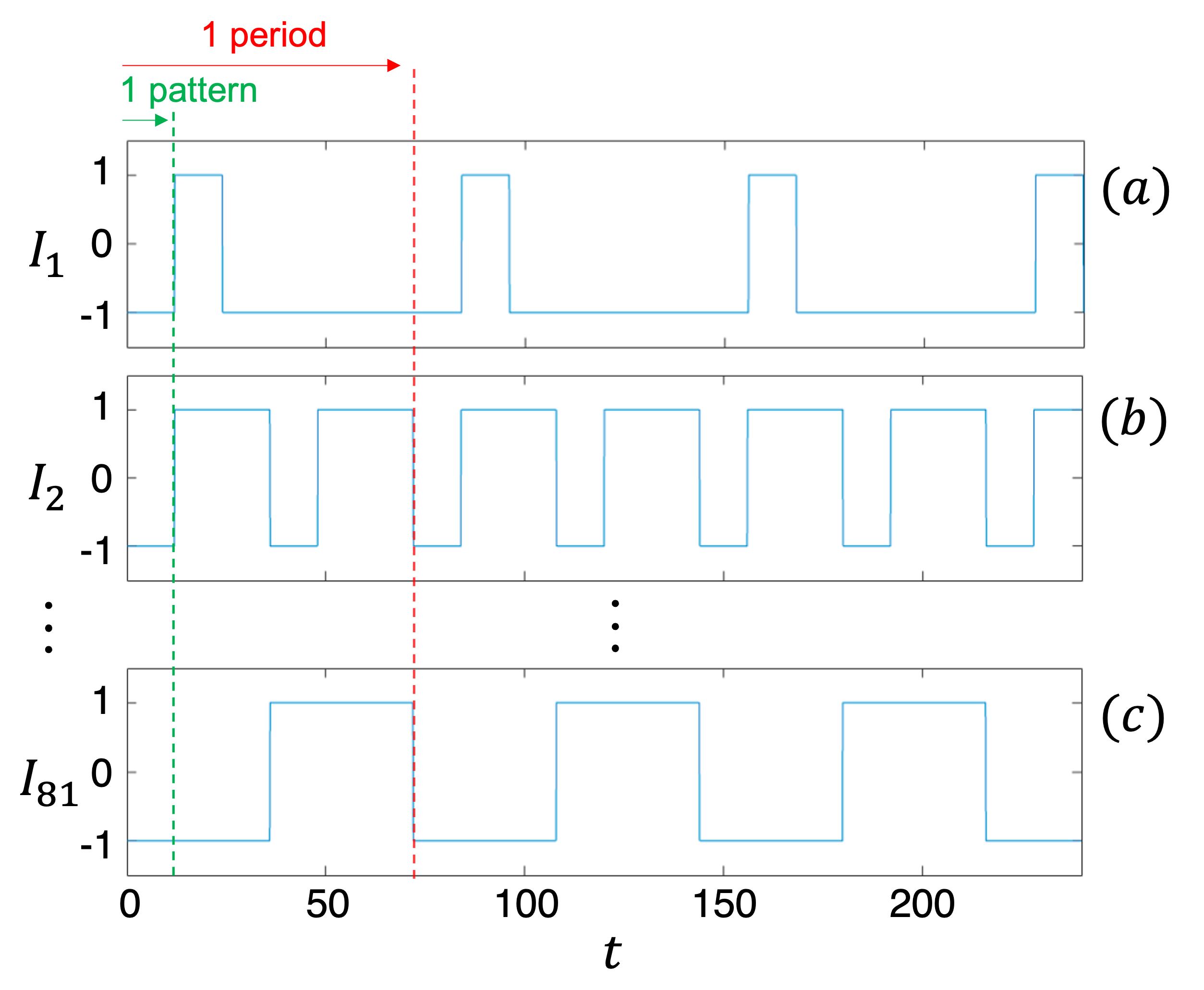}
  \caption{Illustration of construction of input signals $I_{i}(t)$ to the NN (\ref{eqn:eq2.1})--(\ref{eqn:eq2.4}) with  $N$$=$$81$ from training Set~1 of vectors $\mathbf{I}^k$ ($k$$=$$1,$$\ldots$$,6$) given in Tabs.~\ref{supp-tab:tab1-1}--\ref{supp-tab:tab1-2} of Supplementary Note. Vertical green dashed line indicates the time  $t_s$ during which a single vector $\mathbf{I}^k$ is applied to the NN. Red dashed line indicates one full period of stimulus $\mathbf{I}(t)$ (duration of ``training epoch''), which is equal to $6t_s$. Panels (a), (b) and (c) show inputs $I_{i}(t)$ to the $1^{\textrm{st}}$, $2^{\textrm{nd}}$ and $81^{\textrm{st}}$ neurons, respectively. 
  }
  \label{fig:InputSig}
\end{figure}

The local minima of $V$ occur at the stable fixed points, which either by themselves\cite{Amari1972, Hopfield82, Hopfield_NN}, or together with their attraction basins \cite{Barack_Two_views_on_cognitive_brain_NatRevNeurosci21}, represent stored patterns, or memories that the network has learned.  For (\ref{eqn:eq2.1})--(\ref{eqn:eq2.4})  with  $N$$=$$81$, this is illustrated in Fig.~\ref{supp-fig:Energy}(a)--(b) of Supplementary Note, showing two different cross-sections (coloured surfaces) of the same post-training $V$   by two different three-dimensional flat surfaces  in $(81$$+$$1)$-dimensional space, each going through a different attractor. The respective attractors (black circles) are visible at the bottoms of cross-sections of $V$. 
Notably,  $V$ could in some sense be regarded as an energy function discussed in Refs.~\cite{Hopfield82, Hopfield_NN}.  

We study (\ref{eqn:eq2.1})--(\ref{eqn:eq2.4}) with the following parameter values:  $N$$=$$81$, $A$$=$$30$, $B_{ij}$$=$$B$$=$$300$ $(\forall i,j)$, $g$$=$$0.3$  (as in Ref.~\cite{Hopfield2_Dong_Connections}) and $\lambda$$=$$1.4$ (as in Ref.~\cite{Hopfield_NN}). 

Below we describe the procedure, which in ML literature is referred to as  training an ANN on a \emph{single} task, i.e. on a single fixed collection of input patterns, each represented by a vector. 
 By analogy with Ref.~\cite{Hopfield2_Dong_Connections},  for (\ref{eqn:eq2.1})--(\ref{eqn:eq2.4})
we design the input signal  $\mathbf{I}(t)$ corresponding to a single task from a training set of six $N$-dimensional vectors, each assumed to code some non-specified pattern. The components $I_i^k$ ($k$$=$$1,$$\ldots$$,6$)  of these vectors are chosen at random from the set $\{1, -1\}$, so that $1$ or $-1$ have equal probability of occurrence, and the values of the components are statistically independent of each other. Two different training sets, Sets~1 and 2, used in our study for the network of $N$$=$$81$ neurons are given in Tabs.~\ref{supp-tab:tab1-1}--\ref{supp-tab:tab1-2} of Supplementary Note.

Starting from $k$$=$$1$, the signal $\mathbf{I}(t)$ is equal to a vector $\mathbf{I}^k$$=$$(I_1^k,\ldots,I_{N}^k)$ from the training set during the initial $t_{s}$$=$$12$ time units. For the subsequent $t_{s}$ time units, $\mathbf{I}(t)$ becomes equal to $\mathbf{I}^{k+1}$. In the same manner $\mathbf{I}(t)$ goes through all vectors $\mathbf{I}^k$ until $k$ reaches $6$. After that, the signal known as the ``training epoch'' in 
ML literature \cite{Alzubaidi_training_epoch_NN_JBD21}, and here lasting $6 t_s$ time units,   is repeated periodically throughout the maximal training time of $6000$ time units. This way, inputs to individual neurons take the form of telegraph signals shown in Fig.~\ref{fig:InputSig}. 

At the start of learning, the initial conditions (ICs) of (\ref{eqn:eq2.1})--(\ref{eqn:eq2.4})  were randomly and uniformly distributed in $[-1,1]$ for all $x_i$, and in $[-0.01,0.01]$ for all $\omega_{ij}$, the latter ensuring the absence of pre-existing memories.  

\begin{figure*}
 \raggedright 
\includegraphics[width=0.8\textwidth]{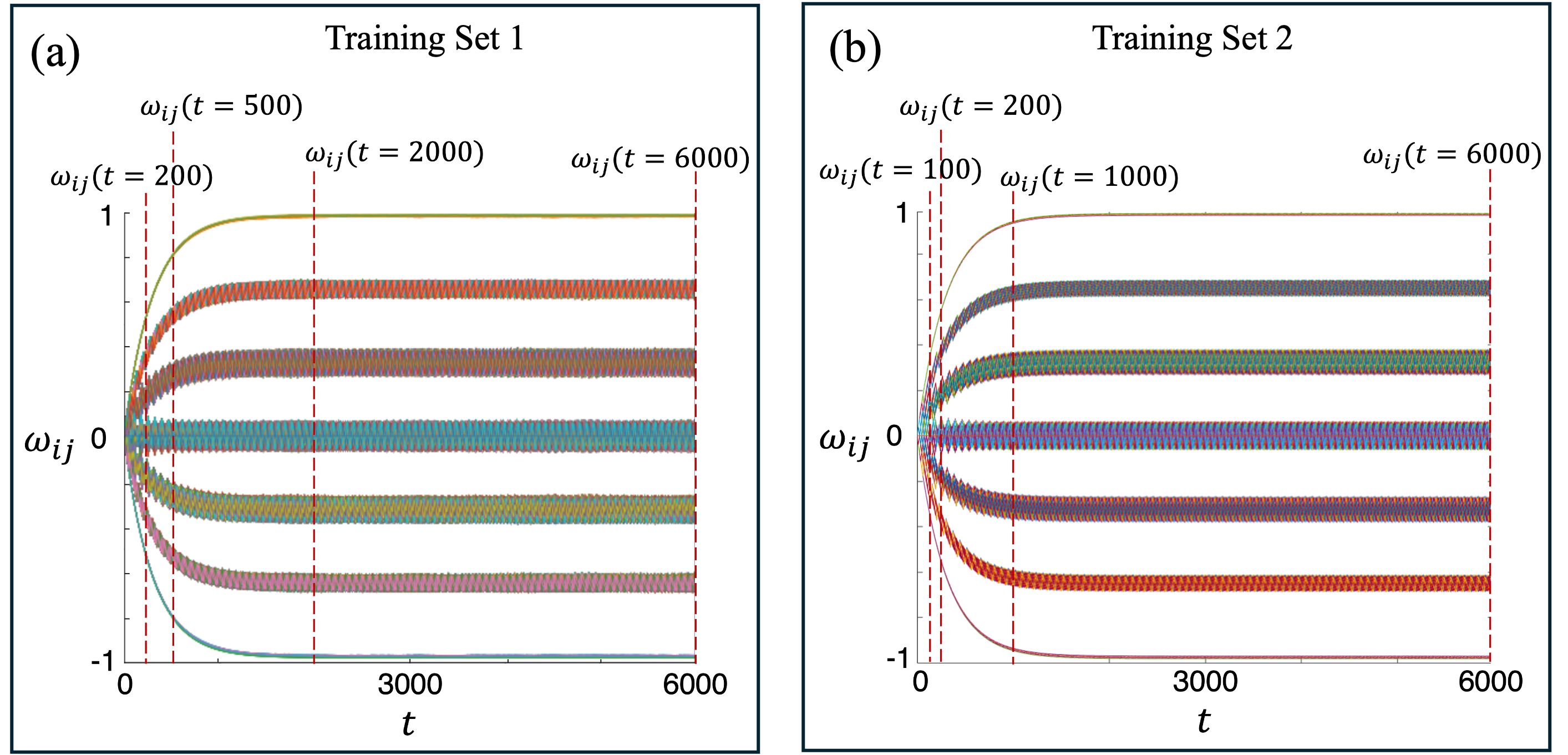}
\caption{Evolution of connection weights $\omega_{ij}(t)$ (solid lines in various colours) during learning by the NN (\ref{eqn:eq2.1})--(\ref{eqn:eq2.4}) with $N$$=$$81$, $A$$=$$30$, $B_{ij}$$=$$B$$=$$300$ $(\forall i,j)$, $g$$=$$0.3$ and $\lambda$$=$$1.4$. Panels show learning from two different training sets (see Tabs.~\ref{supp-tab:tab1-1}--\ref{supp-tab:tab1-2}): (a) Set~1 and (b) Set~2. Vertical dashed lines mark four different stages of learning corresponding to times $t$ given in brackets, at which instantaneous weights $\omega_{ij}$ are collected to reveal memories formed, as illustrated in Fig.~\ref{fig:fig3.1}. 
}
\label{fig:omega_t}
\end{figure*}

The training duration of $6000$ time units was chosen experimentally by ensuring that the weights $\omega_{ij}$ converge to small-amplitude oscillations about some fixed values (see Fig.~\ref{fig:omega_t}) for a sufficiently long time before the end of learning. Small oscillations of individual $\omega_{ij}$, which occur during the learning phase, are due to disturbance of (\ref{eqn:eq2.1})--(\ref{eqn:eq2.4}) by a periodic signal $\mathbf{I}(t)$, and their period is equal to $6 t_s$$=$$72$. 

~\vspace{-7mm}

\section{Memory formation: observation of attractors}
\label{sec:attr}

After the NN completes its training and becomes ready for operation, it takes the form of an autonomous dynamical system (DS):
\begin{equation}
\label{eqn:eq3.1}
  \frac{dx_{i}}{dt'} = -x_{i} + g\sum_{j=1}^{N}{F\left(\omega_{ij}(t) \right)F\left(x_{j}\right)} = u_{i}(\mathbf{x},t).
\end{equation} 
Equation~(\ref{eqn:eq3.1}) is obtained from (\ref{eqn:eq2.1}) by setting $A$$=$$0$,  introducing a new time $t'$ different from $t$ and corresponding to the memory retrieval phase, and fixing $\omega_{ij}$ at the values they took at some time $t$ of the learning system  (\ref{eqn:eq2.1})--(\ref{eqn:eq2.4}). In  (\ref{eqn:eq3.1}) we use the notation $\omega_{ij}(t)$, rather than  $\omega_{ij}$, to emphasise that all the weights are functions of the same parameter $t$, which is different from the time $t'$ of the given system. 
 
Note, that the model (\ref{eqn:eq2.1})--(\ref{eqn:eq2.4}) of a learning NN represents a non-autonomous DS. In such systems it can be possible to define and detect non-autonomous attractors  \cite{Kloeden_b2020}. However, doing so is not practical here, since we wish to study memories existing in the \emph{autonomous} system (\ref{eqn:eq3.1}) at various stages of learning, i.e. when (\ref{eqn:eq3.1}) is considered with various sets of \emph{fixed} $\omega_{ij}$ occurring at various times $t$ during the learning phase described by (\ref{eqn:eq2.1})--(\ref{eqn:eq2.4}). Such memories will be associated with conventional attractors and their basins. 

In (\ref{eqn:eq3.1}), $\mathbf{u}$$=$$(u_1,$$\ldots$$,u_N)$ is the phase velocity vector field of the NN formed at the given stage  $t$ of learning, which is parametrised by the set of $\omega_{ij}(t)$. For the NN of the size $N$$=$$81$, the time evolution of all $\omega_{ij}$ during the learning  phase is illustrated  in Fig.~\ref{fig:omega_t}. Namely, Figs.~\ref{fig:omega_t}(a) and (b) illustrate training with two different training sets, Set~1 and Set~2  (Tabs.~\ref{supp-tab:tab1-1}--\ref{supp-tab:tab1-2}), respectively. The two graphs are are similar to each other and to the graph of weights in Ref.~\cite{Hopfield2_Dong_Connections}, despite the fact that in the latter paper a different form of non-linear activation function $F$ was used as compared to (\ref{eqn:eq2.2}), namely, a non-smooth piece-wise linear function. We observed similar behaviours of weights also in smooth Hopfield NN in an alternative form (\ref{supp-eqn:Hop_Alt})--(\ref{supp-eqn:Hop_Alt3}) described in Supplementary Note, and also in numerous preliminary studies of both versions of smooth Hopfield NN with a variety of parameter values.

It is worth noting that, as follows from the properties of standard non-autonomous DSs, as long as the periodic stimulus $\mathbf{I}(t)$ is present in (\ref{eqn:eq2.1})--(\ref{eqn:eq2.4}), its variables -- including the weights -- do not and cannot converge to fixed values at any time $t$. This is clearly visible in Fig.~\ref{fig:omega_t} where most weights keep oscillating with relatively large amplitudes at any stage of learning. Note, that the weights close to $0$ and $1$ oscillate, too, only with much smaller amplitudes relative to others.

\begin{figure*}
\includegraphics[width=1.0\textwidth]{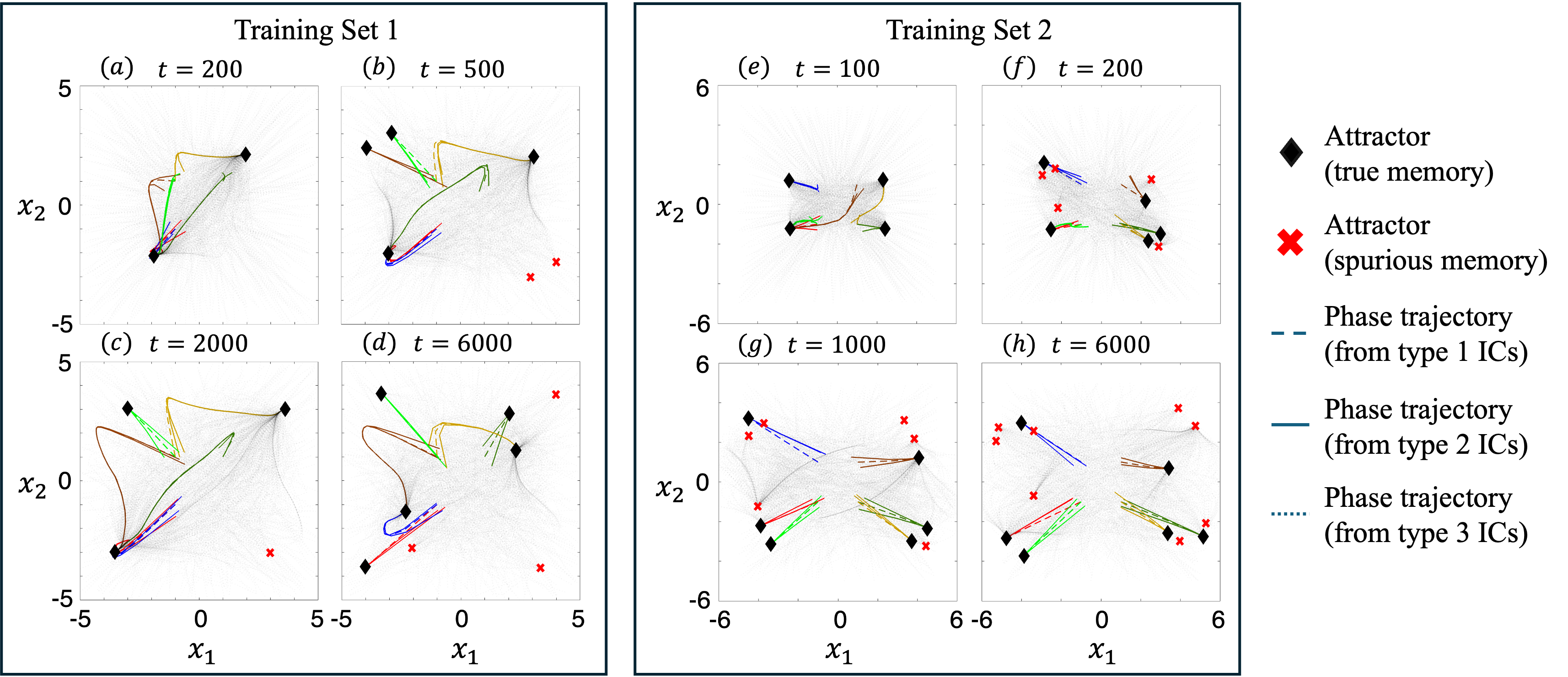}
\caption{Memories formed at various stages of learning  by the NN (\ref{eqn:eq2.1})--(\ref{eqn:eq2.4}) with $N$$=$$81$, $A$$=$$30$, $B_{ij}$$=$$B$$=$$300$ $(\forall i,j)$, $g$$=$$0.3$  and $\lambda$$=$$1.4$. Panels illustrate stages of learning from two different training sets (see Tabs.~\ref{supp-tab:tab1-1}--\ref{supp-tab:tab1-2}): (a)--(d) Set~1 and (e)--(h) Set~2. \\
Panels display projections onto the $(x_1,x_2)$-plane of phase portraits  of the NN (\ref{eqn:eq3.1}) with sets of fixed values of $\omega_{ij}$ taken at various stages of learning corresponding to times $t$ of  (\ref{eqn:eq2.1})--(\ref{eqn:eq2.4}) indicated above each panel,  also marked in Fig.~\ref{fig:omega_t}. Each panel shows attractors associated with true memories of input vectors $\mathbf{I}^k$ (black diamonds), attractors representing spurious memories (red crosses) not associated with any $\mathbf{I}^k$,  phase trajectories originating from six vectors $\mathbf{I}^k$ and their vicinities (Type 1 and 2 ICs, dashed and solid lines with one colour corresponding to one $\mathbf{I}^k$), and phase trajectories launched from random ICs (Type 3 ICs, grey dashed lines). 
}
\label{fig:fig3.1}
\end{figure*}

Despite the similarity of graphs in Figs.~\ref{fig:omega_t}(a) and (b),  for different training sets, both the processes of learning, and the post-learning memories formed, are very different. This is illustrated by Fig.~\ref{fig:fig3.1}, where for two training sets, Set~1 (a)--(d) and Set~2 (e)--(h), the projections of phase portraits of (\ref{eqn:eq3.1}) onto the $(x_1,x_2)$-plane are given at various stages of learning. The respective stages are marked  in Fig.~\ref{fig:omega_t}(a) and (b), see values of $t$ in brackets next to vertical dashed lines. 

It is desirable that, as a result of learning, attractors representing vectors $\mathbf{I}^k$  either coincide with them, or are  found at predictable locations, 
ideally being directly proportional to  $\mathbf{I}^k$. However, for continuous-time and continuous-space NNs this does not happen because of their non-linearity.  This is clearly seen in Figs.~\ref{fig:fig3.1}(d) and (h): among various attractors formed by the end of learning (black diamonds and red crosses), none have components equal to $\pm 1$, $\pm A$ or $\pm C$, where $C$ is some constant. Thus, to associate attractors with vectors $\mathbf{I}^k$, a different approach is needed, as described below. 

In Fig.~\ref{fig:fig3.1} attractors associated with input vectors $\mathbf{I}^k$ (``true memories'') are shown as black diamonds, and attractors representing spurious memories as red crosses. To reveal attractors and to distinguish between them, we considered three types of ICs specified  below. Descriptions in brackets refer to trajectories launched from these ICs. 

\begin{enumerate}

    \item Type~1: Identical to six vectors $\mathbf{I}^k$ from the training set  (dashed coloured lines). 
   
    \item Type~2: Within small vicinities of six vectors  $\mathbf{I}^k$ from the training set. Namely, to each $\mathbf{I}^k$, small vectors were added, whose components were randomly and uniformly distributed in $[-0.5,0.5]$ (solid lines coloured as in Type 1 for respective $\mathbf{I}^k$). 
    
      \item  Type~3: 1000 randomly and uniformly distributed within an $N$-dimensional hyper-cube of side length 10 and the centre at the origin (grey dashed lines). 
    
     \end{enumerate}
     
If (\ref{eqn:eq3.1}) converges to the given attractor from Type~1 and 2  ICs corresponding to a single $\mathbf{I}^k$, we deem this attractor and its basin a true memory of the respective $\mathbf{I}^k$. 
If the system converges to the same attractor from multiple input patterns $\mathbf{I}^k$ and their respective neighborhoods, we interpret this attractor and its basin as representing a blended memory of the corresponding $\mathbf{I}^k$s.  Any attractor whose basin does not encompass any $\mathbf{I}^k$ or its vicinity is considered a spurious memory.

Type~2 ICs are important for identification of true memories, since to ensure that their recall is robust to noise, the respective attractor basins should not be negligibly small. 

Let us observe stages of memory formation as the NN (\ref{eqn:eq2.1})--(\ref{eqn:eq2.4}) is processing the repeatedly applied six vectors from Set~1. 
 By $t$$=$$200$ (Fig.~\ref{fig:fig3.1}(a)) the NN develops only two attractors (black diamonds).  The one on the right is associated with the true memory of $\mathbf{I}^3$, since it attracts phase trajectories from the respective Type~1 and 2 ICs (yellow lines).  The one on the left is associated with a blended memory of five patterns, since it attracts all trajectories from ICs of Types~1 and 2 for $k$$\ne$$3$ (other coloured lines). 

By  $t$$=$$500$ (Fig.~\ref{fig:fig3.1}(b)), in addition to attractors present at $t$$=$$200$ (see Fig.~\ref{fig:fig3.1}(a)), four new attractors are formed. With this, the old attractors have moved away from each other (compare black diamonds along the main diagonal in (b) and (a)), and one of them is still associated with $\mathbf{I}^3$. Among the four new attractors, two are spurious memories (red crosses). However, two of the new attractors with basins become true memories of $\mathbf{I}^5$ and $\mathbf{I}^6$, as evidenced by the trajectories from the respective ICs (brown and light-green lines). Thus, one of the old attractors is now associated with a blend of $\mathbf{I}^1$, $\mathbf{I}^2$ and $\mathbf{I}^4$ (see red, dark green and blue lines, respectively). 

By $t$$=$$2000$ (Fig.~\ref{fig:fig3.1}(c)), a separate memory of $\mathbf{I}^5$ disappears, and the old attractor on the left is now associated with a blended memory of $\mathbf{I}^1$, $\mathbf{I}^2$, $\mathbf{I}^4$ (red, dark green and blue lines, respectively) and of $\mathbf{I}^5$ (brown lines). Thus, $\mathbf{I}^5$ was forgotten. 

By the end of learning at $t$$=$$6000$ (Fig.~\ref{fig:fig3.1}(d)), the NN has formed five memories. Four of them are separate memories of $\mathbf{I}^1$ (red lines), $\mathbf{I}^2$ (dark green lines), $\mathbf{I}^3$ (yellow lines) and $\mathbf{I}^6$ (light green lines). One attractor is associated with a blend of $\mathbf{I}^1$ and $\mathbf{I}^5$  (red and brown lines, respectively). Thus, during learning, attractors were being formed gradually to become associated with some of the individual applied patterns (proper memories) and those associated with more than one applied pattern (incorrectly formed memories). 

Even if the difference between the values of weights at different $t$ can be small, such as at $t$$=$$2000$ (Fig.~\ref{fig:fig3.1}(c)) and $t$$=$$6000$ (Fig.~\ref{fig:fig3.1}(d)), it can be sufficient to induce bifurcations of the vector field of the counterpart autonomous DS (\ref{eqn:eq3.1}). Also, attractors formed at an earlier stage of learning can disappear at later stages, such as a pair of attractors existing at  $t$$=$$500$  (see Fig.~\ref{fig:fig3.1}(b), black diamond in top left corner and red cross in bottom right corner), which disappear by $t$$=$$2000$ (see Fig.~\ref{fig:fig3.1}(d)).

To explore the sensitivity of learning to input data, the NN (\ref{eqn:eq2.1})--(\ref{eqn:eq2.4}) was subjected to the same learning procedure, but with a different training signal $\mathbf{I}(t)$ -- the one formed from training Set~2 (see Tabs.~\ref{supp-tab:tab1-1}--\ref{supp-tab:tab1-2}). The stages of memory formation are illustrated in Fig.~\ref{fig:fig3.1}(e)--(h). Note, that the statistical properties of Sets~1 and ~2  are the same, but the actual values of their vectors $\mathbf{I}^k$ are different. 

Despite the statistical similarity of both sets, they result in very different processes of learning, as is evident from comparing the two boxes in Fig.~\ref{fig:fig3.1}. Attractors appear at very different locations and become associated with true or spurious memories in a seemingly random  manner. These observations are in line with a recent attempt to explain learning by the Hopfield NN using a plastic self-organising vector field, and with the conclusion about its inherent non-explainability \cite{Janson_explainable_plastic_DS_SSRN25}. Unlike with Set~1, with Set~2 the NN develops six proper memories, each associated with an individual vector $\mathbf{I}^k$. 

Our simulations suggest that, although training may lead to creation of attractors whose number equals, or even exceeds, the number of patterns the NN needs to memorise, not all patterns could be represented by a unique attractor. Thus, the development of the full memory critically depends  on the content of the training set. 

\section{Memory formation: bifurcation mechanisms}
\label{sec:bif}

 After observing  stages of memory formation in the NN in the course of its learning, we need to reveal the underlying causes, or \emph{mechanisms}. 
In Sec.~\ref{sec:overview} we quote some results from the area of ML 
with recurrent  ANNs, where those may differ substantially from the continuous-time Hopfield NN considered here, and explain why these results are relevant to our current study.

In Sec.~\ref{sec:challenge} we outline the challenges of bifurcation analysis in \emph{learning} NNs and explain the difficulty of defining such bifurcations within the framework of the standard DS theory. In Sec.~\ref{sec:rig} we propose a rigorous formulation of the problem of finding a bifurcation in a learning NN in terms of standard DS theory. In Sections~\ref{sec:early}, \ref{sec:later} and \ref{sec:learnt} we report the results of bifurcation analysis at various stages of learning. 

\begin{figure*}
\includegraphics[width=1.0\textwidth]{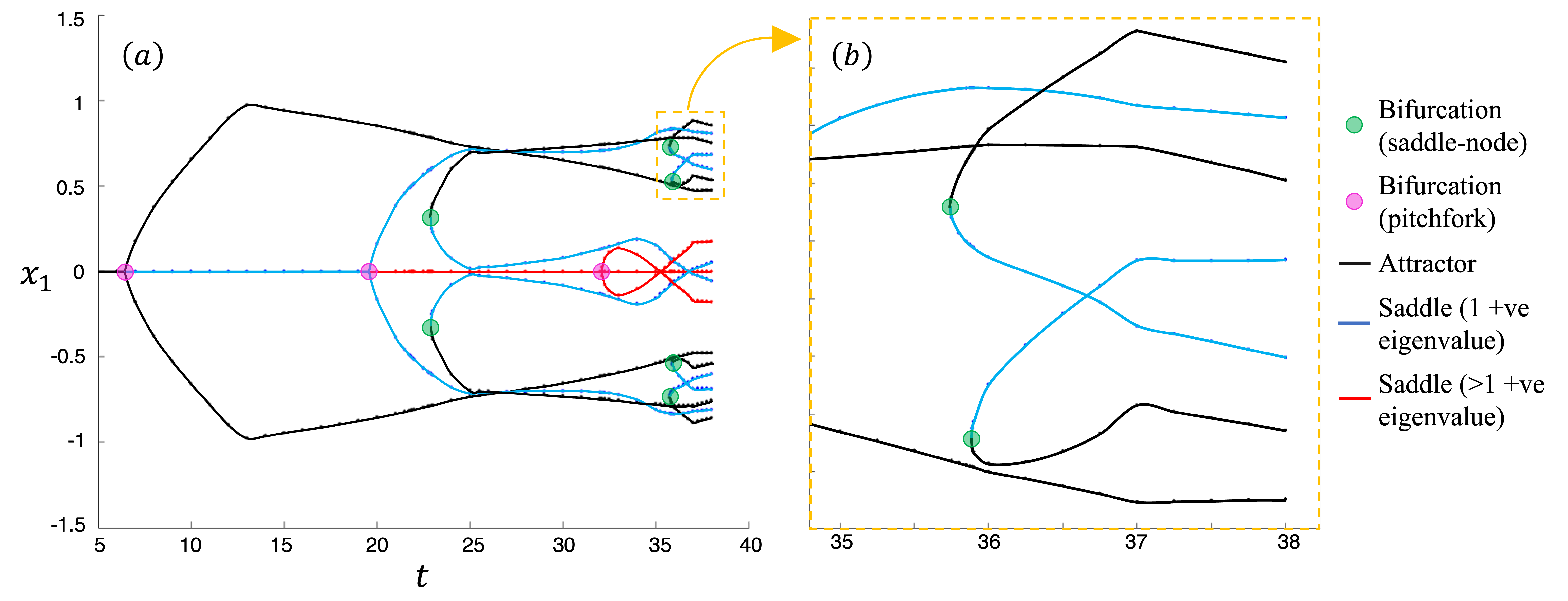}
\caption{One-parameter bifurcation diagram illustrating memory formation in \emph{early} stages of learning by the NN (\ref{eqn:eq2.1})--(\ref{eqn:eq2.4}) from training Set~1 (Tabs.~\ref{supp-tab:tab1-1}--\ref{supp-tab:tab1-2}) with $N$$=$$81$. Panels show $x_1$-coordinates of fixed points of (\ref{eqn:eq3.1}) (solid lines) as functions of control parameter $t$, which coincides with time $t$ in (\ref{eqn:eq2.1})--(\ref{eqn:eq2.4}). Stability of fixed points is indicated by the line colour: stable point, i.e. attractor  potentially associable with memory (black), ``useful''  saddle point with a single positive eigenvalue (blue), and other saddle points (red). Dots along each branch indicate values of $t$ at which the respective fixed point was numerically found and analysed; the dots are connected by interpolating cubic splines. Translucent circles mark bifurcation points: pitchfork at $t$$=$$6.5$ (pink), saddle-node at $t$$=$$22.8783$ and at $t$$\approx$$35.9$ (green). All saddle-node bifurcations occur in pairs due to symmetry in (\ref{eqn:eq3.1}).  Panel (a) shows the bifurcation diagram for $t \in [5,38]$; (b) is a close-up of the rectangular selection in (a). All bifurcations involving attractors are additionally illustrated by phase portraits in Fig.~\ref{fig:Pitchfork}.
}
\label{fig:81N_Bif_Diag}
\end{figure*}

\begin{figure*}
\includegraphics[width=0.9\textwidth]{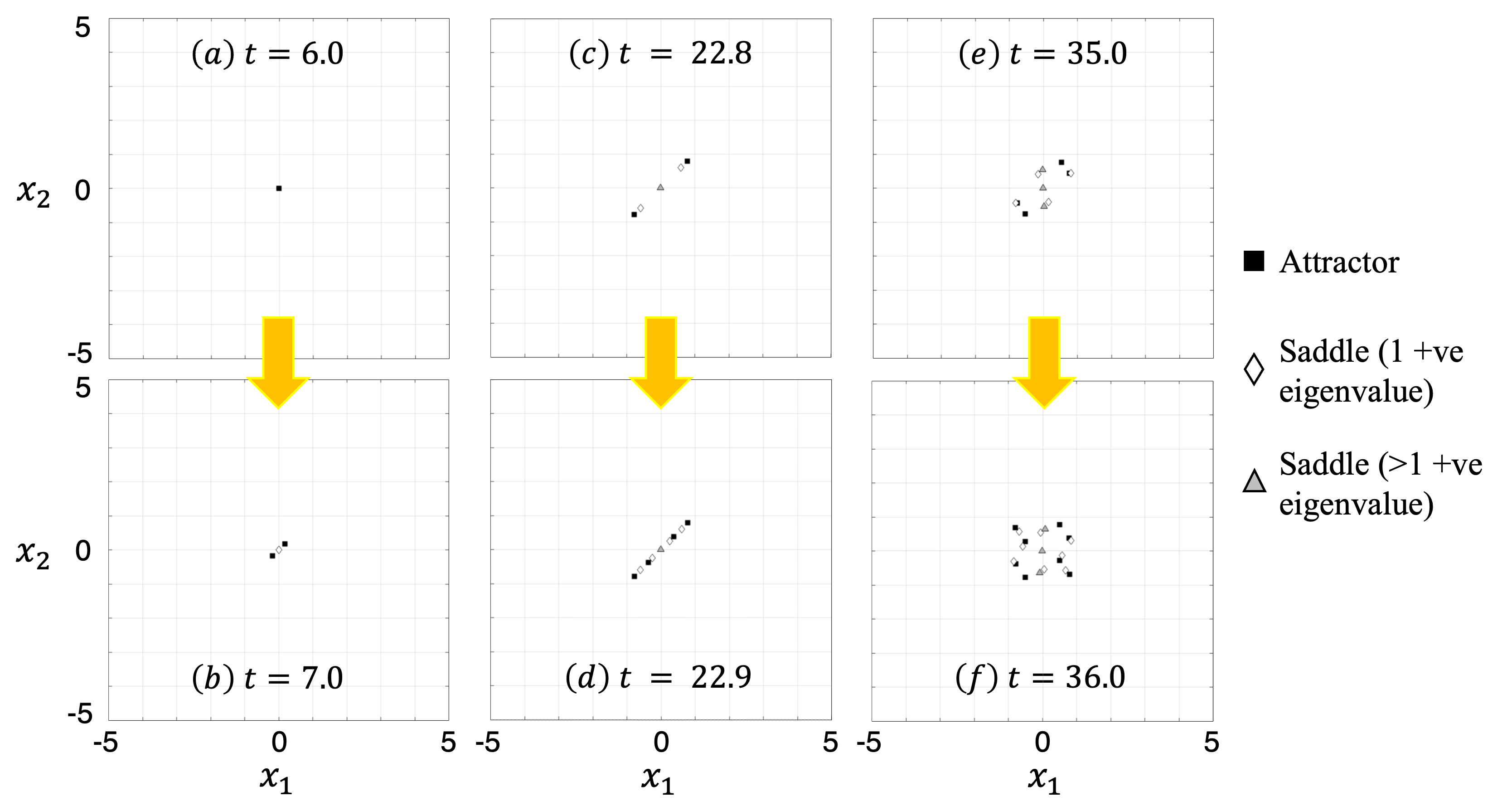}
\caption{Illustration of the first four bifurcations producing new potential memories in the NN (\ref{eqn:eq2.1})--(\ref{eqn:eq2.4}) with $N$$=$$81$ learning from training Set~1, compare with Fig.~\ref{fig:81N_Bif_Diag}. Panels show projections of fixed points on the $(x_1,x_2)$-plane of the NN (\ref{eqn:eq3.1}) before (first row) and after (second row) the respective bifurcations, and provide the value of parameter $t$. Fixed points are marked as: attractors potentially associable with memories  (black boxes), ``useful'' saddle points with one positive eigenvalue (white diamonds), and all other saddle or unstable points (grey triangles). Yellow arrows indicate the flow of time $t$ in (\ref{eqn:eq2.1})--(\ref{eqn:eq2.4}).  Bifurcations illustrated are: (a)--(b) pitchfork bifurcation at $t$$=$$6.5$, which destabilises the point $\mathbf{0}$ and produces two new attractors; (c)--(d) a pair of saddle-node bifurcations at $t$$=$$22.8783$ creating two new attractors (compare with Fig.~\ref{fig:81N_Bif_Diag}(a)), and (e)--(f) two pairs of saddle-node bifurcations at $t$$\approx$$35.9$ creating four new attractors (compare with Fig.~\ref{fig:81N_Bif_Diag}(b)). 
}
\label{fig:Pitchfork}
\end{figure*}

\subsection{Machine learning results relevant to bifurcations}
\label{sec:overview}

ANNs predominantly used in modern ML can  have discrete time and/or a different evolution rule for the neuron states as compared to the Hopfield NN (\ref{eqn:eq2.1}) or (\ref{eqn:eq3.1}). Moreover, modern ML mostly uses  non-Hebbian learning algorithms, since they are more efficient. 
 
However, despite these differences, in all recurrent NNs the \emph{conceptual} principles behind learning  are fundamentally the same. Namely, during  learning the weights need to be adjusted in such a way, that at the end a certain structure of the vector field in the state space of the NN is formed, which ensures an appropriate number and configuration of attractors and their basins -- often metaphorically called ``attractor landscape'' \cite{Saxena_NN_deep_fighting_catastrophic_forgetting_PNAS22}. 
Therefore, better understanding of Hopfield NNs with Hebbian learning will aid understanding of learning in other recurrent NNs. 

In 1988,  a theoretical idea was put forward that if a NN initially  had only one attractor, then to associate different training vectors with distinct attractors, during learning this NN needs to go through bifurcations (or ``catastrophes'') \cite{Pineda_Hopfield_NN_bifurcation_JC88}. Recently, numerical evidence was obtained that in non-small (up to $200$-neuron) discrete-time NNs, bifurcations do occur during learning
by demonstrating that
the numbers or the types of attractors could change between the consecutive  training epochs \cite{Ribeiro_NN_bifurcations_exploding_gradients_conf20,Haputhanthri_Why_NN_learn_bifurcations_24}. In these studies, the  
 types of bifurcations or of attractors involved were not clearly specified. For more details, see Sec.~\ref{sec:challenge}.  
 
Despite the seemingly critical importance of bifurcations for learning, there have been no  \emph{standard} bifurcation analysis of learning NNs, where bifurcations would be expressed in terms of the stage of learning. There were also no rigorous and direct  demonstrations of  how bifurcations
lead to memory formation.  The absence of such studies 
stems from both technical and  conceptual challenges involved in analysing bifurcations in \emph{learning} NNs, as explained below.

\subsection{Challenges of bifurcation analysis in learning neural networks}
\label{sec:challenge}

We remind the reader that for a general autonomous or a periodically perturbed DS depending on several parameters, bifurcations of codimension one can be detected by fixing all parameters 
except one, and allowing this single parameter to monotonously and  \emph{continuously} increase or decrease, while monitoring a certain limit set and its stability.  A local bifurcation \footnote{Global bifurcations are irrelevant to this study.}  is the instant when this stability  changes, thus leading to a dramatic change in the system behaviour \cite{Kuznetsov_bif_theory_98}. A more sophisticated well-established method of bifurcation analysis is continuation \cite{Allgower_continuation_book90,Doedel_AUTO_97,Ermentrout_XPPAUT_guide_book02}, which can detect bifurcations of  codimensions one and higher, when theoretically any number of parameters can be changed simultaneously. This way, one can in principle build bifurcation manifolds in the parameter space, although this is a challenging task, and in practice these methods are mostly used to obtain bifurcation manifolds in the space of just two parameters. Generally,
these approaches work well when the number of control parameters in a DS is relatively small 
 and 
comparable with the number of state variables. 
Such conventional bifurcation analysis was successfully done for relatively small \emph{non-learning}  NNs with two, three, or four neurons \cite{Das_Hopfield_NN_4-neuron_bifurcation_PhD95,Beer_small_NNs_bifurcations_also_some_large_NNs_AB95,Haschke_NN_discrete_2_3_neurons_non-learning_bifurcations_NC05,Huang_Hopfield_NN_3_neuron_bif_chaos_AMC08,Cervantes-Ojeda_NNs_N_2_bifurcation_diagrams_NC17,Njitacke_Hopfield_4_neurons_bifurcations_IJDC19,Hu_Hopfield_NN_3-neuron_forced_multi-scroll_attractors_enctyption_MTA24}. 
 
{\bf Technical challenge.} In NNs of a non-small size $N$, the number of control parameters $\omega_{ij}$ is \emph{much} larger than the number $N$ of state variables. For example, in the Hopfield NN (\ref{eqn:eq3.1}) with $\omega_{ij}$$=$$\omega_{ji}$, $\omega_{ii}$$=$$0$ and $N$$=$$81$, the number $M$ of distinct weights  calculated as $M$$=$$\frac{N(N-1)}{2}$ is $3240$, which is much 
greater than $81$. For such a NN,  it is not feasible to obtain a bifurcation diagram  containing the bifurcation manifolds  in the space of all weights. As an alternative, one might consider singling out a small number of weights and obtaining a bifurcation diagram in their space, while keeping all other weights fixed. However, if one wishes to understand bifurcations occurring in a learning NN, this is not particularly useful, since for the purpose of learning all weights are important and none are more significant than others. 
Thus, the number of control parameters in NNs represents a technical challenge for their bifurcation analysis.

{\bf Conceptual challenge.} In addition, there is a conceptual challenge. Namely, formally, a learning NN (\ref{eqn:eq2.1})--(\ref{eqn:eq2.4}) is a non-autonomous DS with periodic perturbation $\mathbf{I}(t)$, whose state is $(x_1,\ldots,x_N,\omega_{11},\ldots,\omega_{NN})$ and which has \emph{fixed} parameters $A$, $B_{ij}$  and $g$. With fixed parameters, the DS has no bifurcations! 

However, in ML only $(x_1,\ldots,x_N)$ is regarded as the state of the NN, whereas $\omega_{ij}$ are regarded as parameters. With this, during learning \emph{all} parameters $\omega_{ij}$ change simultaneously. To detect a codimension-one bifurcation occurring in the course of learning, one needs a \emph{single} continuously changing bifurcation parameter unambiguously characterising the stage of learning, which has been missing from the available studies of learning NNs. This difficulty could explain the absence of conventional one-parameter bifurcation diagrams of learning NNs in the form of standard attractors (associated with memories)  displayed as functions of the stage of learning.  

Another problem is that, in a learning NN, the parameters $\omega_{ij}$ change together with the state variables $x_i$, which prevents one from detecting and analysing the limit sets and their stability at fixed parameter values, as needed for standard bifurcation analysis.

Paper Ref.~\cite{Haputhanthri_Why_NN_learn_bifurcations_24} attempts to demonstrate bifurcations in a learning ANN, while emphasising their constructive role. However, it employs a non-rigorous methodology departing from that of the standard DS theory. Namely, the term ``bifurcation'' is redefined as a 
``qualitative change in the arrangement of the slow-point (attractor)
landscape'', and the chosen method of its detection is the observation of  
``a sudden change in the loss function'', rather than the change in the stability of a limit set.
In addition, the paper considers the ANN performance as a function of the number of the training epoch, which is a principally discrete quantity and is not an equivalent of a continuously changing bifurcation parameter required to detect standard bifurcations. Also, it does not reveal  standard attractors or other limit sets as required for standard bifurcation analysis, and instead presents various characteristics of the network performance averaged over the training epochs to comply with its own definition of a bifurcation, rather than with a standard one.  Where attractors and bifurcations are claimed to be detected,  they are not identified or classified from the viewpoint of standard bifurcation theory and are not described clearly. For these reasons, the results of Ref.~\cite{Haputhanthri_Why_NN_learn_bifurcations_24} appear incomplete and require further clarification. While presenting a useful illustration of phenomena resembling bifurcations in learning ANNs, they highlight the lack of, and the burning need for, a rigorous approach to bifurcation analysis in such systems.

\subsection{Bifurcations in learning networks: rigorous formulation}
\label{sec:rig}

To reveal bifurcation mechanisms of memory formation in the learning NN (\ref{eqn:eq2.1})--(\ref{eqn:eq2.4}), we start from formulating the problem in a mathematically rigorous way. We note, that during learning the state of the subsystem (\ref{eqn:eq2.4}) follows a one-dimensional continuous path in the $M$-dimensional space of all $\omega_{ij}$, which we refer to as the ``weight trajectory'' and visualise in Sec.~\ref{sec:memfog}. Such a path could be parametrised by a single parameter, such as arc-length \cite{Sharpe_Arc_Length_Parameterisation}. 

However, in our case for the given training signal $\mathbf{I}(t)$ and with the given initial conditions, the state of (\ref{eqn:eq2.4}) at any time $t$ is fully determined by $t$, and therefore the path it follows can be parametrised by $t$. For the DS (\ref{eqn:eq3.1})  we can assume that all $\omega_{ij}$ are functions of the same parameter $t$, i.e. $\omega_{ij}$$=$$\omega_{ij}(t)$. 
Therefore,  $t$ could be treated as a  single continuously changing \emph{bifurcation parameter}  signifying the stage of learning. This way, the problem of finding bifurcations in a non-autonomous periodic DS (\ref{eqn:eq2.1})--(\ref{eqn:eq2.4}) with fixed parameters, which is ill-posed due to the absence of bifurcations in it, is reformulated as finding bifurcations in an autonomous DS (\ref{eqn:eq3.1}) with time $t'$ and a single bifurcation parameter $t$. Thus, by artificially separating two times, $t$ from (\ref{eqn:eq2.1})--(\ref{eqn:eq2.4}) and $t'$ from (\ref{eqn:eq3.1}), we reduce the problem of detecting bifurcations in a learning NN (\ref{eqn:eq2.1})--(\ref{eqn:eq2.4}) to a standard problem of bifurcation analysis of (\ref{eqn:eq3.1}) with a single continuous bifurcation parameter $t$ representing the stage of learning. The given formulation  allows one to stay within the framework of the standard DS theory and does not require introduction of new attractors or bifurcations, contrary to the proposal of Ref.~\cite{Haputhanthri_Why_NN_learn_bifurcations_24}.

Importantly, bifurcations detected with the increase of a single parameter $t$ will most likely have codimension one. The probability to come across a bifurcation with higher codimension is negligibly small. 

We search for bifurcations while being aware that  in (\ref{eqn:eq3.1}), at any combination of $\omega_{ij}$s, stable fixed points can be the only attractors, and any fixed point (stable or unstable) can have only real eigenvalues and be of a node type.  These facts imply that oscillations in (\ref{eqn:eq3.1}) are impossible, which is supported by the proof involving the concept of an energy function \cite{Hopfield_neurons_graded_response_PNAS84}. 

\subsection{Early stages of learning}
\label{sec:early}

\begin{figure*}
\includegraphics[width=1.0\textwidth]{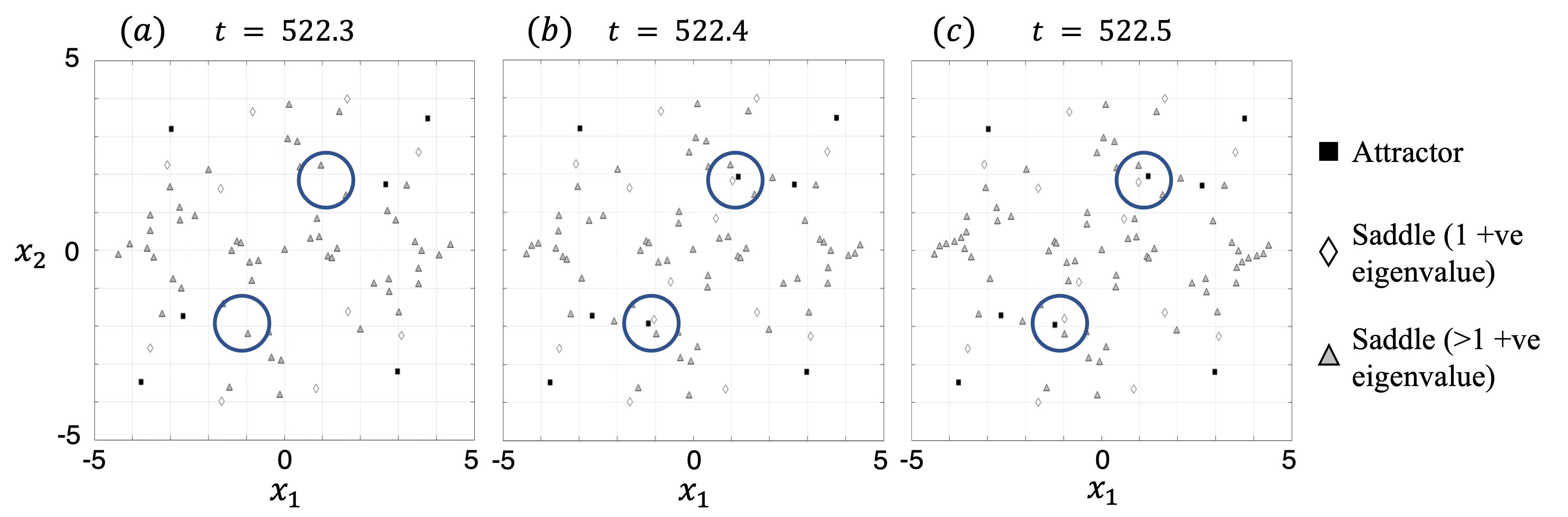}
\caption{Illustration of memory-creating bifurcations at later stages of learning ($t$$>$$38$), specifically, of a pair of saddle-node bifurcations taking place in the NN (\ref{eqn:eq2.1})--(\ref{eqn:eq2.4}) with $N$$=$$81$ that learns from training Set~1 (compare with Fig.~\ref{fig:fig3.1}(a)--(d)). Panels show projections onto the $(x_1,x_2)$-plane of fixed points of the NN (\ref{eqn:eq3.1})  (a) before the bifurcation at $t$$=$$522.3$,  (b) right after the bifurcation at $t$$=$$522.4$,  and  (c) further away from  the bifurcation at $t$$=$$522.5$, with $t$ being the control parameter of (\ref{eqn:eq3.1}) equal to time $t$ of  (\ref{eqn:eq2.1})--(\ref{eqn:eq2.4}). Blue circles highlight areas of the phase space where bifurcations take place.  Fixed points are marked as: attractors potentially associable with memories  (black boxes), ``useful'' saddle points with one positive eigenvalue (white diamonds), and all other saddle or unstable points (grey triangles).
}
\label{fig:fig6.1}
\end{figure*}

Early stages of  learning from training Set~1 are illustrated in Figs.~\ref{fig:81N_Bif_Diag} and \ref{fig:Pitchfork}. Namely, Fig.~\ref{fig:81N_Bif_Diag} presents a one-parameter bifurcation diagram showing components  $x_1$  of all fixed points of (\ref{eqn:eq3.1}) as functions of $t$. Here, $t$ grows from $0$ to $38$, during which the system fully ``processes'' the first three vectors $\mathbf{I}^k$ applied to it ($t $$\in$$ [0,36]$) and starts to process the fourth one ($t $$\in$$ [36,38]$). At the start of learning ($t $$\in$$ [0,6.5]$), the NN (\ref{eqn:eq3.1}) has a single attractor at the origin and no other fixed points (Fig.~\ref{fig:Pitchfork}(a)). 

Figure~\ref{fig:81N_Bif_Diag}  demonstrates the development of new potential memories in the form of stable fixed points (black lines) with their basins. 
Besides attractors, saddle points with a single positive eigenvalue (blue lines) are also critical to memory formation for the following reasons. Firstly, their $(N-1)$-dimensional stable manifolds can serve as boundaries of attractor basins representing memories, as explained in detail in Sec.~\ref{sec:bas}. Secondly, they can participate in saddle-node bifurcations with attractors  causing their deaths, thus potentially destroying memories. We will refer to such saddles as ``useful''. Saddle points with more than one positive eigenvalue (red lines) are not immediately relevant to memory formation. 

Consider a sequence of events induced by the growth of $t$. 
As $t$ reaches  $6.5$ (Fig.~\ref{fig:81N_Bif_Diag}(a)), a pitchfork bifurcation occurs (pink circle) that destabilises the fixed point at $\mathbf{0}$ (blue line starting at $t$$=$$6.5$) and gives rise to two new attractors (black lines extending for $t \in [6.5,38]$). The same effect is illustrated by phase portraits in Fig.~\ref{fig:Pitchfork} as one goes from  (a) a single attractor at $\mathbf{0}$ (black box)  at $t$$=$$6$ to (b) a pair of attractors (black boxes) and a saddle at $\mathbf{0}$ with one positive eigenvalue (white diamond)  at $t$$=$$7$. 

After that,  a cascade of saddle-node bifurcations (green circles in Fig.~\ref{fig:81N_Bif_Diag}) produces more attractors and ``useful'' saddles. Namely, at $t$$=$$22.8783$ the first pair of saddle-node bifurcations give rise to two symmetric pairs of stable and saddle points. This is illustrated in Fig.~\ref{fig:Pitchfork} as one goes from (c)  two pairs of attractors (black boxes) and ``useful'' saddles (white diamonds) 
at $t$$=$$22.8$ to (d) four such pairs at $t$$=$$22.9$ meaning that the NN acquires four potential memories.  

The subsequent two pairs of saddle-node bifurcations occurs almost simultaneously at $t$$\approx$$35.9$, as more clearly seen in Fig.~\ref{fig:81N_Bif_Diag}(b). The transition of the phase portrait from before to  after bifurcations is illustrated in Fig.~\ref{fig:Pitchfork}(e)--(f). As a result, four more pairs of attractors and of ``useful'' saddles appear. Note, that saddle-node bifurcations always occur in pairs at symmetric locations in the phase space due to the symmetry of  $u_i$ in (\ref{eqn:eq3.1}) resulting, in turn, from the symmetry of $F(x)$.

Interestingly, bifurcations giving rise to new attractors occur at the ends of the time intervals of length $t_s$$=$$12$, during which individual vectors $\mathbf{I}^k$ are applied to (\ref{eqn:eq2.1})--(\ref{eqn:eq2.4}). For smaller $t_s$, such bifurcations might not have happened. This is consistent with the idea from psychology that memories need time to ``sink in'' \cite{Smith_exposure-time_memory_formation_visual_PHFS82}. 

A cascade of pitchfork bifurcations  of the origin (pink circles in Fig.~\ref{fig:81N_Bif_Diag}(a)), each endowing  the fixed point at $\mathbf{0}$ with one more positive eigenvalue,  and also producing new saddle fixed points with \emph{more} than one positive eigenvalue, do not affect memory formation and are indicated only for completeness.  

Note, that in (\ref{eqn:eq3.1}) the fixed point at $\mathbf{0}$ is stable only if all $\omega_{ij}$ are close to zero. In that case, the point at $\mathbf{0}$ is the only attractor in the system, which is in agreement with early studies of similar NNs \cite{Doya_Bifurcations_in_the_learning_of_recurrent_neural_networks_IEEE92}. This is true in the beginning of learning illustrated here because of the chosen ICs, i.e. $|\omega_{ij}(0)|$$\leq$$ 0.01$, see Fig.~\ref{fig:omega_t}(a) at small $t$. However, if  in (\ref{eqn:eq2.1})--(\ref{eqn:eq2.4}) the ICs for $\omega_{ij}$  are not small, then already at $t$$=$$0$ the NN (\ref{eqn:eq3.1}) can have an unstable fixed point at $\mathbf{0}$ and multiple stable fixed points away from $\mathbf{0}$. In that case, the pitchfork bifurcation at relatively small $t$ does not occur, and only saddle-node bifurcations take place.

\begin{figure}[H]
\includegraphics[width=0.5\textwidth]{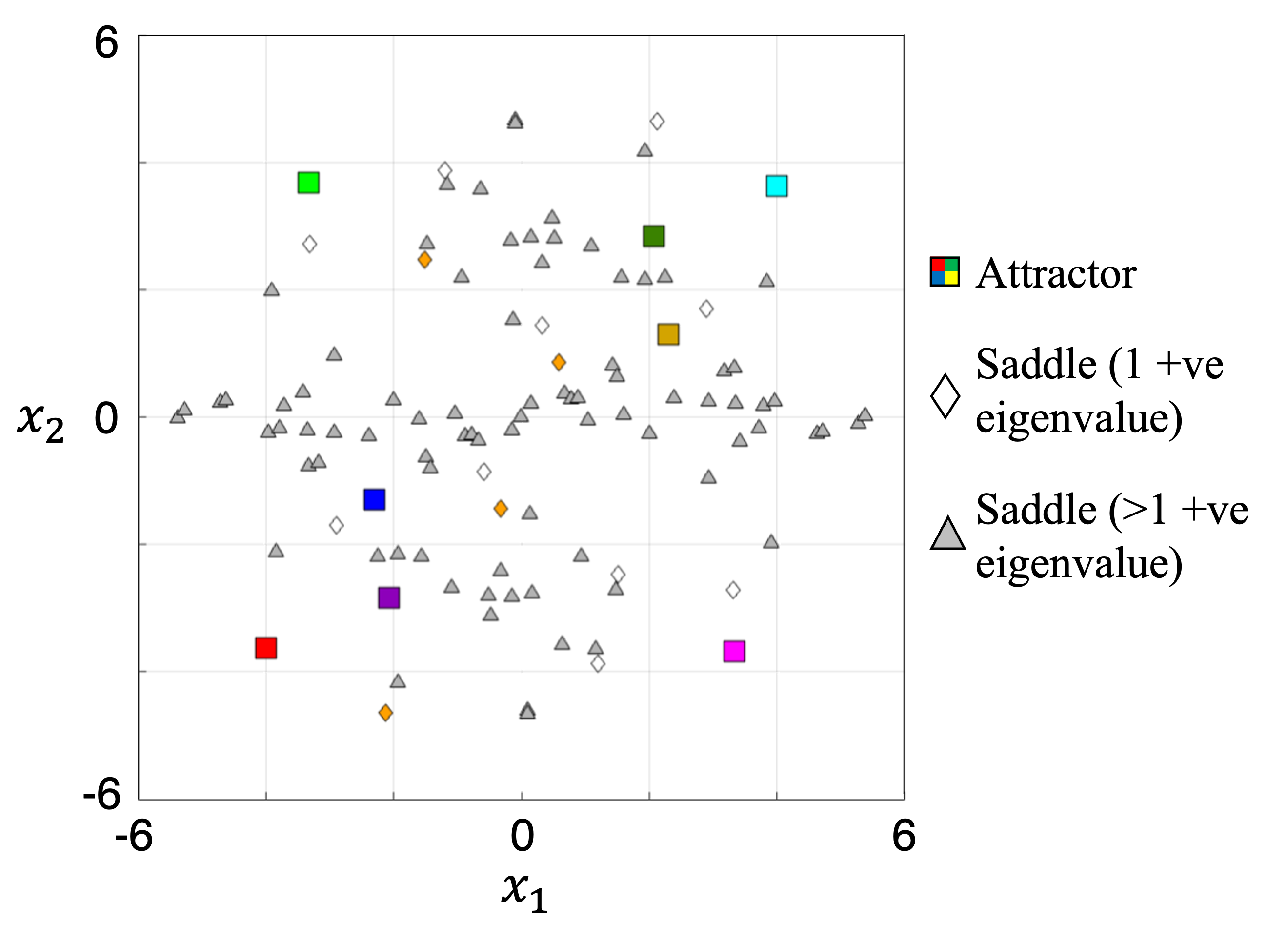}
\caption{All fixed points developed in the $81$-neuron network (\ref{eqn:eq2.1})--(\ref{eqn:eq2.4}), described in caption to Fig.~\ref{fig:fig3.1} and trained on Set~1, at the end of learning at $t$$=$$6000$. Red, dark-green, yellow, blue, and light green boxes represent attractors identical to those shown in Fig.~\ref{fig:fig3.1}(d) by black diamonds  and associated with true memories, their colours matching those of the trajectories in Fig.~\ref{fig:fig3.1}(d) launched from respective $\mathbf{I}^k$s and their vicinities. Purple, pink and cyan boxes show spurious memories identical to red crosses in Fig.~\ref{fig:fig3.1}(d). Diamonds depict ``useful''  saddles with a single positive eigenvalue. Of these, four orange diamonds mark saddles whose manifolds are shown in Fig.~\ref{fig:basins_sad} as basin boundaries. Grey triangles indicate all other saddle or unstable fixed points.} 
 \label{fig:attractor_colours}
\end{figure}

\subsection{Later stages of learning}
\label{sec:later}

At later stages of learning, saddle-node bifurcations remain the prime mechanism of memory formation, as demonstrated in Fig.~\ref{fig:fig6.1}. Namely, as $t$ grows from $522.3$ (Fig.~\ref{fig:fig6.1}(a)) to  $522.4$ (Fig.~\ref{fig:fig6.1}(b)), two pairs of stable and saddle fixed points appear in two symmetrically located areas of the phase space highlighted by blue circles. As $t$ grows to $522.5$, within each newly born pair of fixed points, these points move away from each other  (Fig.~\ref{fig:fig6.1}(c)). The given behaviour of the fixed points, together with the analysis of their eigenvalues, clearly points to saddle-node bifurcations. 

~\vspace{-10mm}

\subsection{End of learning}
\label{sec:learnt}

As a result of the sequence of bifurcations discussed above, the NN  (\ref{eqn:eq2.1})--(\ref{eqn:eq2.4}) with $N$$=$$81$ trained on Set~1, by the time $t$$=$$6000$ develops fixed points shown in Fig.~\ref{fig:attractor_colours} (compare with Fig.~\ref{fig:fig3.1}(d)). Namely, Fig.~\ref{fig:attractor_colours} shows projections onto the $(x_1,x_2)$-plane of \emph{all} fixed points of the NN (\ref{eqn:eq3.1}) considered with the values $\omega_{ij}$ of (\ref{eqn:eq2.4}) at  $t$$=$$6000$. 
 Attractors are depicted with coloured boxes, whose colours are matched with those of their basins shown in Figs.~\ref{fig:basins} and \ref{fig:basins_sad} and discussed in Sec.~\ref{sec:bas}. ``Useful'' saddle fixed points with a single positive eigenvalue are depicted by diamonds. Of these, orange diamonds indicate four saddles whose manifolds forming boundaries of attraction basins are illustrated in Fig.~\ref{fig:basins_sad}. Grey triangles indicate all other saddle or unstable fixed points. 

\section{Forgetting and its relevance to bifurcations}
\label{sec:memfog}

Abrupt forgetting typically and famously occurs  when an ANN, after being trained on a certain task, \emph{switches} to learning the next task. In this situation, all or part of the knowledge of the previously learned task is abruptly lost. This phenomenon is called 
 \emph{catastrophic forgetting} (or catastrophic interference) and remains one of the major deficiencies of modern ANNs \cite{McCloskey_catastrophic_forgetting_PLM89,French_catastrophic_forgetting_in_NNs_TCS99,Kirkpatrick_overcoming_catastrophic_forgetting_in_NNs_PNAS17}. 
 It has been assumed that the knowledge of a task is contained in weights, and catastrophic forgetting has been explained in terms of weights shared by different tasks, and of stability-plasticity dilemma \cite{French_catastrophic_forgetting_in_NNs_TCS99}. 
 
If we interpret this problem in terms of the DS theory, 
the knowledge of a task would be represented by the appropriate configuration of attractors and their basins, 
where each attractor together with its basin is associated with some input vector from the training set and hence represents the individual memory of this vector. 
The task would be deemed learnt after the weights have converged to low-amplitude oscillations centred around some constant values as it happens in Fig.~\ref{fig:omega_t} at $t$ close to $6000$. Information relevant to the new task would consist of a sequence of new training vectors 
repeatedly applied to the network after the weights have converged as described above. Catastrophic forgetting of the previously learnt task would mean the disappearance of any attractors representing any of the vectors from the original training set. 

The existing literature suggests that a full elimination,  rather than the reduction, of catastrophic forgetting is only possible if ANNs expand their sizes to accommodate new knowledge \cite{Rusu_progressive_NNs_against_catastrophic_forgetting_arxiv16,Yoon_Lifelong_Learning_with_Dynamically_Expandable_Netw_arxive18,Fayek_progressive_learning_expanding_eliminate_catastr_forgt_NN20}. With this, the existing research suggests that in ANNs of \emph{fixed} size catastrophic forgetting can only be mitigated, but not eradicated \cite{Kemker_catastrophic_forgetting_deep_NNs_AAAI18}. 
Since the nature of  the underlying mechanisms of catastrophic forgetting have remained largely unknown, 
conventional strategies for mitigating catastrophic forgetting have been predominantly empirical \cite{Robins98,Kirkpatrick_overcoming_catastrophic_forgetting_in_NNs_PNAS17}. 

While in 1988 it was suggested that during learning  a NN could undergo bifurcations \cite{Pineda_Hopfield_NN_bifurcation_JC88}, since 1993 bifurcations have largely been regarded 
as detrimental to learning \cite{Doya_NN_bifurcations_leanring_IEEE93}.  This is because they can cause  gradient ``explosions'' in  loss functions used in popular algorithms governing evolution of weights in the course of learning (which are different from those of the Hebbian learning considered here)\cite{Doya_NN_bifurcations_leanring_IEEE93,Pascanu_NN_learning_bifurcations_conf13}. 

A smaller body of work shows that, as ANNs learn a \emph{single} task, they can also experience memory loss similar to catastrophic forgetting \cite{Toneva_NN_forgetting_single_task_arXiv19}. However, the term ``catastrophic forgetting'' is usually tightly linked to switching between the training tasks, which is regarded as conceptually distinct from single-task training. We argue that this conceptual distinction occurs because the exact mechanisms responsible for the loss of memory in either of the two settings have been unknown. 

In a recent paper Ref.~\cite{Eisenmann_NN_discrete_2-6_neuron_bifurcation_NeurIPS23} it was demonstrated that during the process of training on a \emph{single} task, i.e. with the same set of training data applied repeatedly, at a certain stage of learning the loss function of a piecewise-linear, discrete-time ANN experienced a jump,  which was accompanied by the loss of knowledge about the learnt behaviour. This event was linked to a bifurcation as its cause and referred to as ``catastrophic forgetting'', despite the respective setting being different from the one conventionally associated with the given term.   

In 
mathematical language, a ``bifurcation occurring during learning'' can be interpreted as an event when a bifurcation manifold in the weight space of a non-learning NN is crossed by the weight trajectory of its counterpart learning NN.  
It was proposed that for a good learning, the initial weights of the learning NN should be set to values corresponding to a sufficient number of attractors, and 
the weight trajectory should steer clear of all ``bifurcation boundaries'' \cite{Doya_NN_bifurcations_leanring_IEEE93,Pascanu_NN_learning_bifurcations_conf13}, i.e. of bifurcation manifolds. In practice, this is difficult to achieve because in non-linear high-dimensional multi-parameter  DSs, the locations of bifurcation manifolds are not known in advance. 

Interestingly, Ref.~\cite{Raghavan_NN_learning_as_path_in_weights_space_NatMI24} proposes methods to learn without catastrophic  forgetting by confining the weight trajectories to appropriate parts of the weight space, whose boundaries it does not identify with bifurcation manifolds. Only in Ref.~\cite{Eisenmann_NN_discrete_2-6_neuron_bifurcation_NeurIPS23}  an explicit  link 
 is proposed 
between catastrophic forgetting and bifurcations 
for a special kind of a NN -- a piecewise-linear NN with discrete time  -- which is different from the smooth continuous-time NNs   studied here. Note, that bifurcations in non-smooth DSs may not coincide with those in the smooth ones \cite{Makarenkov_bifurcations_non-smooth_PhD12}. With this, NNs modelled by smooth DSs 
are mathematically closer to both the biological NNs, and the neuromorphic devices. Thus, it is particularly important to reveal the general principles behind all instances of 
forgetting in smooth NNs. 
 
{\bf Forgetting hypothesis.} Here, we propose and test a hypothesis that  in smooth recurrent NNs, forgetting in both settings (i.e. in  single-task training and while training on the next task after learning the previous task) is caused by the weight trajectory  of the \emph{learning} NN crossing some bifurcation manifold in the weight space  of its counterpart \emph{non-learning} NN. We also point out that if two different phenomena have similar manifestations and identical causes, they would appear conceptually linked if not conceptually identical. We therefore argue that, if the above hypothesis is true, the different cases of forgetting in a learning NN would represent essentially the same phenomenon, which could be referred to as "catastrophic forgetting in a broader sense". Here, the term ``catastrophic'' would refer to the abrupt, unpredictable nature of memory loss. 

In Sec.~\ref{sec:forget} we demonstrate the validity of this hypothesis for single-task training of the Hopfield NN with Hebbian learning. In Sec.~\ref{sec_cat} we verify this hypothesis for the same NN, which, after being trained on the first task, starts learning the second task -- i.e. for the conventionally understood catastrophic forgetting.


\section{Forgetting via bifurcation in single-task training}
\label{sec:forget}

Here we reveal the bifurcation mechanisms of forgetting in the Hopfield NN (\ref{eqn:eq2.1})--(\ref{eqn:eq2.4}) learning a single task. 

\subsection{Forgetting in single-task training: observation of attractors}
\label{sec:forget_obs}

\begin{figure*}[ht]
\includegraphics[width=1.0\textwidth]{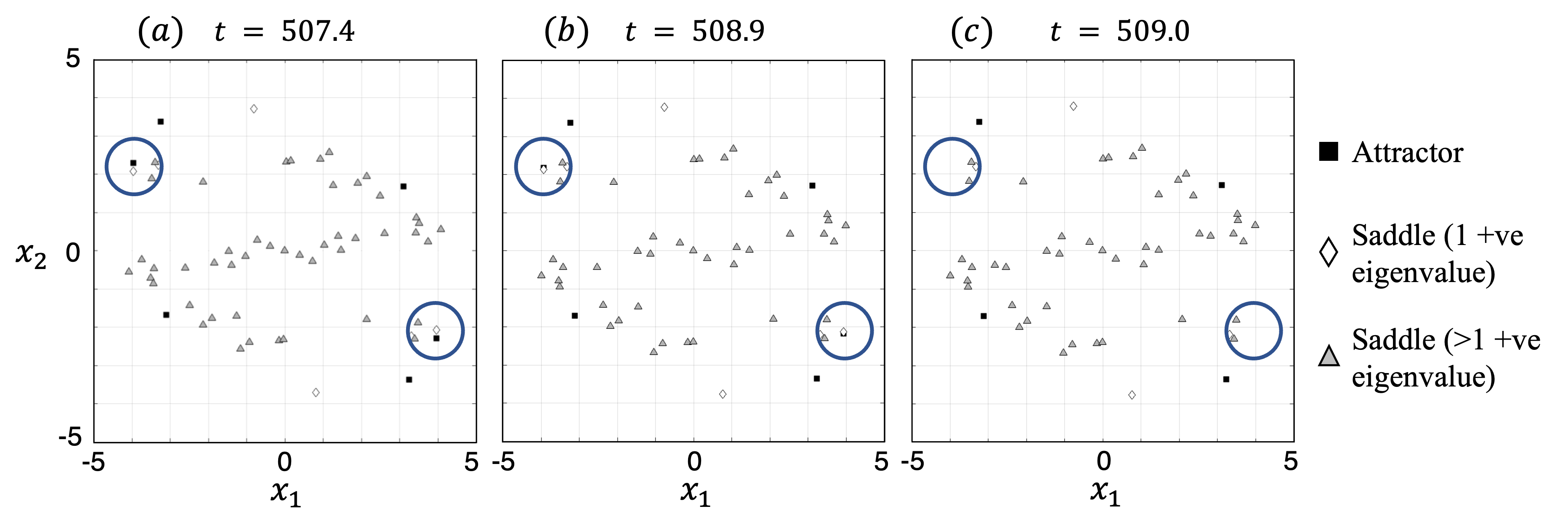}
\caption{Illustration of memory-destroying bifurcations 
during single-task training, specifically, of a pair of saddle-node bifurcations taking place in the NN (\ref{eqn:eq2.1})--(\ref{eqn:eq2.4}) with $N$$=$$81$ learning from training Set~1, compare with Fig.~\ref{fig:fig3.1}(a)--(d). Panels show projections onto the $(x_1,x_2)$-plane of fixed points of the NN (\ref{eqn:eq3.1}) (a) before the bifurcation at $t$$=$$507.4$, (b) immediately before the bifurcation at $t$$=$$508.9$,  and (c) after the bifurcation at $t$$=$$509.0$, with $t$ being the control parameter of (\ref{eqn:eq3.1}) equal to time $t$ of  (\ref{eqn:eq2.1})--(\ref{eqn:eq2.4}). Blue circles highlight areas of the phase space where bifurcations take place.  Fixed points are marked as: attractors potentially associable with memories (black boxes), ``useful'' saddle points with one positive eigenvalue (white diamonds), and all other saddle or unstable points (grey triangles).
}
\label{fig:Forgetting}
\end{figure*}

Firstly, we observe that during the learning phase, as the weights evolve continuously, attractors can not only be born, but can also disappear abruptly. One example of disappearance of a pair of attractors in the NN (\ref{eqn:eq2.1})--(\ref{eqn:eq2.4}) with $N$$=$$81$ learning from training Set~1 is illustrated in the phase portraits of Fig.~\ref{fig:Forgetting}. Namely,  at $t$$=$$507.4$ (Fig.~\ref{fig:Forgetting}(a)), inside the blue circles there are two pairs of stable (black boxes) and saddle (white diamonds) fixed points. At $t$$=$$508.9$ (Fig.~\ref{fig:Forgetting}(b)), these pairs of points come very close to each other, and at $t$$=$$509.0$ (Fig.~\ref{fig:Forgetting}(c)) they no longer exist.  Before the fixed points disappear, their largest eigenvalues approach zero: for the stable fixed points from the negative values, and for the saddle points from the positive values. The phase portraits, together with eigenvalues, suggest that a pair of saddle-node bifurcations took place. 

If any of the attractors discussed above coded a valid category, their sudden death would constitute 
forgetting.

\subsection{Forgetting via bifurcation in a small network: full evidence}
\label{sec:small}

It is well appreciated that a small NN cannot learn efficiently because of the highly limited number of  attractors it could develop. However, before analysing the $81$-neuron network with $3240$ control parameters (weights), in this Subsection we consider a three-neuron network with only three distinct weights, because it permits bifurcation analysis in the space of \emph{all} of its control parameters. This analysis will prepare the ground for the bifurcation analysis of the NN of size $81$, which is done in Sec.~\ref{sec:large}. 

The idea of our analysis is as follows. \emph{First}, we will obtain a three-parameter bifurcation diagram of an autonomous system (\ref{eqn:eq3.1}) with $N$$=$$3$ in the space $(\omega_{12},\omega_{13},\omega_{23})$, which will contain some two-dimensional bifurcation manifolds. Here, we will employ the assumptions of a standard bifurcation analysis, that the weights $\omega_{ij}$ are \emph{not} parametrised by $t$ and can take any values independently of each other. 

 \emph{Second}, we will follow learning by the NN (\ref{eqn:eq2.1})--(\ref{eqn:eq2.4}) with $N$$=$$3$, and will record the weight trajectory  $\pmb{\omega}(t)$$=$$\left( \omega_{12}(t),\omega_{13}(t),\omega_{23}(t) \right)$ of subsystem (\ref{eqn:eq2.4}).  \emph{Thirdly} and finally, in the parameter space $(\omega_{12},\omega_{13},\omega_{23})$ we will superimpose the bifurcation manifolds with the weight trajectory. 
 
For the NN with $N$$=$$3$, the input signal $\mathbf{I}(t)$ was constructed from the training set of three-dimensional vectors $\mathbf{I}^k$ ($k$$=$$1,2,3$) whose coordinates are $\pm 1$ chosen at random with equal probabilities (see Tab.~\ref{supp-tab:tab2} in Supplementary Note).  The  procedures to form $\mathbf{I}(t)$ from the set of $\mathbf{I}^k$, and to generate initial conditions for $x_i$ and  $\omega_{ij}$,  are equivalent to those described in Sec.~\ref{sec:model} for $N$$=$$81$ with suitable adjustments accounting for the different value of $N$ and the number of vectors $\mathbf{I}^k$. 
The values of $A$ and $B$ are also the same as for $N$$=$$81$. We only change the value of $g$ to $g$$=$$5$, because with $g$$=$$0.3$ the small NN produces no bifurcations while processing $\mathbf{I}(t)$. 

\begin{figure*}[ht]
\includegraphics[width=0.9\textwidth]{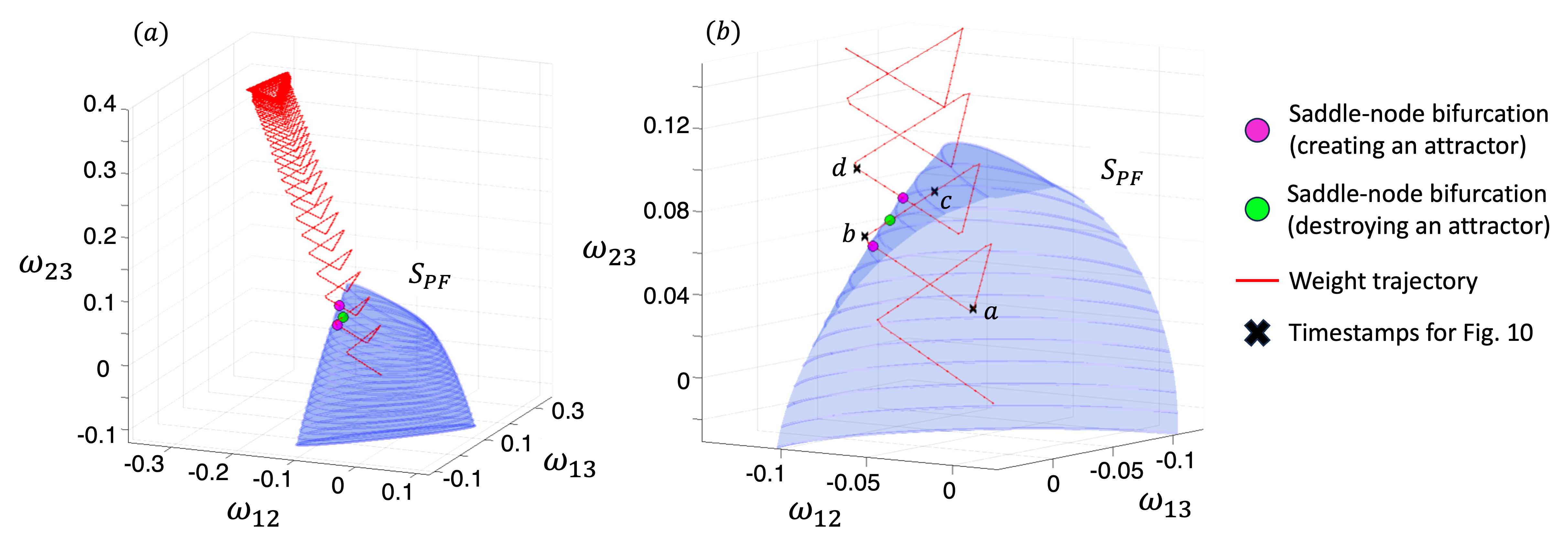}
\caption{Illustration of bifurcation mechanisms of memory formation and forgetting in a \emph{small} NN  (\ref{eqn:eq2.1})--(\ref{eqn:eq2.4}) with $N$$=$$3$, $A$$=$$30$, $B$$=$$300$, $g$$=$$5$ and $\lambda$$=$$1.4$,  trying to memorise vectors $\mathbf{I}^k$ ($k$$=$$1,2,3$, see Tab.~\ref{supp-tab:tab2} of Supplementary Note). Panels show a superposition in the space of weights $\left( \omega_{12},\omega_{13},\omega_{23} \right)$ of a three-parameter bifurcation diagram  of the NN (\ref{eqn:eq3.1}) with $N$$=$$3$ and $g$$=$$5$  with the weight trajectory 
$\pmb{\omega}(t)$$=$$\left( \omega_{12}(t),\omega_{13}(t),\omega_{23}(t) \right)$ (red line) of subsystem (\ref{eqn:eq2.4}) as the NN (\ref{eqn:eq2.1})--(\ref{eqn:eq2.4}) learns. Panel (b) is a zoomed version of (a). \\
$S_{PF}$ (blue surface) is the manifold of pitchfork bifurcation of the fixed point at $\mathbf{0}$. Within the volume bounded by $S_{PF}$, the point $\mathbf{0}$ is stable; outside, the point $\mathbf{0}$ is saddle and there are additionally two stable fixed points.  The trajectory $\pmb{\omega}(t)$ (red line) starts near the origin inside the volume bounded by $S_{PF}$.  Each zag of $\pmb{\omega}(t)$  lasts approximately $t_s$$=$$12$ time units, but see text for a more accurate description. Within $48$ time units, $\pmb{\omega}(t)$ crosses $S_{PF}$ three times by moving out (pink circle near point $b$ in (b)), in (green circle in (a) and (b)) and out again (pink circle near point $d$ in (b)). Panel (b) shows time stamps (black crosses) illustrated by phase portraits in Fig.~\ref{fig:pp_small}.
}
\label{fig:fig9.1}
\end{figure*}

The only bifurcation occurring in (\ref{eqn:eq3.1}) with $N$$=$$3$ is a pitchfork bifurcation of the fixed point at $\mathbf{0}$, which destabilises $\mathbf{0}$ and creates a symmetric pair of attractors. 
To obtain a three-parameter bifurcation diagram, 
we fix $\omega_{23}$ and perform a standard two-parameter bifurcation analysis on the plane $(\omega_{12},\omega_{13})$ using software XPPAUT \cite{xppaut07}. After revealing the curve of pitchfork bifurcation  for the given $\omega_{23}$, we change the value of $\omega_{23}$ and repeat the two-parameter bifurcation analysis. 
We do so for a set of values of $\omega_{23}$ sampled from the interval $[-0.1,0.1]$ with the step $0.01$. As a result,  we obtain a collection of one-dimensional curves of pitchfork bifurcation. We then place these curves in $(\omega_{12},\omega_{13},\omega_{23})$-space and use as a frame to render the manifold (the surface) of pitchfork bifurcation visualised in Fig.~\ref{fig:fig9.1} as the blue surface marked $S_{PF}$. 

Figure~\ref{fig:fig9.1}(a)--(b) shows how $\pmb{\omega}(t)$ (red line) behaves relative to the bifurcation manifold $S_{PF}$ (blue surface) during learning. At $t$$=$$0$ the trajectory starts  inside the region bounded by 
$S_{PF}$, for which there is a single fixed point in the NN -- a stable point at $\mathbf{0}$ 
signifying the absence of memories.  Each ``zag'' (an approximately straight line segment)  of the trajectory lasts  approximately $t_s$$=$$12$ time units, which is equal to the time of application of a single vector $\mathbf{I}^k$ ($k$$=$$1,2,3$)  from the training set given in  Tab.~\ref{supp-tab:tab2} of Supplementary Note. The exception is the very first zag, which lasts approximately $13$ time units. This is due to a delay of about $1$ time unit between the end of the previous $\mathbf{I}^k$ and the abrupt turn of the weight trajectory in response to application of the next $\mathbf{I}^k$. Each $m^{\textrm{th}}$ zag starts at approximately $\left( 1+t_s (m-1) \right)$ and ends at $\left(1+t_s m \right)$ time units. Note, that the period of $\mathbf{I}(t)$ is $3t_s$$=$$36$.

In Fig.~\ref{fig:fig9.1}, small dots on $\pmb{\omega}(t)$ mark the values of $t$ sampled each time unit, and 
panel (b) is a zoomed version of (a) showing the details of $\pmb{\omega}(t)$ crossing $S_{PF}$. 

\begin{figure}
\includegraphics[width=0.5\textwidth]{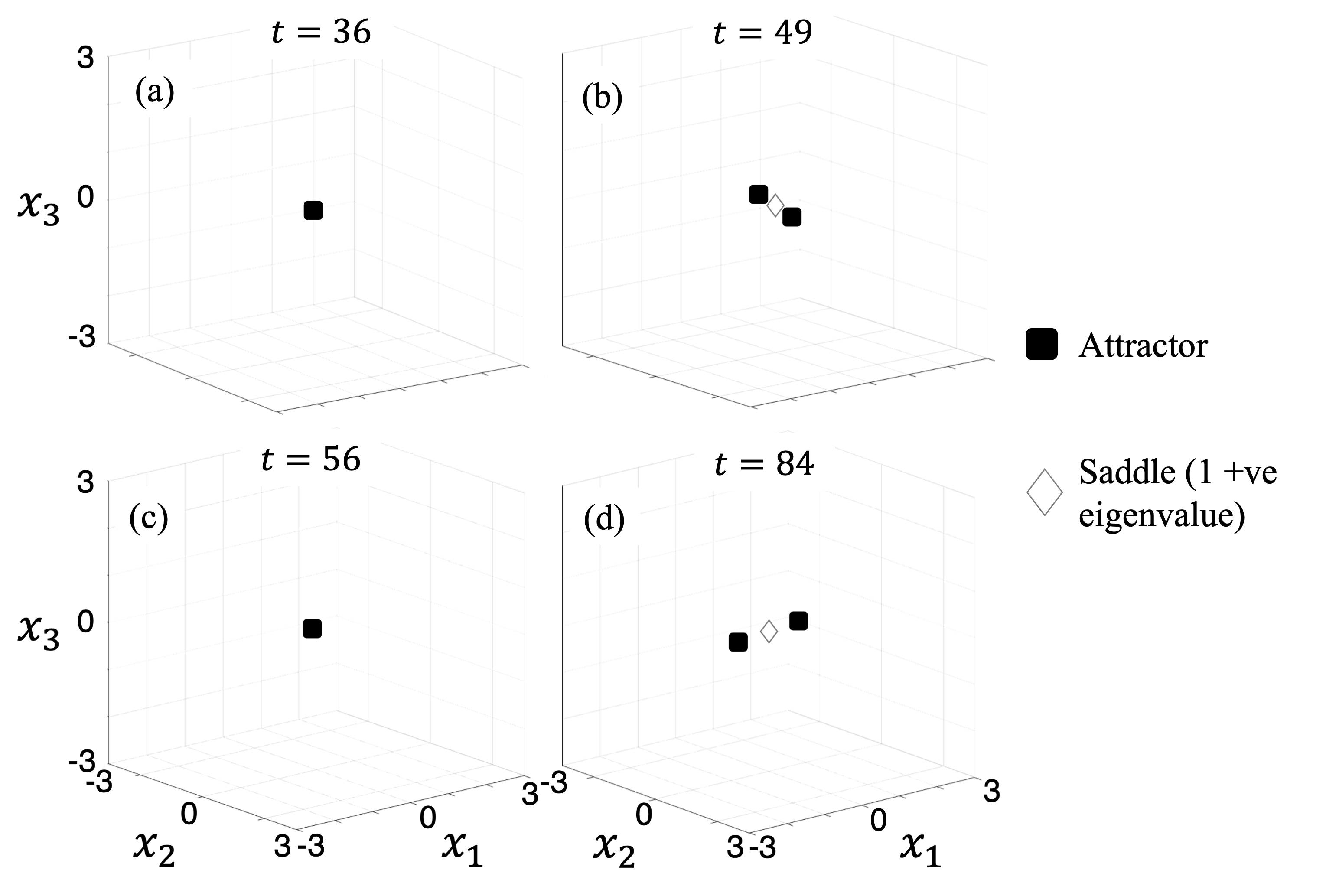}
\caption{Phase portraits illustrating stages of memory formation and forgetting in the small learning NN  (\ref{eqn:eq2.1})--(\ref{eqn:eq2.4}) described in caption to Fig.~\ref{fig:fig9.1}. Panels (a)--(d) show fixed points of (\ref{eqn:eq3.1}) at various values of parameter $t$ corresponding to various stages of learning by the NN  (\ref{eqn:eq2.1})--(\ref{eqn:eq2.4}), marked by black crosses on the weight trajectory in Fig.~\ref{fig:fig9.1}(b) next to matching letters $a$--$d$. Fixed points are: attractors (black boxes) and ``useful'' saddles with a single positive eigenvalue (white diamonds). Between  (a) $t$$=$$36$  and  (b)  $t$$=$$49$ the NN increases the number of potential memories from one to two (black boxes). By (c)  $t$$=$$56$ their number reverts to one, which could constitute forgetting. However, by (d) $t$$=$$84$  the NN has two attractors again (black boxes), which are placed at different locations compared to (b) because their birth occurs at a different combination of weights, as evidenced by different locations of two pink circles in Fig.~\ref{fig:fig9.1}(b).
}
\label{fig:pp_small}
\end{figure}

By $t$$=$$36$, after making three zags, $\pmb{\omega}(t)$ reaches point $a$ (Fig.~\ref{fig:fig9.1}(b)). It is still inside the region bounded by $S_{PF}$, where the stable origin $\mathbf{0}$ is the only fixed point, as illustrated by the phase portrait in Fig.~\ref{fig:pp_small}(a). During the next zag, 
at $t$$\approx$$48$,   $\pmb{\omega}(t)$ crosses $S_{PF}$ (pink circle near point $b$ in Fig.~\ref{fig:fig9.1}(b)), thus destabilising $\mathbf{0}$ and creating a pair of stable fixed points nearby. The phase portrait shortly after the pitchfork bifurcation, corresponding to $t$$=$$49$ and point $b$ in Fig.~\ref{fig:fig9.1}(b), is given in Fig.~\ref{fig:pp_small}(b) and shows two attractors (black boxes) potentially associable with memories. 

During the next zag, at $t$$\approx$$52$, $\pmb{\omega}(t)$ crosses $S_{PF}$ again (green circe in Fig.~\ref{fig:fig9.1}(a)--(b)) and returns to the region bounded by $S_{PF}$, where it travels for a while. A phase portrait typical of this travel, corresponding to $t$$=$$56$ and marked by point $c$ in Fig.~\ref{fig:fig9.1}(b),  is given in Fig.~\ref{fig:pp_small}(c). The fixed point $\mathbf{0}$ (signifying no memories) becomes stable again, and all other attractors (potential memories) disappear, which signifies 
forgetting. 

During the seventh zag that occurs for $t$$\in$$[73,85]$, at $t$$\approx$$79$,  $\pmb{\omega}(t)$ crosses $S_{PF}$ (pink circle near point $c$ in Fig.~\ref{fig:fig9.1}(b)) and brings back two symmetric attractors. The phase portrait after the bifurcation, at $t$$=$$84$ (point $d$ in Fig.~\ref{fig:fig9.1}(b)), is shown in Fig.~\ref{fig:pp_small}(d).  The attractor locations here are different from those at point $b$ because the parameters $\omega_{ij}$ have different respective values.  After that, the trajectory $\pmb{\omega}(t)$ moves away from $S_{PF}$ in an oscillatory manner (Fig.~\ref{fig:fig9.1}(a)), never crosses $S_{PF}$ again, and converges to a closed set resembling a limit cycle. 

\begin{figure*}
\includegraphics[width=1.0\textwidth]{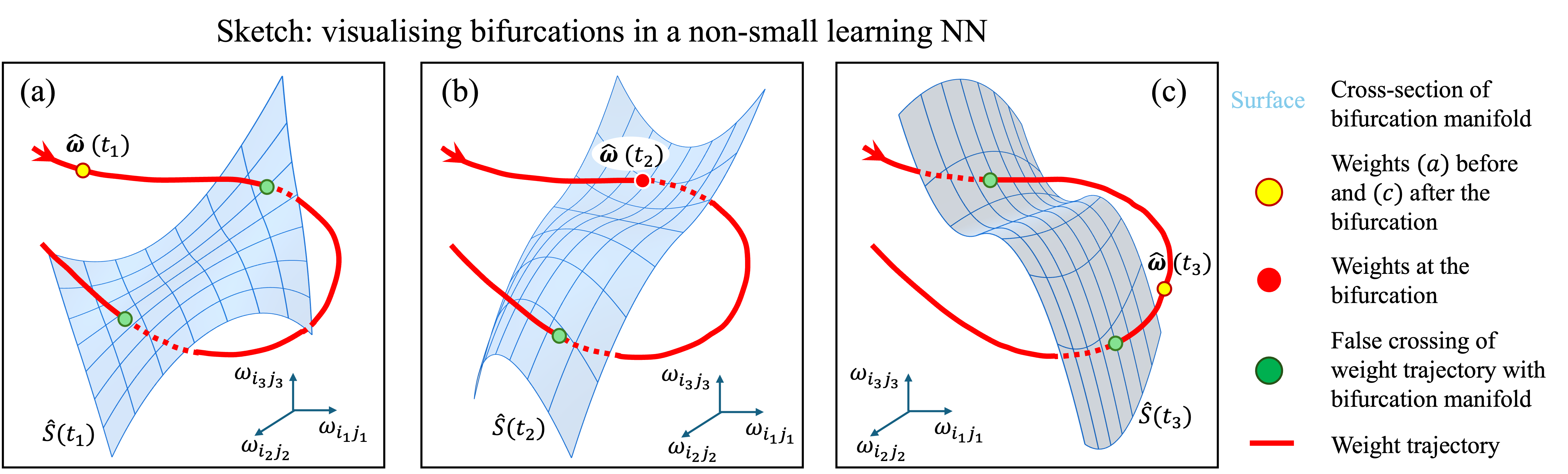}
\caption{A sketch illustrating a method to visualise bifurcation mechanisms of memory formation and forgetting in a non-small learning NN. Panels show cross-sections $\widehat{S}(t_n)$ (blue surfaces)  of the \emph{same} bifurcation manifold $S$, which exists in the $M$-dimensional weight space of (\ref{eqn:eq3.1}) with $N$$>$$3$, 
with different three-dimensional flat surfaces defined by $\omega_{ij}$$=$$\omega_{ij}(t_n)$ for all  $\omega_{ij}$  (except $\omega_{i_1 j_1}$, $\omega_{i_2 j_2}$ and $\omega_{i_3 j_3}$) at the values they take in (\ref{eqn:eq2.1})--(\ref{eqn:eq2.4}) at time $t_n$. 
In the three-dimensional space of selected weights,    for each $t$$=$$t_n$, there is a   \emph{different}  $\widehat{S}(t_n)$. 
Also shown is a projection $\widehat{\pmb{\omega}}(t)$ (red line) of the weight trajectory $\pmb{\omega}(t)$ of (\ref{eqn:eq2.4}). 
Panels illustrate the states (a) before the bifurcation at $t$$=$$t_1$ (yellow circle), (b) at the bifurcation at $t$$=$$t_2$ (red circle), and (c) after the bifurcation at $t$$=$$t_3$ (yellow circle), with $t_1$$<$$t_2$$<$$t_3$. 
The instant of bifurcation at $t$$=$$t_2$ is visualised in (b) as the crossing by $\widehat{\pmb{\omega}}(t)$ (red circle) of $\widehat{S}(t_2)$. In (a)--(c) other visible crossings between $\widehat{\pmb{\omega}}(t)$ and $\widehat{S}(t_n)$ (green circles) 
do \emph{not} represent real crossings between $\pmb{\omega}(t)$ and $S$, since each green circle corresponds to $t$ different from $t_n$, $n$$=$$1,2,3$, for which $\widehat{S}(t_n)$ is constructed.
}
\label{fig:Bif_Surf_sketch}
\end{figure*}

\begin{figure*}
\includegraphics[width=1.0\textwidth]{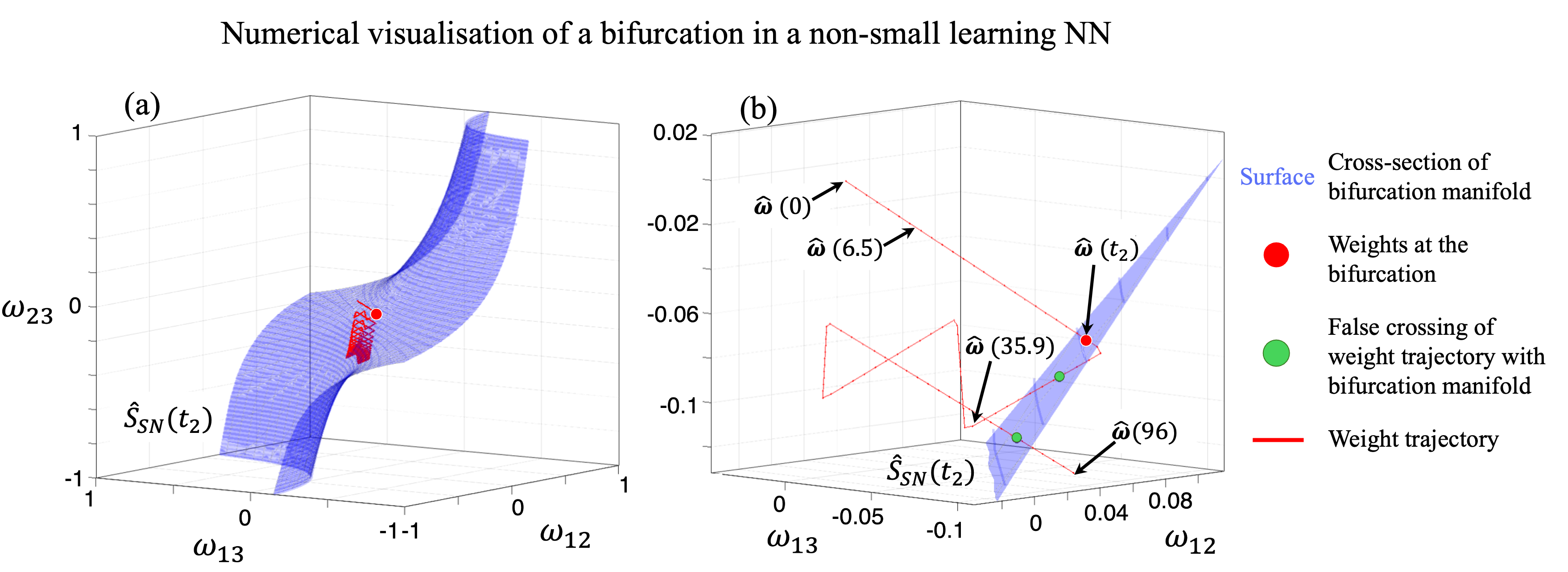}
\caption{Compare with Fig.~\ref{fig:Bif_Surf_sketch}(b). Numerical visualisation  of a bifurcation responsible for the formation of two potential memories (attractors with basins) in a \emph{large} NN  (\ref{eqn:eq2.1})--(\ref{eqn:eq2.4}) with $N$$=$$81$, $g$$=$$0.3$, $B$$=$$300$, $A$$=$$30$ learning from training Set~1 (Tabs.~\ref{supp-tab:tab1-1}--\ref{supp-tab:tab1-2} of Supplementary Note).  $\widehat{S}_{SN}(t_2)$ (blue surface) is a cross-section of a manifold $S_{SN}$ of a pair of saddle-node bifurcations in the NN (\ref{eqn:eq3.1}) with $N$$=$$81$ and $g$$=$$0.3$, visualised in the sub-space of weights $\left( \omega_{12},\omega_{13},\omega_{23} \right)$ at $t$$=$$t_2$$=$$22.8783$. To obtain $\widehat{S}_{SN}(t_2)$, all other  $\omega_{ij}$ were fixed at values they take at  $t$$=$$t_2$ as (\ref{eqn:eq2.1})--(\ref{eqn:eq2.4}) undergoes the first pair of saddle-node bifurcations marked by the first pair of green circles in Fig.~\ref{fig:81N_Bif_Diag}. Panel (b) is a zoomed version of (a). 
Also in (a)--(b) there is a projection $\widehat{\pmb{\omega}}(t)$ of the ``weight trajectory'' $\pmb{\omega}(t)$, $\widehat{\pmb{\omega}}(t)$$=$$\left( \omega_{12}(t),\omega_{13}(t),\omega_{23}(t) \right)$ (red line), of (\ref{eqn:eq2.4}). At $6.5$$<$$t$$<$$t_2$, the NN (\ref{eqn:eq2.1})--(\ref{eqn:eq2.4}) has two attractors (Fig.~\ref{fig:Pitchfork}(c)). At $t_2$$<$$t$$<$$35.9$,  it has four attractors (Fig.~\ref{fig:Pitchfork}(d)), two of which are born from the respective pair of saddle-node bifurcations. Panel (b) shows a portion of  $\widehat{\pmb{\omega}}(t)$ for $t$$\in$$[0,96]$, during which \emph{eight} input vectors $\mathbf{I}^k$ are applied consecutively with $k$$=$$1,2,3,4,5,6,1,2$. In (a)--(b) red circle at intersection between  $\widehat{\pmb{\omega}}(t)$ and $\widehat{S}_{SN}(t_2)$ signifies a real intersection between $\pmb{\omega}(t)$ and $S_{SN}$, i.e. the instant of bifurcation. In (b) green circles at other visible intersections between  $\widehat{\pmb{\omega}}(t)$ and $\widehat{S}_{SN}(t_2)$ do not signify real bifurcations, see Sec.~\ref{sec_method} and Fig.~\ref{fig:Bif_Surf_sketch}. 
}
\label{fig:81N_Bif_Surf}
\end{figure*}

Note, that this small three-neuron network develops only a maximum of two attractors while trying  -- and failing -- to memorise three patterns. This illustrates the well-known fact that a small NN cannot learn efficiently.  Moreover, in the course of learning,  due to a bifurcation, it 
temporarily loses all memories  it manages to form.  In the case when the \emph{full} bifurcation diagram in the space of weights  can be obtained, the above analysis explicitly demonstrates that, while the NN trains on a single task, forgetting is caused by the crossings between the weight trajectory of the learning NN and a bifurcation manifold of the counterpart non-learning NN.

\subsection{Forgetting via bifurcation in a large network: method and some evidence}
\label{sec:large}

To verify a hypothesis, that the causes of 
forgetting are the crossings by the weight trajectory of the bifurcation manifolds, ideally, one needs a bifurcation diagram in the space of \emph{all} weights of  the autonomous NN  (\ref{eqn:eq3.1}). However, while it is theoretically possible to perform a bifurcation analysis of an $81$-neuron network in the $3240$-dimensional space of weights, doing so would not only be extremely challenging computationally, but the results would be impossible to visualise using standard tools or interpret in an intuitively accessible way.

In Sec.~\ref{sec_method} we present a possible form of an explicit evidence that could support the above hypothesis when a complete bifurcation diagram is unavailable. Given that one can feasibly visualise in a three-dimensional space \emph{cross-sections} of the weight space and of the  bifurcation manifolds, consider what these cross-sections would show when a bifurcation takes place during learning. 

\subsubsection{Crossing a bifurcation manifold: a detection method}
\label{sec_method}

Assume that there exists some bifurcation manifold $S$ in the space of \emph{all} weights of (\ref{eqn:eq3.1}). As its counterpart NN (\ref{eqn:eq2.1})--(\ref{eqn:eq2.4}) learns and $t$ evolves, $S$ remains fixed. For any value of $t$, at which the respective bifurcation does \emph{not} occur, the current state $\pmb{\omega}(t)$ is at a finite distance from $S$. Only at the instant of bifurcation, $\pmb{\omega}(t)$ crosses $S$. For a small network, this situation is visualised in Fig.~\ref{fig:fig9.1}. 

Now, consider \emph{cross-sections} of $S$ in a non-small NN (\ref{eqn:eq3.1}) as its counterpart (\ref{eqn:eq2.1})--(\ref{eqn:eq2.4}) undergoes various stages of learning. In order to  visualise how $\pmb{\omega}(t)$ intersects with $S$, one needs to choose the cross-sections of $S$ meaningfully. Figure~\ref{fig:Bif_Surf_sketch} schematically illustrates the method we propose to achieve this. Assume that at some time $t$$=$$t_n$ during the learning phase, we wish to construct a cross-section  $\widehat{S}(t)\big|_{t=t_n}$$=$$\widehat{S}(t_n)$ of $S$ in a subspace of weights $\left( \omega_{i_1 j_1},\omega_{i_2 j_2},\omega_{i_3 j_3} \right)$. We propose to consider an intersection between $S$ and a three-dimensional flat surface in the $M$-dimensional weight space defined by the set of equations $\omega_{ij}$$=$$\omega_{ij}(t_n)$ for all  $\omega_{ij}$ except the three selected weights  $\omega_{i_1 j_1}$, $\omega_{i_2 j_2}$ and $\omega_{i_3 j_3}$. Alternatively, $\widehat{S}(t_n)$ can be understood as the result of slicing  $S$ by $(M$$-$$3)$ mutually orthogonal hyper-planes  of dimension $(M$$-$$1)$, each defined by  $\omega_{ij}$$=$$\omega_{ij}(t_n)$ for appropriate $\omega_{ij}$. 

Thus, as learning progresses, for  every new value $t$$=$$t_n$, there is a new  secant surface and therefore a \emph{new} cross-section $\widehat{S}(t_n)$ of the same $S$. This is illustrated in Fig.~\ref{fig:Bif_Surf_sketch}, where different cross-sections $\widehat{S}(t_n)$ (blue and grey surfaces) are schematically shown for three consecutive times $t_n$: (a) $t_1$, (b) $t_2$ (bifurcation), and (c) $t_3$, with $t_1$$<$$t_2$$<$$t_3$. 

Now, consider a projection $\widehat{\pmb{\omega}}(t)$ (red line in Fig.~\ref{fig:Bif_Surf_sketch}, the same in (a)--(c)) of the weight trajectory $\pmb{\omega}(t)$ onto the $\left( \omega_{i_1 j_1},\omega_{i_2 j_2},\omega_{i_3 j_3} \right)$-space relative to $\widehat{S}(t_n)$. At each $t$$=$$t_n$ when there is no bifurcation (Fig.~\ref{fig:Bif_Surf_sketch}(a) and (c)), the current state $\widehat{\pmb{\omega}}(t_n)$ (yellow circle) is at a finite distance from the respective $\widehat{S}(t_n)$. With this,  there can exist intersections between the \emph{full} $\widehat{\pmb{\omega}}(t)$ and  the instantaneous cross-section $\widehat{S}(t_n)$ (green circles in Fig.~\ref{fig:Bif_Surf_sketch}(a) and (c)). However, these intersections do not signify the actual crossings in the space of \emph{all} weights between $\pmb{\omega}(t)$ and $S$, and hence do not mark bifurcations. This is because for any point $\widehat{\pmb{\omega}}(t^*)$ on $\widehat{\pmb{\omega}}(t)$, at which the latter crosses $\widehat{S}(t_n)$ with  $t^*$$\ne$$t_n$, its counterpart cross-section $\widehat{S}(t^*)$ is generally \emph{different} from $\widehat{S}(t_n)$. 

Only at an instant of the bifurcation, the \emph{current} state $\widehat{\pmb{\omega}}(t_n)$ lands on the respective $\widehat{S}(t_n)$ (red circle in Fig.~\ref{fig:Bif_Surf_sketch}(b) corresponding to $t$$=$$t_2$). Other possible intersections between $\widehat{\pmb{\omega}}(t)$ and $\widehat{S}(t_2)$ (green circle in Fig.~\ref{fig:Bif_Surf_sketch}(b)) would not indicate bifurcations, as explained above. 

\subsubsection{Crossing a bifurcation manifold in a large NN: numerical evidence}
\label{sec_num_evid}

In Sec.~\ref{sec:bif} we observed that in the NN (\ref{eqn:eq2.1})--(\ref{eqn:eq2.4}) with $N$$=$$81$, the primary mechanism of memory formation is the saddle-node bifurcation. Thus, for the NN learning from Set~1, as parameter $t$ increases, we focus on the first pair of saddle-node bifurcations marked by the first pair of green circles in Fig.~\ref{fig:81N_Bif_Diag}, illustrated by phase portraits in Fig.~\ref{fig:Pitchfork}(c)--(d), and occurring  at $t$$=$$22.8783$. 
 
Following the method described in Sec.~\ref{sec_method}, below we visualise the instant when the weight trajectory $\pmb{\omega}(t)$ of the large NN crosses the manifold $S_{SN}$ of the saddle-node bifurcation. In other words, we obtain a graph conveying the message of Fig.~\ref{fig:Bif_Surf_sketch}(b) for the NN under study. To do this, in (\ref{eqn:eq3.1}) we  fix all $\omega_{ij}$, except $\omega_{12}$, $\omega_{13}$ and $\omega_{23}$, at the values they take in (\ref{eqn:eq2.1})--(\ref{eqn:eq2.4}) at $t$$=$$22.8783$$=$$t_2$.  We then perform a three-parameter bifurcation analysis of  (\ref{eqn:eq3.1}) as $\omega_{12}$, $\omega_{13}$ and $\omega_{23}$ are allowed to vary freely and independently of each other. 
 
This way, in the space $(\omega_{12},\omega_{13},\omega_{23})$ we obtain the cross-section $\widehat{S}_{SN}(t_2)$ (Fig.~\ref{fig:81N_Bif_Surf}, blue surface) between $S_{SN}$ and a three-dimensional flat surface in the $M$-dimensional space (with $M$$=$$3240$) defined by the set of equations $\omega_{ij}$$=$$\omega_{ij}(t_2)$ for all $\omega_{ij}$  excluding the three selected weights.
 
 Specifically, to construct the cross-section $\widehat{S}_{SN}(t_2)$, we sample the values of $\omega_{23}$ from $[-1,1]$ with a step $0.1$. For each fixed $\omega_{23}$, a one-dimensional bifurcation curve was obtained in $(\omega_{12},\omega_{13})$-plane using XPPAUT \cite{xppaut07}. A collection of these curves was placed in  $(\omega_{12},\omega_{13},\omega_{23})$-space and used as a frame to render the surface $\widehat{S}_{SN}(t_2)$. 

Next, we register the weight trajectory of (\ref{eqn:eq2.4}), i.e. $\pmb{\omega}(t)$$=$$\left(\omega_{12}(t), \ldots,\omega_{(N-1)N}\right(t))$, such that $\omega_{ii}$$=$$0$,  $\omega_{ij}$$=$$\omega_{ji}$, $\forall i,j$$=$$1,$$\dots$$N$. The  components of $\pmb{\omega}(t)$ are given in Fig.~\ref{fig:omega_t}(a), and this $\pmb{\omega}(t)$ was used to produce Fig.~\ref{fig:81N_Bif_Diag}. Finally, we superimpose the projection $\widehat{\pmb{\omega}}(t)$ of $\pmb{\omega}(t)$ onto the $(\omega_{12},\omega_{13},\omega_{23})$-space (red line in Fig.~\ref{fig:81N_Bif_Surf}(a)--(b)) with $\widehat{S}_{SN}(t_2)$.

Along the portion of $\widehat{\pmb{\omega}}(t)$ visible slightly to the left of  $\widehat{S}_{SN}(t_2)$ in Fig.~\ref{fig:81N_Bif_Surf}(b), i.e. at $6.5$$<$$t$$<$$t_2$, the NN (\ref{eqn:eq3.1}) has two attractors (see Fig.~\ref{fig:Pitchfork}(c)). At $t$$=$$22.8783$$=$$t_2$, $\widehat{\pmb{\omega}}(t)$ crosses $\widehat{S}_{SN}(t_2)$ (red circle in Fig.~\ref{fig:81N_Bif_Surf}), thus signifying the crossing between $\pmb{\omega}(t)$ and $S_{SN}$ and the occurrence of the pair of saddle-node bifurcations giving rise to two more attractors.  At $t_2$$<$$t$$<$$35.9$, the NN (\ref{eqn:eq3.1}) has four attractors altogether  (see Fig.~\ref{fig:Pitchfork}(d)). 
  
At $t$$=$$35.9$ two pairs of saddle-node bifurcations occur -- to some other fixed points. Note, that the location of $\widehat{\pmb{\omega}}(35.9)$ in Fig.~\ref{fig:81N_Bif_Surf}(b) provides a rough idea of where the cross-section of the manifold of the latter bifurcation is located at $t$$=$$35.9$.
  
The additional crossings between $\widehat{\pmb{\omega}}(t)$ and $\widehat{S}_{SN}(t_2)$ (green circles in Fig.~\ref{fig:81N_Bif_Surf}(b)) do not signify bifurcations, as explained in Sec.~\ref{sec_method}.

As $t$ grows to very large values, $\widehat{\pmb{\omega}}(t)$ converges to a closed curve, see condensation of red lines in the lower part of $\widehat{\pmb{\omega}}(t)$ in Fig.~\ref{fig:81N_Bif_Surf}(a). This is qualitatively similar to the behaviour of  $\pmb{\omega}(t)$ in the small NN illustrated by Fig.~\ref{fig:fig9.1}(a). Thus, towards the end of learning, $\omega_{23}$ oscillates around $-\frac{1}{3}$, whereas $\omega_{12}$ and $\omega_{13}$ oscillate around zero (compare with Fig.~\ref{fig:omega_t}(a)).

\subsubsection{Forgetting via bifurcation in single-task training: evidence-based inference}
\label{sec_inference}

When a complete bifurcation diagram of a NN in the space of its weights is unavailable, an explicit numerical proof, that 
forgetting occurs when the weight trajectory crosses a bifurcation surface, can in principle be obtained by constructing sequences of appropriate cross-sections of the bifurcation manifolds. The method to achieve this is described in Sec.~\ref{sec_method} and illustrated for an $81$-neuron network in Sec.~\ref{sec_num_evid}.  

Recall that the Hopfield NN considered here starts from zero memories and acquires new memories as its weight trajectory goes through the space of weights and crosses various bifurcation manifolds. One example of such a crossing is explicitly demonstrated in Fig.~\ref{fig:81N_Bif_Surf}. 
Forgetting means that attractors, already formed in (\ref{eqn:eq2.1})--(\ref{eqn:eq2.4}) by a certain stage of learning, should suddenly disappear as $t$ keeps increasing continuously. This implies that as $t$ grows,  the weight trajectory must cross the \emph{same} bifurcation manifolds it crossed previously, only in the opposite direction, as demonstrated for the small NN  in Fig.~\ref{fig:fig9.1}. 

To explicitly demonstrate this effect in a large NN, even for a single instance of 
forgetting, one needs to construct a series of cross-sections of a certain bifurcation manifold $S$ at several consecutive values of $t$. These values of $t$ should include two instants of the respective bifurcation (when the same attractors are first born and then die),  and in addition  instants before the first bifurcation, in between the two bifurcations, and after the second bifurcation. 

Numerical calculations of such cross-sections for non-small NNs are technically challenging, laborious and computationally demanding. Here we calculate only one of the required cross-sections in Fig.~\ref{fig:81N_Bif_Surf}, whereas obtaining the complete and explicit numerical proof of the suggested mechanism of 
 forgetting while the NN trains on a single task could be done in the future with the same approach. 

\begin{figure}
\includegraphics[width=0.45\textwidth]{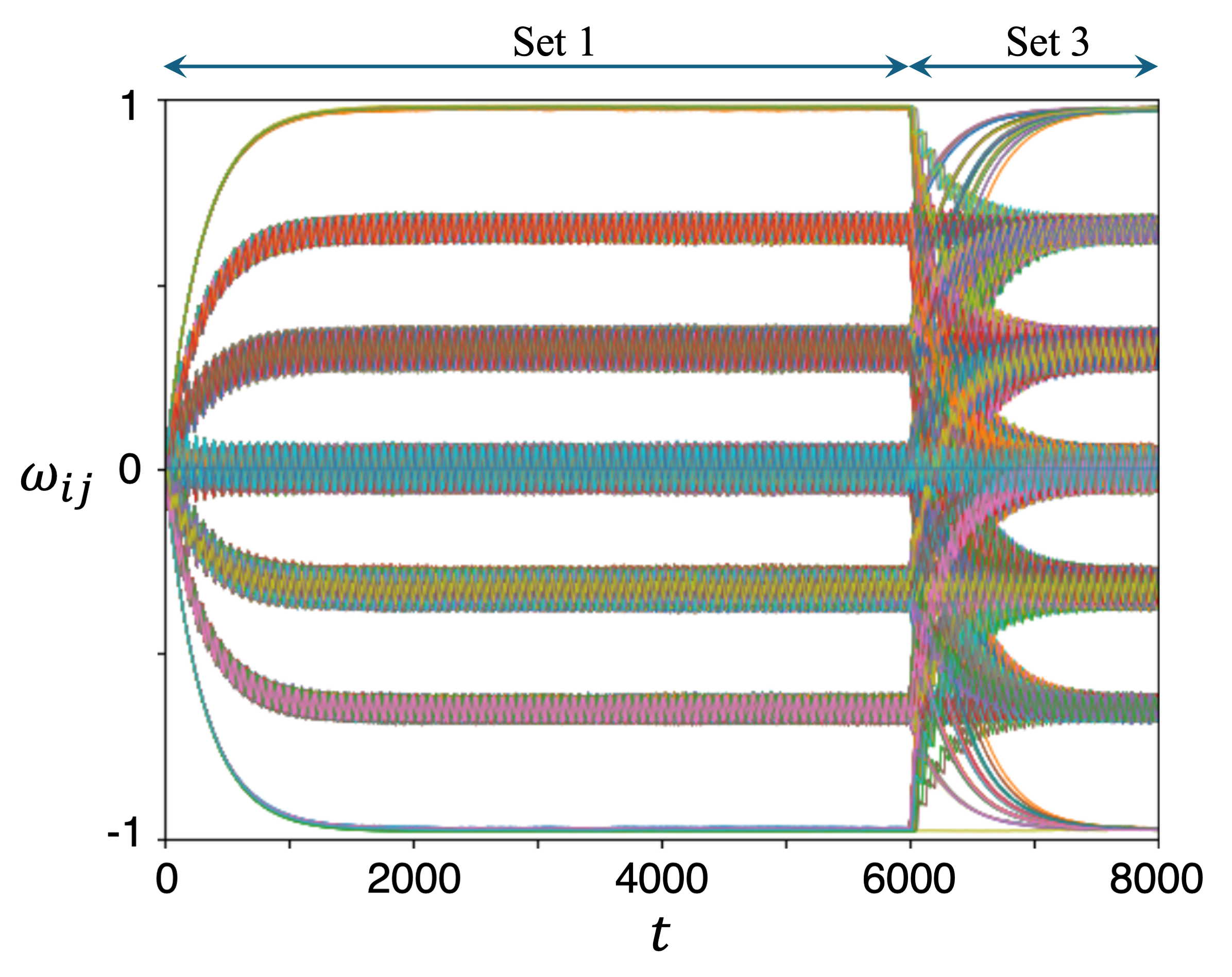}
\caption{Evolution of connection weights $\omega_{ij}(t)$ (solid lines in various colours) during learning by the NN (\ref{eqn:eq2.1})--(\ref{eqn:eq2.4}) 
as it initially learns the task represented by training Set~1 (see Tabs.~\ref{supp-tab:tab1-1}--\ref{supp-tab:tab1-2}), and then \emph{switches} to learning the task represented by Set~3 (see Tabs.~\ref{supp-tab_new_task_1}--\ref{supp-tab_new_task_2}). 
Learning from Set~1 occurs for $t\in[0,6000]$, and the respective part of the graph is identical to Fig.~\ref{fig:omega_t}(a). Learning from Set~3 occurs for $t$$\in$$[6000,8000]$. Parameters (same as in Fig.~\ref{fig:omega_t}) are: $N$$=$$81$, $A$$=$$30$, $B_{ij}$$=$$B$$=$$300$ $(\forall i,j)$, $g$$=$$0.3$ and $\lambda$$=$$1.4$. \emph{Catastrophic forgetting} of memories, which were formed during training on Set~1  and then destroyed while training on Set~3, are illustrated in Figs.~\ref{fig:pp_cat} and \ref{fig:bd_cat}. 
}
\label{fig:weights_cat}
\end{figure}

However, the numerical evidence collected here allows us to make an inference regarding the validity of our hypothesis. Firstly, we discover that throughout learning, the weight trajectory oscillates (zig-zags) with a considerable amplitude due to the applied periodic stimulus $\mathbf{I}(t)$. Secondly, we observe that in the space of weights, the bifurcation manifolds are probably located quite close to each other. The latter is evidenced by the rather small distance in the subspace of weights between cross-sections of two different bifurcation manifolds: $\widehat{S}_{SN}(t_2)$ and of another saddle-node bifurcation, which must go through the point $\widehat{\pmb{\omega}}(35.9)$ in Fig.~\ref{fig:81N_Bif_Surf}(b). Thirdly, 
Fig.~\ref{fig:omega_t} illustrates the well-known fact that,  if the initial weights are inside the box  $\left| \omega_{ij} \right| <1$, the whole weight trajectory is confined to the same box. Fourthly, by the end of learning, the weight trajectory almost reaches some 
sides
of the box $\left| \omega_{ij} \right| <1$, since some weights tend to $\pm 1$ (see Fig.~\ref{fig:omega_t}). Fifthly,  bifurcation manifolds can stretch beyond the boundaries of the box $\left| \omega_{ij} \right| <1$ confining the weight trajectory, as seen in Fig.~\ref{fig:81N_Bif_Surf}(a). 

\begin{figure*}[ht]
\includegraphics[width=1.0\textwidth]{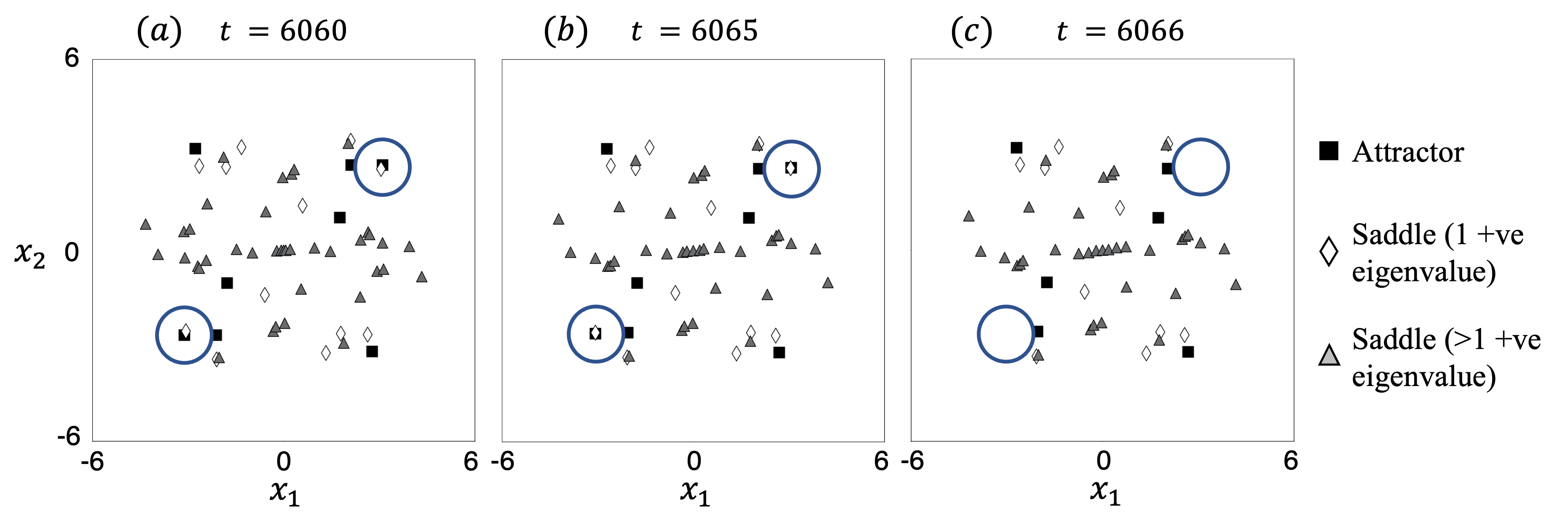}
\caption{Illustration of catastrophic forgetting of memories, which were formed in the NN (\ref{eqn:eq2.1})--(\ref{eqn:eq2.4}) after it was trained on Set~1, while it undergoes training on Set~3 (the respective weight evolution is shown in Fig.~\ref{fig:weights_cat}). 
Panels show projections onto the $(x_1,x_2)$-plane of fixed points of the NN (\ref{eqn:eq3.1}) at various values of the control parameter $t$ equal to time $t$ of  (\ref{eqn:eq2.1})--(\ref{eqn:eq2.4}). Blue circles highlight areas of the phase space where bifurcations take place.  Fixed points are marked as: attractors (black boxes), ``useful'' saddle points with one positive eigenvalue (white diamonds), and all other saddle or unstable points (grey triangles). (a) At $t$$=6060$, a true memory (black box inside the lower blue circle, also shown as the black diamond at the bottom left of Fig.~\ref{fig:fig3.1}(d)) and a spurious memory (black box inside the upper blue circle, also shown as the red cross at the top right of Fig.~\ref{fig:fig3.1}(d)) move close to saddle points with a single positive eigenvalue (white diamonds). (b) At $t$$=$$6065$ these two pairs of points collide in a pair of simultaneous saddle-node bifurcations (black boxes and white diamonds overlapping inside blue circles). (c) By $t$$=$$6066$ both attractors have disappeared abruptly, thus signifying catastrophic forgetting of the respective memories. 
 }
\label{fig:pp_cat}
\end{figure*}

The five observations above suggest that there is quite a high probability that a zig-zagging weight trajectory repeatedly crosses  in both directions the same bifurcation manifolds, which are rather tightly packed in the bounded volume of the space of weights. 
A backward crossing that destroys a previously created attractor signifies forgetting of the previously created memory.  Thus, the numerical results already available for the relatively large Hopfield NN with Hebbian learning provide a good evidence in favour of the hypothesis that, 
in the course of single-task training, forgetting  is caused by intersections between the weight trajectory of the learning NN and the bifurcation manifolds of the counterpart non-learning NN. 

\begin{figure*}[ht]
\includegraphics[width=0.85\textwidth]{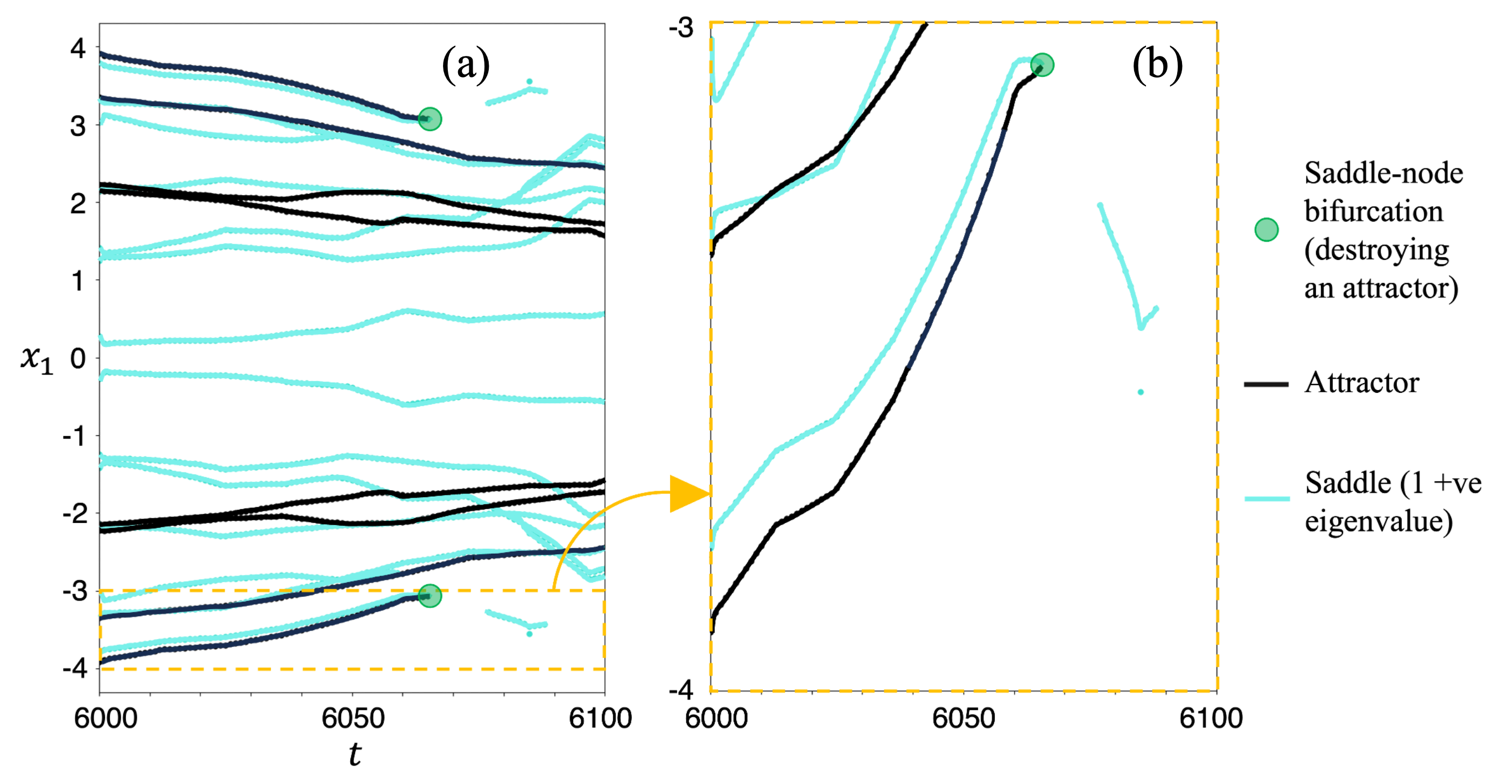}
\caption{One-parameter bifurcation diagram demonstrating the mechanism of catastrophic forgetting in the NN (\ref{eqn:eq2.1})--(\ref{eqn:eq2.4}) with $N$$=$$81$, while it is being trained on Set~3 after training on Set~1 (the respective weight evolution is shown in Fig.~\ref{fig:weights_cat}). Panels show $x_1$-coordinates of fixed points of (\ref{eqn:eq3.1}) (solid lines) as functions of control parameter $t$, which coincides with time $t$ in (\ref{eqn:eq2.1})--(\ref{eqn:eq2.4}). Stability of fixed points is indicated by the line colour: stable point, i.e. attractor  potentially associable with memory (black), and ``useful''  saddle point with a single positive eigenvalue (blue). 
Green circles mark a pair of saddle-node bifurcations at $t$$\approx$$6065$ between attractors and saddles with a single positive eigenvalue, which occur simultaneously due to symmetry in (\ref{eqn:eq3.1}).  Panel (a) shows the bifurcation diagram for $t \in [6000,6100]$. Panel (b) is a close-up of the rectangular selection in (a) and shows the bifurcation destroying the attractor representing a true memory learnt during training on Set~1. These bifurcations are additionally illustrated by phase portraits in Fig.~\ref{fig:pp_cat}. 
 }
\label{fig:bd_cat}
\end{figure*}

\section{Catastrophic forgetting via bifurcation}

\label{sec_cat} 

Here we consider the Hopfield NN with Hebbian learning experiencing catastrophic forgetting in the conventional sense, namely, the abrupt loss of memories when the network switches to being trained on a new task after the end of training on the previous task.
We demonstrate that, just like forgetting during single-task training considered in Sec.~\ref{sec:forget}, conventional catastrophic forgetting is also caused by bifurcations.

Consider the NN (\ref{eqn:eq2.1})--(\ref{eqn:eq2.4}) with $N$$=$$81$ after it has undergone training on Set~1, given in Tabs.~\ref{supp-tab:tab1-1}--\ref{supp-tab:tab1-2}, and developed a number of true and spurious memories, as demonstrated in Sec.~\ref{sec:bif} (see Fig.~\ref{fig:omega_t}(a) for the evolution of weights and Fig.~\ref{fig:fig3.1}(d) for the memories formed). Now launch (\ref{eqn:eq2.1})--(\ref{eqn:eq2.4}) at $t$$=$$t_1$$=$$6000$  from initial conditions $x_i(t_1)$ and $\omega_{ij}(t_1)$, and apply the stimulus $\mathbf{I}(t)$ formed from a new training Set~3 given in Tabs.~\ref{supp-tab_new_task_1}--\ref{supp-tab_new_task_2}. The evolution of weights as the NN is being consecutively trained on Set~1 and then on Set~3 is illustrated in Fig.~\ref{fig:weights_cat}.

As the NN processes Set~3, let us follow the evolution of a selected memory acquired from training on Set~1 and represented by an attractor (see Fig.~\ref{fig:fig3.1}(d), black diamond in the bottom left corner) with its basin. Note, that a symmetrically positioned attractor (see Fig.~\ref{fig:fig3.1}(d), red cross in the top right corner) with its own basin represents a spurious memory. In Fig.~\ref{fig:pp_cat}(a) the same attractors are shown at $t$$=$$6060$ as black boxes inside blue circles -- however, their locations have slightly shifted as compared to Fig.~\ref{fig:fig3.1}(d), since they correspond to a different set of $\omega_{ij}$. 

Figure~\ref{fig:pp_cat}(a) illustrates a situation just before a pair of saddle-node bifurcations involving these selected attractors, when the attractors (black boxes) come close to the saddle fixed points with a single positive eigenvalue (white diamonds inside blue circles). At the saddle-node bifurcation at $t$$=$$6065$ (Fig.~\ref{fig:pp_cat}(b)) the two pairs of stable and saddle fixed points collide, and after the bifurcation at $t$$=$$6066$ (Fig.~\ref{fig:pp_cat}(c)) they no longer exist. The disappearance of these attractors signifies the conventionally understood \emph{catastrophic forgetting} of two memories: a true and a spurious one. 

For the completeness of the argument, Figure~\ref{fig:bd_cat}(a) shows a segment of a one-parameter bifurcation diagram of the NN (\ref{eqn:eq2.1})--(\ref{eqn:eq2.4}) while, after the end of training on Set~1, it is trained on Set~3 for $t$$ \in$$ [6000,6100]$. Here, black lines show locations of attractors, and blue lines locations of ``useful'' saddle points with a single positive eigenvalue. A pair of green circles indicate the pair of simultaneous saddle-node bifurcations occurring to the attractors of our interest (compare with Fig.~\ref{fig:pp_cat}).  Figure~\ref{fig:bd_cat}(b) shows an enlarged segment of the diagram in (a) bounded by yellow dashed line, which focuses on the disappearance of the true memory. 

Note, that for $t$$ \in$$ [6000,6100]$ there exist  saddle fixed points with more than one positive eigenvalue, but they are not shown in Fig.~\ref{fig:bd_cat} in order to avoid cluttering. Fixed points with two positive eigenvalues can participate in saddle-node bifurcations with ``useful'' saddles (blue lines).  Evidence of two pairs of these bifurcations can be seen in Fig.~\ref{fig:bd_cat}(a), where two symmetric short segments of blue lines exist only for $t \in [6077,6088]$. Namely, the respective saddles appear at $t$$=$$6077$ as if out of nowhere, and vanish at $t$$=$$6089$. However, the actual reasons for their appearances and disappearances are two pairs of saddle-node bifurcations with saddles having two positive eigenvalues each, which are not shown in  Fig.~\ref{fig:bd_cat}.

To summarise, in this Section we explicitly demonstrate that the conventionally understood catastrophic forgetting is caused by the bifurcations in the learning NN. Therefore, there is a conceptual connection between the two cases of forgetting in NNs, namely, forgetting during single-task training and the conventional catastrophic forgetting.

\section{Basins of attraction}
\label{sec:bas}

\begin{figure*}[ht]
\includegraphics[width=0.85\textwidth]{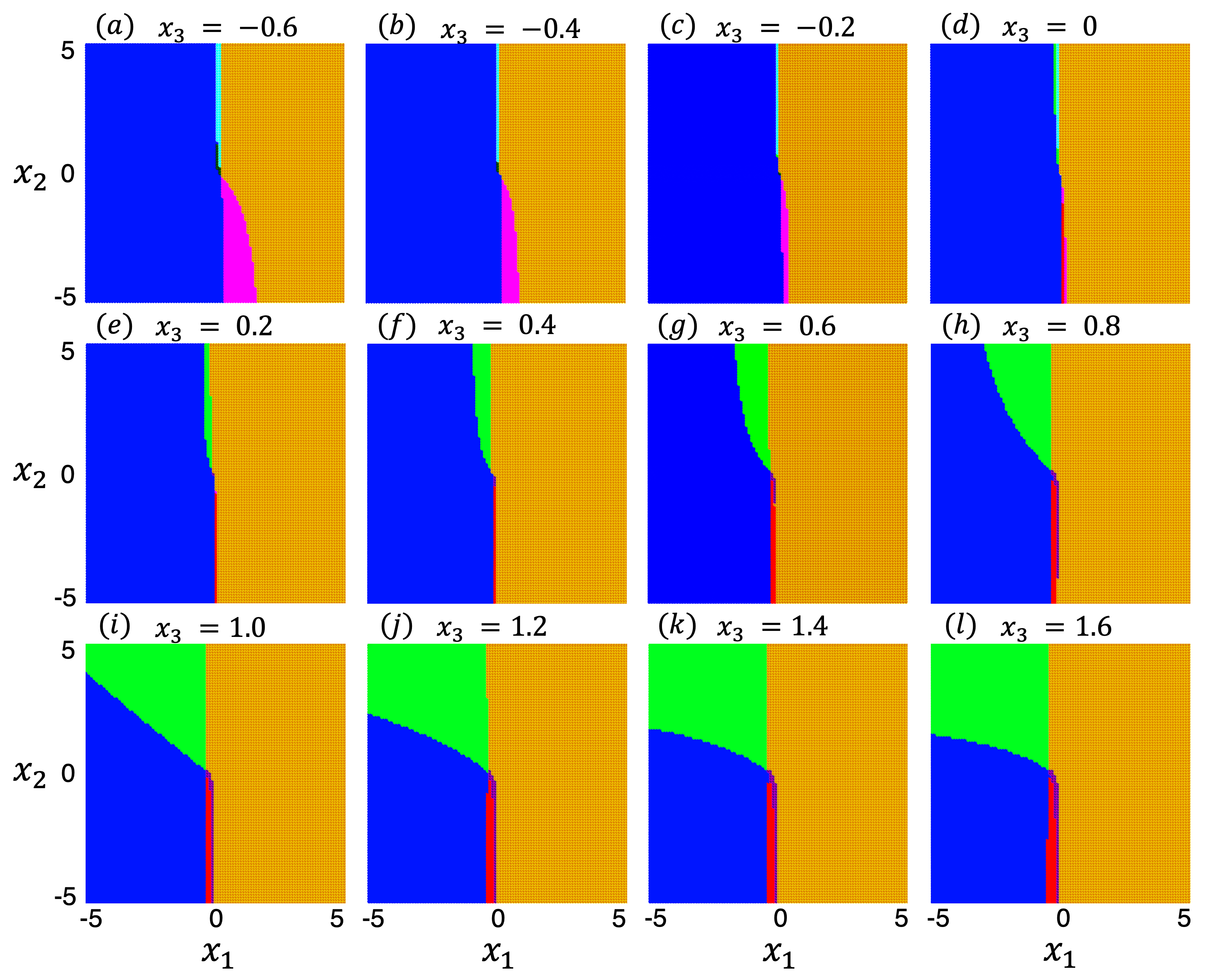}
\caption{Attractor basins (shaded regions in various colours) developed in the NN (\ref{eqn:eq2.1})--(\ref{eqn:eq2.4}) with 
$N$$=$$81$, $g$$=$$0.3$, $B$$=$$300$ and  $A$$=$$30$, by the end of learning from Set~1 at $t$$=$$6000$. Twelve different two-dimensional cross-sections (a)--(l)  of the same set of basins are shown as the phase space intersects with different planes defined by (\ref{eq_sec}), where $x_3^*$ is displayed above each panel. Colours of attractor basins match those of respective attractors in Fig.~\ref{fig:attractor_colours}.  Attractors and their basins are found for the autonomous version (\ref{eqn:eq3.1}) of the NN, where $\omega_{ij}$ are set to values achieved by (\ref{eqn:eq2.1})--(\ref{eqn:eq2.4}) at the end of learning.
 }
\label{fig:basins}
\end{figure*}

\begin{figure}[ht]
\includegraphics[width=0.5\textwidth]{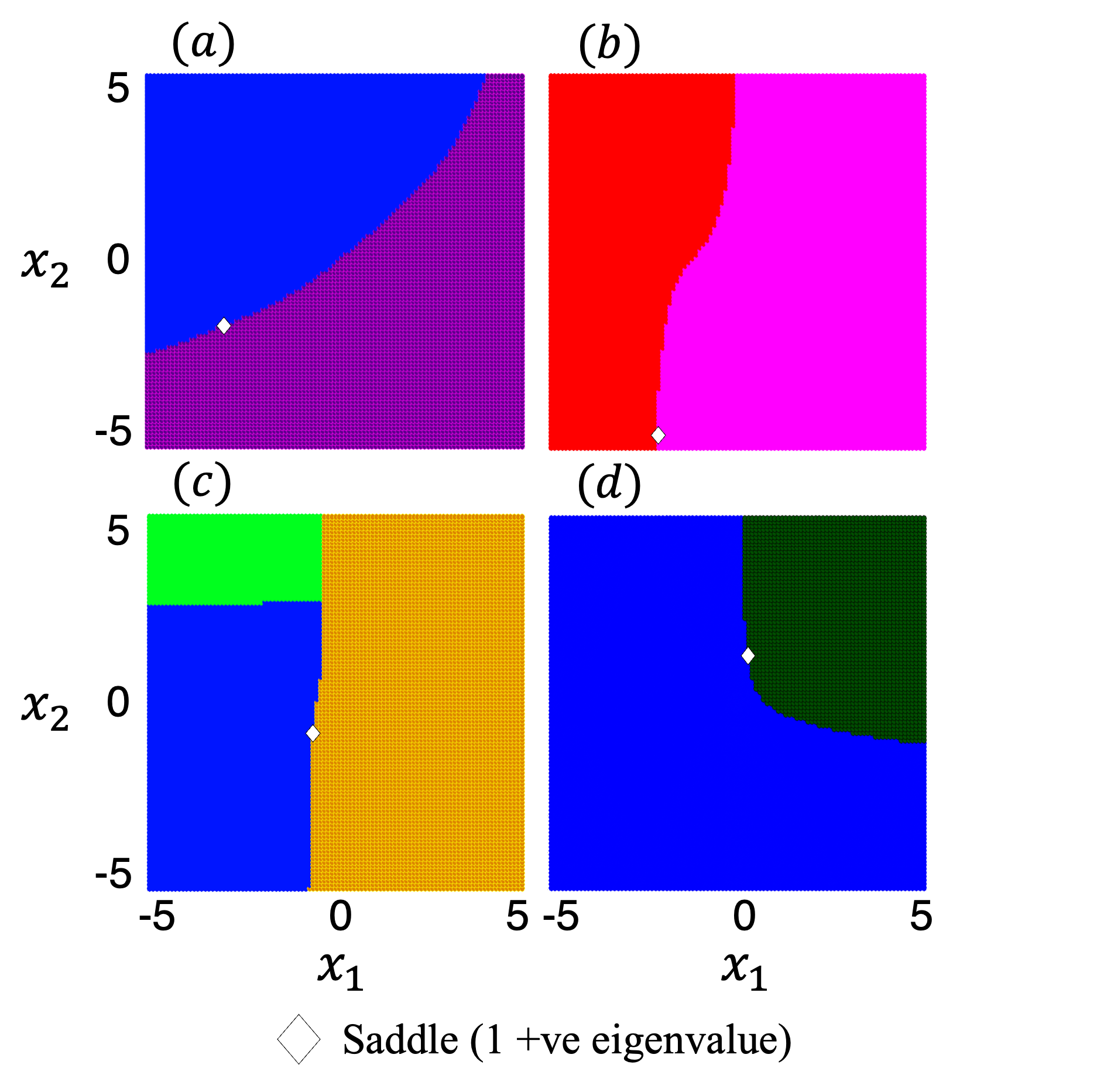}
\caption{Evidence that the boundaries of attractor basins are formed by the stable manifolds of the saddle fixed points with one positive eigenvalue (white diamonds). Here, more cross-sections are shown of the same set of attractor basins as those illustrated in Fig.~\ref{fig:basins}. 
Namely, every panel shows a cross-section with a different plane, each going through a different saddle point (white diamond), see details in text. Note, that the saddles lie exactly on the basin boundaries. 
Colours of attractor basins (shaded regions) match those of respective attractors in Fig.~\ref{fig:attractor_colours} and of their basins in Fig.~\ref{fig:basins}.
 }
\label{fig:basins_sad}
\end{figure}

The question about representation of memories in NNs is fundamental, but remains largely unresolved. To date the most popular idea about the way an ANN represents a \emph{single} memory is that it is represented by an attractor in the phase space of the NN  \cite{Amari1972, Hopfield82, Hopfield_NN}. However, when the ANN is used to recognise a previously unseen input as belonging to a certain \emph{category}, it is more appropriate to associate this category with the whole basin of attraction   \cite{Pineda_Hopfield_NN_bifurcation_JC88,Barack_Two_views_on_cognitive_brain_NatRevNeurosci21}. 

Therefore, the question about the size and the shape of attractor basins formed in a NN during learning becomes very important. Some investigations on the boundaries of attraction basins in small two- or three-neuron networks have been reported in Refs.~\cite{Atencia:2003aa,Krauth_basins_attraction_discrete_perceptron_NN_CS88}. For large discrete-time and discrete-state  NNs, it can be possible to evaluate some statistical characteristics of the radiuses of attraction basins \cite{Storkey_Hopfield_NN_dsicrete_basins_stat_NN99,Davey_Hopfield_NN_discrete_basin_size_99,Zhang_Hopfield_NN_discrete_basin_radius_08,Lin_Hopfield_NN_basin_radius_stat_conf23,Sampath_NN_discrete_basins_PLOS20}. However,  detailed analysis of the sizes and shapes of attraction basins in continuous-state NNs with more than three neurons has not been done to the best of our knowledge. 

Here we present an insight into the structure of the attraction basins of all attractors (Fig.~\ref{fig:attractor_colours}) formed by the end of learning in the NN (\ref{eqn:eq2.1})--(\ref{eqn:eq2.4}) with $N$$=$$81$ using training Set~1. We visualise the attractor basins in several \emph{two-dimensional} cross-sections of the 
phase space of (\ref{eqn:eq3.1})  with the planes defined by 
\begin{equation}
\label{eq_sec}
x_3=x_3^*, \ x_4=0, \ \ldots, \ x_{81}=0, 
\end{equation}
where $x_3^*$ is a value in $[-0.6,1.6]$ sampled  with step $0.2$.

In practice, to numerically calculate the basin cross-sections, we prepare a grid of initial conditions (ICs) on the secant plane representing the plane (\ref{eq_sec}). Namely, on $(x_1,x_2)$-plane we split a region $x_{1,2}$$\in$$[-5,5]$ into square boxes with the side $0.1$ to obtain $101$$\times$$101$ nodes. The vectors of ICs are then formed by taking $(x_1,x_2)$-coordinates of these nodes  and augmenting them by the values of $x_3$, $\ldots$, $x_{81}$ defined by (\ref{eq_sec}). 

For each IC,  we register an attractor to which the phase trajectory eventually converges. Then on the secant plane represented by  $(x_1,x_2)$-coordinates (since all other coordinates are fixed as in (\ref{eq_sec})) we mark this IC by the colour matching the colour of the respective attractor in Fig.~\ref{fig:attractor_colours}. 

The cross-sections of the same set of  attractor basins are shown in Fig.~\ref{fig:basins}, where different panels correspond to secant planes (\ref{eq_sec}) that differ only in the values of $x_3$, as indicated above each panel. Note, that attractors themselves cannot be seen in the cross-sections considered, since they do not belong to the respective secant planes. The gradual increase of $x_3$ from (a) to (l) allows one to reveal the structure of the basins in a \emph{three-dimensional} cross-section of the $81$-dimensional phase space. The incremental reshaping of the two-dimensional cross-sections of the basin boundaries with the gradual increase of $x_3$ from (a) to (l) is evidence of their continuity in the full space. One can also observe the complex shapes of these boundaries, and notice different basin sizes. 

Interestingly, the basins with the largest sizes in the chosen cross-sections (blue and yellow shades) belong to the symmetric pair of attractors born first in the course of training from the pitchfork bifurcation marked by the first pink circle in Fig.~\ref{fig:81N_Bif_Diag}(a). Note the symmetry of their basins visible in Fig.~\ref{fig:basins}(d), which corresponds to $x_3$$=$$0$, implying that the secant plane goes through the origin.  The largest visible sizes of these attractor basins could be explained by the fact that our secant planes are quite  close to the origin, in whose vicinity the first attractors emerged. Other attractors are born later and further away from the origin, so their basins occupy much smaller areas in the given cross-sections. 

It is important to understand how the \emph{boundaries} of attraction basins (or separatrices) are formed. 
They should be  represented by some $(N$$-$$1)$-dimensional  manifolds. Given that the only limit sets in the NN (\ref{eqn:eq3.1}) are fixed points, suitable manifolds should belong to the saddle fixed points with a single positive eigenvalue marked by diamonds in Fig.~\ref{fig:attractor_colours}.

To verify this, we calculate a few more cross-sections of the same basins of attraction as those illustrated in Fig.~\ref{fig:basins}, but this time obtained with the secant planes going directly through (as close as the numerical error allows) various saddle points with exactly one positive eigenvalue.  Namely, we consider cross-sections of the $81$-dimensional space with four different planes defined by $x_i$$=$$s^m_i$, $i$$=$$3,$$\ldots$$,81$, $m$$=$$1,$$\ldots$$,4$, where $s^m_i$ is the $i^{\textrm{th}}$ coordinate of the $m^{\textrm{th}}$ saddle point chosen from points shown in  Fig.~\ref{fig:attractor_colours} as orange diamonds. If our assumption is true, in the given cross-sections, the respective saddle points should belong to the basin boundaries.  

In Fig.~\ref{fig:basins_sad}(a)--(d) the resultant four cross-sections of attractor basins are shown together with the respective saddle points $\mathbf{s}^m$ (white diamonds). The colour-coding of the basins is the same as in Fig.~\ref{fig:basins}. One can see that, within the numerical error, the saddles lie on the basin boundaries, which therefore must coincide with the stable manifolds of these saddles.

In the Supplementary Note, in Sec.~\ref{supp-sec:alt}, we show the attractors (Fig.~\ref{supp-fig:mem_alt}(d)) and attractor basins (Figs.~\ref{supp-fig:basins_alt} and \ref{supp-fig:basins_sad}) formed  in the Hopfield NN in an \emph{alternative} form (\ref{supp-eqn:Hop_Alt})--(\ref{supp-eqn:Hop_Alt3}) after it completes Hebbian learning from Set~1. For both forms of the Hopfield NN, basins are shown in the same cross-sections (compare Figs.~\ref{fig:basins} and \ref{supp-fig:basins_alt}), and appear quite similar in shape, albeit not identical. Figure~\ref{supp-fig:basins_sad} provides additional evidence that the basin boundaries are formed by the stable manifolds of saddle fixed points with a single positive eigenvalue.

\section{Summary}
\label{sec:summary}

We comprehensively and systematically studied  the
mechanisms of memory formation 
and memory loss in an archetypal relatively \emph{large} recurrent artificial neural network (ANN) --  the continuous-time Hopfield neural network (NN) with Hebbian learning. 
For this NN, we verified the previously formulated idea that memories in learning NNs can be formed thanks to bifurcations.   We also hypothesised and proved that the abrupt loss of previously acquired memories either in the course of its training on a single task, or when the NN switches to training on the next task after being trained on the previous task -- known in machine learning (ML) as catastrophic forgetting -- is also caused by bifurcations.

To verify all of the above hypotheses,  we needed to rigorously formulate the problem of detecting bifurcations in a \emph{learning} NN in terms of a continuously changing parameter representing the stage of learning. A bifurcation occurs when the weight trajectory of the learning  NN (non-autonomous dynamical system (DS)) crosses a bifurcation manifold of the non-learning NN (autonomous DS).

During single-task training, the weight trajectory can cross the same bifurcation manifolds in both directions. Crossing in one direction produces new attractors, whereas crossing in the other direction destroys the previously formed attractors. Our analysis demonstrates that, when the NN starts learning with no attractors potentially associable with memories, the \emph{same} bifurcation mechanism underlies two opposite effects: memory formation and memory loss.

In the context of catastrophic forgetting, at the end of training on the previous task, the weights settle into  small fluctuations around some point in the weight space. When training on the next task starts, application of the new stimuli pushes the weight trajectory away from that point and towards the bifurcation manifolds. Crossings of some bifurcation manifolds induce the birth of new attractors, which could be associated with new memories.  Crossings of other bifurcation manifolds destroy the attractors formed while training on the previous task.

We revealed the possible mechanism behind the birth of spurious memories. Namely, 
with identical neurons and symmetric weights, applications of \emph{individual} training vectors induce the simultaneous birth of \emph{several} attractors (pairs, quadruples, etc.), thus highlighting the fact that association of a single training vector with a single attractor is not straightforward. Unfortunately, a desirable scenario when application of one training vector could give rise to strictly one attractor, which could be immediately and forever associated with this vector alone, is not, and cannot be, realised, and we can now see why. 

Also, our results offer an insight into the reasons behind the well-known phenomenon in ML
that before the attractors are formed in required numbers and with the required separation in the phase space, and are associated with individual training vectors, the same sequence of training vectors  must be applied many times. Indeed, we show that the first training epoch initiates only the first cascade of bifurcations (see Fig.~\ref{fig:81N_Bif_Diag}) and does not create enough attractors with enough space separation. Obtaining more attractors and achieving their good spatial separation requires more and more applications of the same training vectors in order to push the weight trajectory towards other bifurcation manifolds and at the end to the location at a sufficient distance from all of them, so that the resultant NN (\ref{eqn:eq3.1}) can be capable of a robust performance.

We also revealed the structure of attractor \emph{basins} formed as a result of learning, and their link to the saddle fixed points.  This research highlights the importance for memory formation not only of the attractors, but of the less visible saddle objects.  

Our work is a rare example of studies and visualisations of truly high-dimensional and multi-parameter non-linear DSs.
The findings of this paper agree well with the mechanisms of real and spurious memory formation considered from the view point of dynamical systems with plastic self-organising vector fields (termed ``plastic dynamical systems'' for brevity)\cite{Janson_explainable_plastic_DS_SSRN25}.

\section{Discussion}
\label{sec:disc}

\paragraph{\textit{\textbf{Significance of the method for bifurcation analysis in learning NNs.}}}  Although memory formation in  learning ANNs has long been associated with bifurcations theoretically \cite{Pineda_Hopfield_NN_bifurcation_JC88,Doya_NN_bifurcations_leanring_IEEE93}, to date this connection has not been rigorously verified. 
One of the reasons 
is the dimensionality of the parameter (weight) space of NNs, which is abnormally and intractably high for the standard 
bifurcation analysis with standard analytical or numerical tools. Another reason has been the absence of a rigorous mathematical  formulation of what, in a learning NN, constitutes a bifurcation as an event occurring with the continuous change of a parameter signifying the stage of learning. 

Our proposal of a rigorous method to analyse bifurcations in a learning recurrent NN is based on formally introducing its two versions. The first version is a \emph{learning} NN (such as (\ref{eqn:eq2.1})--(\ref{eqn:eq2.4})) whose state and weights evolve with continuous time $t$, where $t$ is the stage of learning.  The second version is a \emph{non-learning} NN with fixed weights (such as (\ref{eqn:eq3.1})), which evolves with a different time $t'$, and for which all weights are formally functions of a single parameter $t$ that can serve as the bifurcation parameter signifying the stage of learning. In this setting, the bifurcation analysis of a learning NN is formally reduced to the standard bifurcation analysis of a non-learning NN depending on a single parameter $t$. 

This reduction is  in principle applicable to any recurrent NN of any dimension employing any learning algorithm.  Bifurcation analysis can be a powerful  tool to explain operations of learning NNs, even in cases where energy functions are not defined, such as Hopfield-type NNs with asymmetric weights.  

\paragraph{\textit{\textbf{Broader significance of results.}}}
Our study takes a significant step towards demystifying the ``black box'' of recurrent learning ANNs  by integrating isolated aspects of their operation as fragments of a larger unified picture. Indeed, although the Hopfield NN with Hebbian learning is different from recurrent ANNs that are more popular in ML, they share the same learning principle: guiding the weight trajectory to a point in the weight space at which the ANN has the desired configuration of attractors and basins. 

There are many ways to choose the start, the end, and the interim steps of this journey. With Hebbian learning, the weight trajectory is formed as part of self-organisation of the whole NN being a non-autonomous non-linear dissipative DS.  

 In other ML algorithms, the choice of steps is justified by a variety of
criteria often based on the loss function.  The prevailing preference is to initialise the weights at values ensuring a good number of attractors, and to avoid bifurcations during learning in order to prevent gradient explosions of the loss 
\cite{Doya_NN_bifurcations_leanring_IEEE93,Pascanu_NN_learning_bifurcations_conf13}. 
However, the appropriateness of the initial vector of weights  depends on its location relative to bifurcation manifolds, which for non-small non-linear continuous-state NNs  are generally impossible to find analytically and  are typically very difficult to find even numerically.   Also, the mechanism to associate a new stimulus with one of the existing attractors is unclear.
Moreover, such learning pre-conditions the existence of attractors  not associated with stimuli (spurious memories) during much or all of the training phase. 

It is reasonable to suggest that in other recurrent ANNs used in ML, memory formation and destruction are also caused by bifurcations.  The latter would need to be verified, and our method offers a tool to achieve this directly and rigorously. 
With this, the knowledge of the mechanisms of memory formation could
be used to prepare a network for conventional learning in a relatively predictable and controlled way. One  could guide the NN through the necessary number of attractor-creating bifurcations. If there are not enough attractors at any stage of learning, one could increase their number by forcing the weight trajectory to cross bifurcation manifolds in appropriate directions.  One could also avoid any instances of 
forgetting, including catastrophic forgetting,  by preventing unwanted crossings of bifurcation manifolds. 

Overall, by exposing the relationship between memory formation, 
forgetting, and bifurcations in a sample recurrent learning  ANN, our findings together with the existing indirect evidence suggest that similar relationships could exist in other recurrent ANNs. Assuming that this is the case, the knowledge of such relationships  
could help create learning algorithms targeting the \emph{causes} of unwanted or desirable effects. 

In practice, the development of algorithms utilising or avoiding bifurcations requires as the first step a rigorous demonstration of such bifurcations and of the possible bifurcation scenarios taking place in typical ANNs with standard learning algorithms.  Our current work represents this first step as applied to the historically and conceptually foundational learning algorithm – Hebbian one – operating within the archetypal continuous-time continuous-state Hopfield NN. Our work demonstrates that making this first step is not a trivial task both conceptually, and computationally, and requires a comprehensive and systematic assessment of the situation from different angles. The design of the learning algorithms based on the use and on the avoidance of bifurcations  is outside the scope of this paper. 
It might be worth noting that such algorithms could utilise methods from the theory of non-linear DSs for the control  of attractors \cite{Parmananda96, Claussen98, Epureanu00,Ahlborn04} and of bifurcations \cite{Abed86, Abed87, Chen00}.

\paragraph{\textit{\textbf{Future work.}}}  The hypothesis about the causes of  any type of 
forgetting cannot be verified by simply detecting crossings between the weight trajectory and bifurcation manifolds, since these manifolds cannot realistically be calculated (except for small NNs 
as in Sec.~\ref{sec:small}). Our approach to do this using a series of cross-sections of these manifolds is still technically  demanding and laborious, but is feasible with the computational tools available and could be used in the future for a more thorough verification of this hypothesis. 

Also, bifurcation analysis may help reveal the learning mechanisms in deep ANNs -- the most powerful yet poorly understood AI available to date \cite{Ali_AI_explainable_review_IF23}. 

With regard to Hebbian learning in Hopfield NNs, it could be useful to systematically explore its dependence on the values of control parameters and on the form and the steepness of the activation function $F$. Note, that $t_s$ could also be an important parameter affecting learning particularly in the context of bifurcations. Indeed, the amplitudes of weight oscillations in the course of learning depend monotonously on, and are almost directly proportional to, the length of time $t_s$ during which individual training vectors $\mathbf{I}^k$ are applied as part of the stimulus $\mathbf{I}(t)$. This is also visible in Figs. \ref{fig:fig9.1} and \ref{fig:81N_Bif_Surf}, where each ``zag’’ lasts as long as the given training vector $\mathbf{I}^k$ is applied. The given observation suggests that, whereas too small values of $t_s$ might postpone the crossing of the bifurcation manifolds to later training epochs, excessively large values of $t_s$ might lead to an undesirably large weight amplitudes making it more difficult to define the end of learning.

It would be interesting to know the memory capacity of the learning ANN considered here. However, whereas reliable estimates of memory capacity exist for discrete-time discrete-space ANNs, these results cannot be straightforwardly translated to continuous-time continuous-space networks.  Estimation of memory capacity of such networks is the considerable and generally unsolved problem in ML due to the high non-linearity of such systems and the impossibility to analytically predict the locations and the number of attractors. In the context of our research on bifurcations enabling learning, to predict the maximal number of attractors that can exist in the given NN, one would need to analytically predict the number of bifurcations along the weight trajectory, which is currently impossible to achieve. A paper Ref.~\cite{Ramsauer_Hopfield_NN_continuous_memory_capacity_arxiv21} studies memory capacity of a continuous-time continuous-state modified Hopfield NN, but its results are linked to a specific learning algorithm, which is different to Hebbian learning studied here, and are not applicable to our case.  Estimation of memory capacity linked to Hebbian learning could be the subject of future research.

\section{Supplementary Material}

Supplementary online Material in the form of ``Supplementary Note'' provides the following information. In Sec.~\ref{supp-sec:Energy_Land} energy landscape for a special form of Hopfield NN is shown. In Sec.~\ref{supp-sec:Train_Pat} data used for training the NNs are described.  In Sec.~\ref{supp-sec:alt}, for the Hopfield NN in an alternative form, basins of attraction at the end of training are shown.

\section{Author Contributions} 

AEE performed all numerical calculations reported in paper, wrote the first draft and edited the final draft. NBJ conceived this research. RN produced the first preliminary results, not illustrated in paper, under the supervision of NBJ.  AEE, NBJ and AGB  analysed the data and interpreted the results.   NBJ and AGB supervised and conceptualised the work, and co-wrote the final version of the paper. AGB coordinated the paper preparation. 

\section{Acknowledgements} AEE received funding from EPSRC (UK) to support his PhD studies under grant EP/W523987/1. 

\section{Data Availability Statement} No data was generated during this research. 


\begin{thebibliography}{75}%
\makeatletter
\providecommand \@ifxundefined [1]{%
 \@ifx{#1\undefined}
}%
\providecommand \@ifnum [1]{%
 \ifnum #1\expandafter \@firstoftwo
 \else \expandafter \@secondoftwo
 \fi
}%
\providecommand \@ifx [1]{%
 \ifx #1\expandafter \@firstoftwo
 \else \expandafter \@secondoftwo
 \fi
}%
\providecommand \natexlab [1]{#1}%
\providecommand \enquote  [1]{``#1''}%
\providecommand \bibnamefont  [1]{#1}%
\providecommand \bibfnamefont [1]{#1}%
\providecommand \citenamefont [1]{#1}%
\providecommand \href@noop [0]{\@secondoftwo}%
\providecommand \href [0]{\begingroup \@sanitize@url \@href}%
\providecommand \@href[1]{\@@startlink{#1}\@@href}%
\providecommand \@@href[1]{\endgroup#1\@@endlink}%
\providecommand \@sanitize@url [0]{\catcode `\\12\catcode `\$12\catcode
  `\&12\catcode `\#12\catcode `\^12\catcode `\_12\catcode `\%12\relax}%
\providecommand \@@startlink[1]{}%
\providecommand \@@endlink[0]{}%
\providecommand \url  [0]{\begingroup\@sanitize@url \@url }%
\providecommand \@url [1]{\endgroup\@href {#1}{\urlprefix }}%
\providecommand \urlprefix  [0]{URL }%
\providecommand \Eprint [0]{\href }%
\providecommand \doibase [0]{https://doi.org/}%
\providecommand \selectlanguage [0]{\@gobble}%
\providecommand \bibinfo  [0]{\@secondoftwo}%
\providecommand \bibfield  [0]{\@secondoftwo}%
\providecommand \translation [1]{[#1]}%
\providecommand \BibitemOpen [0]{}%
\providecommand \bibitemStop [0]{}%
\providecommand \bibitemNoStop [0]{.\EOS\space}%
\providecommand \EOS [0]{\spacefactor3000\relax}%
\providecommand \BibitemShut  [1]{\csname bibitem#1\endcsname}%
\let\auto@bib@innerbib\@empty
\bibitem [{\citenamefont {LeCun}, \citenamefont {Bengio},\ and\ \citenamefont
  {Hinton}(2015)}]{LeCun:2015aa}%
  \BibitemOpen
  \bibfield  {author} {\bibinfo {author} {\bibfnamefont {Y.}~\bibnamefont
  {LeCun}}, \bibinfo {author} {\bibfnamefont {Y.}~\bibnamefont {Bengio}},\ and\
  \bibinfo {author} {\bibfnamefont {G.}~\bibnamefont {Hinton}},\ }\bibfield
  {title} {\enquote {\bibinfo {title} {Deep learning},}\ }\href
  {https://doi.org/10.1038/nature14539} {\bibfield  {journal} {\bibinfo
  {journal} {Nature}\ }\textbf {\bibinfo {volume} {521}},\ \bibinfo {pages}
  {436--444} (\bibinfo {year} {2015})}\BibitemShut {NoStop}%
\bibitem [{\citenamefont {Alzubaidi}\ \emph
  {et~al.}(2021{\natexlab{a}})\citenamefont {Alzubaidi} \emph
  {et~al.}}]{Alzubaidi:2021aa}%
  \BibitemOpen
  \bibfield  {author} {\bibinfo {author} {\bibfnamefont {L.}~\bibnamefont
  {Alzubaidi}} \emph {et~al.},\ }\bibfield  {title} {\enquote {\bibinfo {title}
  {Review of deep learning: concepts, cnn architectures, challenges,
  applications, future directions},}\ }\href
  {https://doi.org/10.1186/s40537-021-00444-8} {\bibfield  {journal} {\bibinfo
  {journal} {J. Big Data}\ }\textbf {\bibinfo {volume} {8}},\ \bibinfo {pages}
  {53} (\bibinfo {year} {2021}{\natexlab{a}})}\BibitemShut {NoStop}%
\bibitem [{\citenamefont {Amari}(1972)}]{Amari1972}%
  \BibitemOpen
  \bibfield  {author} {\bibinfo {author} {\bibfnamefont {S.-I.}\ \bibnamefont
  {Amari}},\ }\bibfield  {title} {\enquote {\bibinfo {title} {Learning patterns
  and pattern sequences by self-organizing nets of threshold elements},}\
  }\href {https://doi.org/10.1109/T-C.1972.223477} {\bibfield  {journal}
  {\bibinfo  {journal} {IEEE Transactions on Computers}\ }\textbf {\bibinfo
  {volume} {C-21}},\ \bibinfo {pages} {1197--1206} (\bibinfo {year}
  {1972})}\BibitemShut {NoStop}%
\bibitem [{\citenamefont {Hopfield}(1982)}]{Hopfield82}%
  \BibitemOpen
  \bibfield  {author} {\bibinfo {author} {\bibfnamefont {J.}~\bibnamefont
  {Hopfield}},\ }\bibfield  {title} {\enquote {\bibinfo {title} {Neural
  networks and physical systems with emergent collective computational
  abilities},}\ }\href@noop {} {\bibfield  {journal} {\bibinfo  {journal}
  {Proc. Natl. Acad. Sci. U.S.A.}\ }\textbf {\bibinfo {volume} {79}},\ \bibinfo
  {pages} {2554--2558} (\bibinfo {year} {1982})}\BibitemShut {NoStop}%
\bibitem [{\citenamefont {Hopfield}(1984{\natexlab{a}})}]{Hopfield_NN}%
  \BibitemOpen
  \bibfield  {author} {\bibinfo {author} {\bibfnamefont {J.}~\bibnamefont
  {Hopfield}},\ }\bibfield  {title} {\enquote {\bibinfo {title} {Neurons with
  graded response have collective computational properties like those of
  two-state neurons},}\ }\href@noop {} {\bibfield  {journal} {\bibinfo
  {journal} {Proc. Natl. Acad. Sci. U.S.A.}\ }\textbf {\bibinfo {volume}
  {81}},\ \bibinfo {pages} {3088--3092} (\bibinfo {year}
  {1984}{\natexlab{a}})}\BibitemShut {NoStop}%
\bibitem [{\citenamefont {Pineda}(1988)}]{Pineda_Hopfield_NN_bifurcation_JC88}%
  \BibitemOpen
  \bibfield  {author} {\bibinfo {author} {\bibfnamefont {F.~J.}\ \bibnamefont
  {Pineda}},\ }\bibfield  {title} {\enquote {\bibinfo {title} {Dynamics and
  architecture for neural computation},}\ }\href
  {https://doi.org/https://doi.org/10.1016/0885-064X(88)90021-0} {\bibfield
  {journal} {\bibinfo  {journal} {J. Complex.}\ }\textbf {\bibinfo {volume}
  {4}},\ \bibinfo {pages} {216--245} (\bibinfo {year} {1988})}\BibitemShut
  {NoStop}%
\bibitem [{\citenamefont {Barack}\ and\ \citenamefont
  {Krakauer}(2021)}]{Barack_Two_views_on_cognitive_brain_NatRevNeurosci21}%
  \BibitemOpen
  \bibfield  {author} {\bibinfo {author} {\bibfnamefont {D.}~\bibnamefont
  {Barack}}\ and\ \bibinfo {author} {\bibfnamefont {J.}~\bibnamefont
  {Krakauer}},\ }\bibfield  {title} {\enquote {\bibinfo {title} {Two views on
  the cognitive brain},}\ }\href@noop {} {\bibfield  {journal} {\bibinfo
  {journal} {Nat. Rev. Neurosci.}\ }\textbf {\bibinfo {volume} {22}},\ \bibinfo
  {pages} {359--371} (\bibinfo {year} {2021})}\BibitemShut {NoStop}%
\bibitem [{\citenamefont {Hertz}, \citenamefont {Krogh},\ and\ \citenamefont
  {Palmer}(1991)}]{Hertz_book_91}%
  \BibitemOpen
  \bibfield  {author} {\bibinfo {author} {\bibfnamefont {J.}~\bibnamefont
  {Hertz}}, \bibinfo {author} {\bibfnamefont {A.}~\bibnamefont {Krogh}},\ and\
  \bibinfo {author} {\bibfnamefont {R.}~\bibnamefont {Palmer}},\ }\href@noop {}
  {\emph {\bibinfo {title} {Introduction to the theory of neural
  computation}}}\ (\bibinfo  {publisher} {Addison-Wesley},\ \bibinfo {year}
  {1991})\BibitemShut {NoStop}%
\bibitem [{\citenamefont {Dong}\ and\ \citenamefont
  {Hopfield}(1992)}]{Hopfield2_Dong_Connections}%
  \BibitemOpen
  \bibfield  {author} {\bibinfo {author} {\bibfnamefont {D.}~\bibnamefont
  {Dong}}\ and\ \bibinfo {author} {\bibfnamefont {J.}~\bibnamefont
  {Hopfield}},\ }\bibfield  {title} {\enquote {\bibinfo {title} {Dynamic
  properties of neural networks with adapting synapses},}\ }\href@noop {}
  {\bibfield  {journal} {\bibinfo  {journal} {Network: Computation in Neural
  Systems}\ }\textbf {\bibinfo {volume} {3}},\ \bibinfo {pages} {267--283}
  (\bibinfo {year} {1992})}\BibitemShut {NoStop}%
\bibitem [{\citenamefont {Mandic}\ and\ \citenamefont
  {Chambers}(2001)}]{Mandic_NN_recurrent_learning_algirithms_book01}%
  \BibitemOpen
  \bibfield  {author} {\bibinfo {author} {\bibfnamefont {D.~P.}\ \bibnamefont
  {Mandic}}\ and\ \bibinfo {author} {\bibfnamefont {J.~A.}\ \bibnamefont
  {Chambers}},\ }\href@noop {} {\emph {\bibinfo {title} {Recurrent Neural
  Networks for Prediction: Learning Algorithms, Architectures and Stability}}}\
  (\bibinfo  {publisher} {Wiley},\ \bibinfo {address} {Chichester, UK},\
  \bibinfo {year} {2001})\BibitemShut {NoStop}%
\bibitem [{\citenamefont {Amit}(1989)}]{Amit_89}%
  \BibitemOpen
  \bibfield  {author} {\bibinfo {author} {\bibfnamefont {D.}~\bibnamefont
  {Amit}},\ }\href@noop {} {\emph {\bibinfo {title} {Modeling brain function:
  The world of attractor neural networks}}}\ (\bibinfo  {publisher} {Cambridge
  University Press},\ \bibinfo {year} {1989})\BibitemShut {NoStop}%
\bibitem [{\citenamefont {McCloskey}\ and\ \citenamefont
  {Cohen}(1989)}]{McCloskey_catastrophic_forgetting_PLM89}%
  \BibitemOpen
  \bibfield  {author} {\bibinfo {author} {\bibfnamefont {M.}~\bibnamefont
  {McCloskey}}\ and\ \bibinfo {author} {\bibfnamefont {N.~J.}\ \bibnamefont
  {Cohen}},\ }\bibfield  {title} {\enquote {\bibinfo {title} {Catastrophic
  interference in connectionist networks: The sequential learning problem},}\
  }\href@noop {} {\bibfield  {journal} {\bibinfo  {journal} {Psychol. Learn.
  Motiv.}\ }\textbf {\bibinfo {volume} {24}},\ \bibinfo {pages} {109--165}
  (\bibinfo {year} {1989})}\BibitemShut {NoStop}%
\bibitem [{\citenamefont
  {French}(1999)}]{French_catastrophic_forgetting_in_NNs_TCS99}%
  \BibitemOpen
  \bibfield  {author} {\bibinfo {author} {\bibfnamefont {R.}~\bibnamefont
  {French}},\ }\bibfield  {title} {\enquote {\bibinfo {title} {Catastrophic
  forgetting in connectionist networks},}\ }\href@noop {} {\bibfield  {journal}
  {\bibinfo  {journal} {Trends Cogn. Sci.}\ }\textbf {\bibinfo {volume} {3}},\
  \bibinfo {pages} {128--135} (\bibinfo {year} {1999})}\BibitemShut {NoStop}%
\bibitem [{\citenamefont {Kirkpatrick}\ \emph {et~al.}(2017)\citenamefont
  {Kirkpatrick} \emph
  {et~al.}}]{Kirkpatrick_overcoming_catastrophic_forgetting_in_NNs_PNAS17}%
  \BibitemOpen
  \bibfield  {author} {\bibinfo {author} {\bibfnamefont {J.}~\bibnamefont
  {Kirkpatrick}} \emph {et~al.},\ }\bibfield  {title} {\enquote {\bibinfo
  {title} {Overcoming catastrophic forgetting in neural networks},}\
  }\href@noop {} {\bibfield  {journal} {\bibinfo  {journal} {PNAS}\ }\textbf
  {\bibinfo {volume} {114}},\ \bibinfo {pages} {3521--3526} (\bibinfo {year}
  {2017})}\BibitemShut {NoStop}%
\bibitem [{\citenamefont {Kemker}\ \emph {et~al.}(2018)\citenamefont {Kemker},
  \citenamefont {Mcclure}, \citenamefont {Abitino}, \citenamefont {Hayes},\
  and\ \citenamefont {Kanan}}]{Kemker_catastrophic_forgetting_deep_NNs_AAAI18}%
  \BibitemOpen
  \bibfield  {author} {\bibinfo {author} {\bibfnamefont {R.}~\bibnamefont
  {Kemker}}, \bibinfo {author} {\bibfnamefont {M.}~\bibnamefont {Mcclure}},
  \bibinfo {author} {\bibfnamefont {A.}~\bibnamefont {Abitino}}, \bibinfo
  {author} {\bibfnamefont {T.}~\bibnamefont {Hayes}},\ and\ \bibinfo {author}
  {\bibfnamefont {C.}~\bibnamefont {Kanan}},\ }\bibfield  {title} {\enquote
  {\bibinfo {title} {Measuring catastrophic forgetting in neural networks},}\
  }\href {https://doi.org/10.1609/aaai.v32i1.11651} {\bibfield  {journal}
  {\bibinfo  {journal} {Proceedings of the AAAI Conference on Artificial
  Intelligence}\ }\textbf {\bibinfo {volume} {32}},\ \bibinfo {pages}
  {3390--3398} (\bibinfo {year} {2018})}\BibitemShut {NoStop}%
\bibitem [{\citenamefont {Yoon}\ \emph {et~al.}(2018)\citenamefont {Yoon},
  \citenamefont {Yang}, \citenamefont {Lee},\ and\ \citenamefont
  {Hwang}}]{Yoon_Lifelong_Learning_with_Dynamically_Expandable_Netw_arxive18}%
  \BibitemOpen
  \bibfield  {author} {\bibinfo {author} {\bibfnamefont {J.}~\bibnamefont
  {Yoon}}, \bibinfo {author} {\bibfnamefont {E.}~\bibnamefont {Yang}}, \bibinfo
  {author} {\bibfnamefont {J.}~\bibnamefont {Lee}},\ and\ \bibinfo {author}
  {\bibfnamefont {S.~J.}\ \bibnamefont {Hwang}},\ }\href
  {https://arxiv.org/abs/1708.01547} {\enquote {\bibinfo {title} {Lifelong
  learning with dynamically expandable networks},}\ } (\bibinfo {year}
  {2018}),\ \Eprint {https://arxiv.org/abs/1708.01547} {arXiv:1708.01547
  [cs.LG]} \BibitemShut {NoStop}%
\bibitem [{\citenamefont
  {Parisi}(2019)}]{Parisi_lifelong_learning_NN_review_NN19}%
  \BibitemOpen
  \bibfield  {author} {\bibinfo {author} {\bibfnamefont {G.}~\bibnamefont
  {Parisi}},\ }\bibfield  {title} {\enquote {\bibinfo {title} {Continual
  lifelong learning with neural networks: A review},}\ }\href@noop {}
  {\bibfield  {journal} {\bibinfo  {journal} {Neural Netw.}\ }\textbf {\bibinfo
  {volume} {113}},\ \bibinfo {pages} {54--71} (\bibinfo {year}
  {2019})}\BibitemShut {NoStop}%
\bibitem [{\citenamefont {Ali}\ \emph {et~al.}(2023)\citenamefont {Ali},
  \citenamefont {Abuhmed}, \citenamefont {El-Sappagh}, \citenamefont
  {Muhammad}, \citenamefont {Alonso-Moral}, \citenamefont {Confalonieri},
  \citenamefont {Guidotti}, \citenamefont {{Del Ser}}, \citenamefont
  {D{\'\i}az-Rodr{\'\i}guez},\ and\ \citenamefont
  {Herrera}}]{Ali_AI_explainable_review_IF23}%
  \BibitemOpen
  \bibfield  {author} {\bibinfo {author} {\bibfnamefont {S.}~\bibnamefont
  {Ali}}, \bibinfo {author} {\bibfnamefont {T.}~\bibnamefont {Abuhmed}},
  \bibinfo {author} {\bibfnamefont {S.}~\bibnamefont {El-Sappagh}}, \bibinfo
  {author} {\bibfnamefont {K.}~\bibnamefont {Muhammad}}, \bibinfo {author}
  {\bibfnamefont {J.~M.}\ \bibnamefont {Alonso-Moral}}, \bibinfo {author}
  {\bibfnamefont {R.}~\bibnamefont {Confalonieri}}, \bibinfo {author}
  {\bibfnamefont {R.}~\bibnamefont {Guidotti}}, \bibinfo {author}
  {\bibfnamefont {J.}~\bibnamefont {{Del Ser}}}, \bibinfo {author}
  {\bibfnamefont {N.}~\bibnamefont {D{\'\i}az-Rodr{\'\i}guez}},\ and\ \bibinfo
  {author} {\bibfnamefont {F.}~\bibnamefont {Herrera}},\ }\bibfield  {title}
  {\enquote {\bibinfo {title} {Explainable artificial intelligence ({XAI}):
  What we know and what is left to attain trustworthy artificial
  intelligence},}\ }\href@noop {} {\bibfield  {journal} {\bibinfo  {journal}
  {Information Fusion}\ }\textbf {\bibinfo {volume} {99}},\ \bibinfo {pages}
  {101805} (\bibinfo {year} {2023})}\BibitemShut {NoStop}%
\bibitem [{\citenamefont {Ratcliff}(1990)}]{Ratcliff90}%
  \BibitemOpen
  \bibfield  {author} {\bibinfo {author} {\bibfnamefont {R.}~\bibnamefont
  {Ratcliff}},\ }\bibfield  {title} {\enquote {\bibinfo {title} {Connectionist
  models of recognition memory: Constraints imposed by learning and forgetting
  functions},}\ }\href@noop {} {\bibfield  {journal} {\bibinfo  {journal}
  {Psychological Review}\ }\textbf {\bibinfo {volume} {97}},\ \bibinfo {pages}
  {285--308} (\bibinfo {year} {1990})}\BibitemShut {NoStop}%
\bibitem [{\citenamefont {Robins}\ and\ \citenamefont
  {McCallum}(1998)}]{Robins98}%
  \BibitemOpen
  \bibfield  {author} {\bibinfo {author} {\bibfnamefont {A.}~\bibnamefont
  {Robins}}\ and\ \bibinfo {author} {\bibfnamefont {S.}~\bibnamefont
  {McCallum}},\ }\bibfield  {title} {\enquote {\bibinfo {title} {Catastrophic
  forgetting and the pseudorehearsal solution in {Hopfield}-type networks},}\
  }\href@noop {} {\bibfield  {journal} {\bibinfo  {journal} {Connection
  Science}\ }\textbf {\bibinfo {volume} {10}},\ \bibinfo {pages} {121--135}
  (\bibinfo {year} {1998})}\BibitemShut {NoStop}%
\bibitem [{\citenamefont {Rusu}\ \emph {et~al.}(2016)\citenamefont {Rusu},
  \citenamefont {Rabinowitz}, \citenamefont {Desjardins}, \citenamefont
  {Soyer}, \citenamefont {Kirkpatrick}, \citenamefont {Kavukcuoglu},
  \citenamefont {Pascanu},\ and\ \citenamefont
  {Hadsell}}]{Rusu_progressive_NNs_against_catastrophic_forgetting_arxiv16}%
  \BibitemOpen
  \bibfield  {author} {\bibinfo {author} {\bibfnamefont {A.}~\bibnamefont
  {Rusu}}, \bibinfo {author} {\bibfnamefont {N.}~\bibnamefont {Rabinowitz}},
  \bibinfo {author} {\bibfnamefont {G.}~\bibnamefont {Desjardins}}, \bibinfo
  {author} {\bibfnamefont {H.}~\bibnamefont {Soyer}}, \bibinfo {author}
  {\bibfnamefont {J.}~\bibnamefont {Kirkpatrick}}, \bibinfo {author}
  {\bibfnamefont {K.}~\bibnamefont {Kavukcuoglu}}, \bibinfo {author}
  {\bibfnamefont {R.}~\bibnamefont {Pascanu}},\ and\ \bibinfo {author}
  {\bibfnamefont {R.}~\bibnamefont {Hadsell}},\ }\bibfield  {title} {\enquote
  {\bibinfo {title} {Progressive neural networks},}\ }\href
  {https://doi.org/10.48550/arXiv.1606.04671} {\bibfield  {journal} {\bibinfo
  {journal} {arXiv:1606.04671}\ } (\bibinfo {year} {2016})}\BibitemShut
  {NoStop}%
\bibitem [{\citenamefont {Fayek}, \citenamefont {Cavedon},\ and\ \citenamefont
  {Wu}(2020)}]{Fayek_progressive_learning_expanding_eliminate_catastr_forgt_NN20}%
  \BibitemOpen
  \bibfield  {author} {\bibinfo {author} {\bibfnamefont {H.}~\bibnamefont
  {Fayek}}, \bibinfo {author} {\bibfnamefont {L.}~\bibnamefont {Cavedon}},\
  and\ \bibinfo {author} {\bibfnamefont {H.}~\bibnamefont {Wu}},\ }\bibfield
  {title} {\enquote {\bibinfo {title} {Progressive learning: A deep learning
  framework for continual learning},}\ }\href@noop {} {\bibfield  {journal}
  {\bibinfo  {journal} {Neural Networks}\ }\textbf {\bibinfo {volume} {128}},\
  \bibinfo {pages} {345--357} (\bibinfo {year} {2020})}\BibitemShut {NoStop}%
\bibitem [{\citenamefont {Nauta}\ \emph {et~al.}(2023)\citenamefont {Nauta},
  \citenamefont {Trienes}, \citenamefont {Pathak}, \citenamefont {Nguyen},
  \citenamefont {Peters}, \citenamefont {Schmitt}, \citenamefont
  {Schl\"{o}tterer}, \citenamefont {van Keulen},\ and\ \citenamefont
  {Seifert}}]{Nauta_NN_Explainability_Methods_Review}%
  \BibitemOpen
  \bibfield  {author} {\bibinfo {author} {\bibfnamefont {M.}~\bibnamefont
  {Nauta}}, \bibinfo {author} {\bibfnamefont {J.}~\bibnamefont {Trienes}},
  \bibinfo {author} {\bibfnamefont {S.}~\bibnamefont {Pathak}}, \bibinfo
  {author} {\bibfnamefont {E.}~\bibnamefont {Nguyen}}, \bibinfo {author}
  {\bibfnamefont {M.}~\bibnamefont {Peters}}, \bibinfo {author} {\bibfnamefont
  {Y.}~\bibnamefont {Schmitt}}, \bibinfo {author} {\bibfnamefont
  {J.}~\bibnamefont {Schl\"{o}tterer}}, \bibinfo {author} {\bibfnamefont
  {M.}~\bibnamefont {van Keulen}},\ and\ \bibinfo {author} {\bibfnamefont
  {C.}~\bibnamefont {Seifert}},\ }\bibfield  {title} {\enquote {\bibinfo
  {title} {From anecdotal evidence to quantitative evaluation methods: A
  systematic review on evaluating explainable {AI}},}\ }\href@noop {}
  {\bibfield  {journal} {\bibinfo  {journal} {ACM Comput. Surv.}\ } (\bibinfo
  {year} {2023})}\BibitemShut {NoStop}%
\bibitem [{\citenamefont {Saleem}\ \emph {et~al.}(2022)\citenamefont {Saleem},
  \citenamefont {Yuan}, \citenamefont {Kurugollu}, \citenamefont {Anjum},\ and\
  \citenamefont {Liu}}]{Saleem_DNN_Explainability_Overview}%
  \BibitemOpen
  \bibfield  {author} {\bibinfo {author} {\bibfnamefont {R.}~\bibnamefont
  {Saleem}}, \bibinfo {author} {\bibfnamefont {B.}~\bibnamefont {Yuan}},
  \bibinfo {author} {\bibfnamefont {F.}~\bibnamefont {Kurugollu}}, \bibinfo
  {author} {\bibfnamefont {A.}~\bibnamefont {Anjum}},\ and\ \bibinfo {author}
  {\bibfnamefont {L.}~\bibnamefont {Liu}},\ }\bibfield  {title} {\enquote
  {\bibinfo {title} {Explaining deep neural networks: A survey on the global
  interpretation methods},}\ }\href
  {https://doi.org/https://doi.org/10.1016/j.neucom.2022.09.129} {\bibfield
  {journal} {\bibinfo  {journal} {Neurocomputing}\ }\textbf {\bibinfo {volume}
  {513}},\ \bibinfo {pages} {165--180} (\bibinfo {year} {2022})}\BibitemShut
  {NoStop}%
\bibitem [{\citenamefont {Doya}(1993)}]{Doya_NN_bifurcations_leanring_IEEE93}%
  \BibitemOpen
  \bibfield  {author} {\bibinfo {author} {\bibfnamefont {K.}~\bibnamefont
  {Doya}},\ }\bibfield  {title} {\enquote {\bibinfo {title} {Bifurcations of
  recurrent neural networks in gradient descent learning},}\ }\href@noop {}
  {\bibfield  {journal} {\bibinfo  {journal} {IEEE Transactions on neural
  networks}\ }\textbf {\bibinfo {volume} {1}},\ \bibinfo {pages} {218}
  (\bibinfo {year} {1993})}\BibitemShut {NoStop}%
\bibitem [{\citenamefont {Pascanu}, \citenamefont {Mikolov},\ and\
  \citenamefont {Bengio}(2013)}]{Pascanu_NN_learning_bifurcations_conf13}%
  \BibitemOpen
  \bibfield  {author} {\bibinfo {author} {\bibfnamefont {R.}~\bibnamefont
  {Pascanu}}, \bibinfo {author} {\bibfnamefont {T.}~\bibnamefont {Mikolov}},\
  and\ \bibinfo {author} {\bibfnamefont {Y.}~\bibnamefont {Bengio}},\
  }\bibfield  {title} {\enquote {\bibinfo {title} {On the difficulty of
  training recurrent neural networks},}\ }in\ \href
  {https://proceedings.mlr.press/v28/pascanu13.html} {\emph {\bibinfo
  {booktitle} {Proceedings of the 30th International Conference on Machine
  Learning}}},\ \bibinfo {series} {Proceedings of Machine Learning Research},
  Vol.\ \bibinfo {volume} {28(3)},\ \bibinfo {editor} {edited by\ \bibinfo
  {editor} {\bibfnamefont {S.}~\bibnamefont {Dasgupta}}\ and\ \bibinfo {editor}
  {\bibfnamefont {D.}~\bibnamefont {McAllester}}}\ (\bibinfo  {publisher}
  {PMLR},\ \bibinfo {address} {Atlanta, Georgia, USA},\ \bibinfo {year}
  {2013})\ pp.\ \bibinfo {pages} {1310--1318}\BibitemShut {NoStop}%
\bibitem [{\citenamefont {Ribeiro}\ \emph {et~al.}(2020)\citenamefont
  {Ribeiro}, \citenamefont {Tiels}, \citenamefont {Aguirre},\ and\
  \citenamefont
  {Sch\"on}}]{Ribeiro_NN_bifurcations_exploding_gradients_conf20}%
  \BibitemOpen
  \bibfield  {author} {\bibinfo {author} {\bibfnamefont {A.}~\bibnamefont
  {Ribeiro}}, \bibinfo {author} {\bibfnamefont {K.}~\bibnamefont {Tiels}},
  \bibinfo {author} {\bibfnamefont {L.}~\bibnamefont {Aguirre}},\ and\ \bibinfo
  {author} {\bibfnamefont {T.}~\bibnamefont {Sch\"on}},\ }\bibfield  {title}
  {\enquote {\bibinfo {title} {Beyond exploding and vanishing gradients:
  analysing rnn training using attractors and smoothness},}\ }in\ \href
  {https://proceedings.mlr.press/v108/ribeiro20a.html} {\emph {\bibinfo
  {booktitle} {Proceedings of the Twenty Third International Conference on
  Artificial Intelligence and Statistics}}},\ \bibinfo {series} {Proceedings of
  Machine Learning Research}, Vol.\ \bibinfo {volume} {108},\ \bibinfo {editor}
  {edited by\ \bibinfo {editor} {\bibfnamefont {S.}~\bibnamefont {Chiappa}}\
  and\ \bibinfo {editor} {\bibfnamefont {R.}~\bibnamefont {Calandra}}}\
  (\bibinfo  {publisher} {PMLR},\ \bibinfo {year} {2020})\ pp.\ \bibinfo
  {pages} {2370--2380}\BibitemShut {NoStop}%
\bibitem [{\citenamefont {Haputhanthri}\ \emph {et~al.}(2024)\citenamefont
  {Haputhanthri}, \citenamefont {Storan}, \citenamefont {Jiang}, \citenamefont
  {Shai}, \citenamefont {Orhun}, \citenamefont {Schnitzer}, \citenamefont
  {Dinc},\ and\ \citenamefont
  {Tanaka}}]{Haputhanthri_Why_NN_learn_bifurcations_24}%
  \BibitemOpen
  \bibfield  {author} {\bibinfo {author} {\bibfnamefont {U.}~\bibnamefont
  {Haputhanthri}}, \bibinfo {author} {\bibfnamefont {L.}~\bibnamefont
  {Storan}}, \bibinfo {author} {\bibfnamefont {Y.}~\bibnamefont {Jiang}},
  \bibinfo {author} {\bibfnamefont {A.}~\bibnamefont {Shai}}, \bibinfo {author}
  {\bibfnamefont {H.~A.}\ \bibnamefont {Orhun}}, \bibinfo {author}
  {\bibfnamefont {M.}~\bibnamefont {Schnitzer}}, \bibinfo {author}
  {\bibfnamefont {F.}~\bibnamefont {Dinc}},\ and\ \bibinfo {author}
  {\bibfnamefont {H.}~\bibnamefont {Tanaka}},\ }\bibfield  {title} {\enquote
  {\bibinfo {title} {Why do recurrent neural networks suddenly learn?
  {B}ifurcation mechanisms in neuro- inspired short-term memory tasks},}\ }in\
  \href {https://openreview.net/pdf?id=njmXdqzHJq} {\emph {\bibinfo {booktitle}
  {ICML 2024 Workshop on Mechanistic Interpretability}}}\ (\bibinfo {year}
  {2024})\BibitemShut {NoStop}%
\bibitem [{\citenamefont {Hebb}(1949)}]{Hebb49}%
  \BibitemOpen
  \bibfield  {author} {\bibinfo {author} {\bibfnamefont {D.}~\bibnamefont
  {Hebb}},\ }\href@noop {} {\emph {\bibinfo {title} {The organization of
  behavior; a neuropsychological theory}}}\ (\bibinfo  {publisher} {Wiley},\
  \bibinfo {address} {Oxford, England},\ \bibinfo {year} {1949})\BibitemShut
  {NoStop}%
\bibitem [{\citenamefont {Gerstner}\ and\ \citenamefont
  {Kistler}(2002)}]{Gerstner_Hebbian_learning_BC02}%
  \BibitemOpen
  \bibfield  {author} {\bibinfo {author} {\bibfnamefont {W.}~\bibnamefont
  {Gerstner}}\ and\ \bibinfo {author} {\bibfnamefont {W.}~\bibnamefont
  {Kistler}},\ }\bibfield  {title} {\enquote {\bibinfo {title} {Mathematical
  formulations of {H}ebbian learning},}\ }\href@noop {} {\bibfield  {journal}
  {\bibinfo  {journal} {Biol. Cybern.}\ }\textbf {\bibinfo {volume} {87}},\
  \bibinfo {pages} {404--415} (\bibinfo {year} {2002})}\BibitemShut {NoStop}%
\bibitem [{\citenamefont {Toneva}\ \emph {et~al.}(2019)\citenamefont {Toneva},
  \citenamefont {Sordoni}, \citenamefont {des Combes}, \citenamefont
  {Trischler}, \citenamefont {Bengio},\ and\ \citenamefont
  {Gordon}}]{Toneva_NN_forgetting_single_task_arXiv19}%
  \BibitemOpen
  \bibfield  {author} {\bibinfo {author} {\bibfnamefont {M.}~\bibnamefont
  {Toneva}}, \bibinfo {author} {\bibfnamefont {A.}~\bibnamefont {Sordoni}},
  \bibinfo {author} {\bibfnamefont {R.~T.}\ \bibnamefont {des Combes}},
  \bibinfo {author} {\bibfnamefont {A.}~\bibnamefont {Trischler}}, \bibinfo
  {author} {\bibfnamefont {Y.}~\bibnamefont {Bengio}},\ and\ \bibinfo {author}
  {\bibfnamefont {G.~J.}\ \bibnamefont {Gordon}},\ }\bibfield  {title}
  {\enquote {\bibinfo {title} {An empirical study of example forgetting during
  deep neural network learning},}\ }in\ \href
  {https://openreview.net/forum?id=BJlxm30cKm} {\emph {\bibinfo {booktitle}
  {International Conference on Learning Representations}}}\ (\bibinfo {year}
  {2019})\BibitemShut {NoStop}%
\bibitem [{\citenamefont {Eisenmann}\ \emph {et~al.}(2023)\citenamefont
  {Eisenmann}, \citenamefont {Monfared}, \citenamefont {G\"{o}ring},\ and\
  \citenamefont
  {Durstewitz}}]{Eisenmann_NN_discrete_2-6_neuron_bifurcation_NeurIPS23}%
  \BibitemOpen
  \bibfield  {author} {\bibinfo {author} {\bibfnamefont {L.}~\bibnamefont
  {Eisenmann}}, \bibinfo {author} {\bibfnamefont {Z.}~\bibnamefont {Monfared}},
  \bibinfo {author} {\bibfnamefont {N.}~\bibnamefont {G\"{o}ring}},\ and\
  \bibinfo {author} {\bibfnamefont {D.}~\bibnamefont {Durstewitz}},\ }\bibfield
   {title} {\enquote {\bibinfo {title} {Bifurcations and loss jumps in {RNN}
  training},}\ }in\ \href@noop {} {\emph {\bibinfo {booktitle} {Proceedings of
  the 37th International Conference on Neural Information Processing
  Systems}}},\ \bibinfo {series and number} {NIPS '23}\ (\bibinfo  {publisher}
  {Curran Associates Inc.},\ \bibinfo {address} {Red Hook, NY, USA},\ \bibinfo
  {year} {2023})\BibitemShut {NoStop}%
\bibitem [{\citenamefont {Yu}\ \emph {et~al.}(2020)\citenamefont {Yu},
  \citenamefont {Abdulghani}, \citenamefont {Zahid}, \citenamefont {Heidari},
  \citenamefont {Imran},\ and\ \citenamefont {Abbasi}}]{Yu2020}%
  \BibitemOpen
  \bibfield  {author} {\bibinfo {author} {\bibfnamefont {Z.}~\bibnamefont
  {Yu}}, \bibinfo {author} {\bibfnamefont {A.}~\bibnamefont {Abdulghani}},
  \bibinfo {author} {\bibfnamefont {A.}~\bibnamefont {Zahid}}, \bibinfo
  {author} {\bibfnamefont {H.}~\bibnamefont {Heidari}}, \bibinfo {author}
  {\bibfnamefont {M.}~\bibnamefont {Imran}},\ and\ \bibinfo {author}
  {\bibfnamefont {O.}~\bibnamefont {Abbasi}},\ }\bibfield  {title} {\enquote
  {\bibinfo {title} {An overview of neuromorphic computing for artificial
  intelligence enabled hardware-based hopfield neural network},}\ }\href
  {https://doi.org/10.1109/ACCESS.2020.2985839} {\bibfield  {journal} {\bibinfo
   {journal} {IEEE Access}\ }\textbf {\bibinfo {volume} {8}},\ \bibinfo {pages}
  {67085--67099} (\bibinfo {year} {2020})}\BibitemShut {NoStop}%
\bibitem [{\citenamefont {Krotov}(2023)}]{Krotov:2023aa}%
  \BibitemOpen
  \bibfield  {author} {\bibinfo {author} {\bibfnamefont {D.}~\bibnamefont
  {Krotov}},\ }\bibfield  {title} {\enquote {\bibinfo {title} {A new frontier
  for {H}opfield networks},}\ }\href@noop {} {\bibfield  {journal} {\bibinfo
  {journal} {Nat. Rev. Phys.}\ }\textbf {\bibinfo {volume} {5}},\ \bibinfo
  {pages} {366--367} (\bibinfo {year} {2023})}\BibitemShut {NoStop}%
\bibitem [{\citenamefont
  {K\"{o}berle}(1989)}]{Koberle_NN_content-addressable_memory_CPC89}%
  \BibitemOpen
  \bibfield  {author} {\bibinfo {author} {\bibfnamefont {R.}~\bibnamefont
  {K\"{o}berle}},\ }\bibfield  {title} {\enquote {\bibinfo {title} {Neural
  networks as content addressable memories and learning machines},}\
  }\href@noop {} {\bibfield  {journal} {\bibinfo  {journal} {Comput. Phys.
  Commun.}\ }\textbf {\bibinfo {volume} {56}},\ \bibinfo {pages} {43--50}
  (\bibinfo {year} {1989})}\BibitemShut {NoStop}%
\bibitem [{\citenamefont {Janson}\ and\ \citenamefont
  {Marsden}(2017)}]{Janson:2017aa}%
  \BibitemOpen
  \bibfield  {author} {\bibinfo {author} {\bibfnamefont {N.}~\bibnamefont
  {Janson}}\ and\ \bibinfo {author} {\bibfnamefont {C.}~\bibnamefont
  {Marsden}},\ }\bibfield  {title} {\enquote {\bibinfo {title} {Dynamical
  system with plastic self-organized velocity field as an alternative
  conceptual model of a cognitive system},}\ }\href@noop {} {\bibfield
  {journal} {\bibinfo  {journal} {Sci. Rep.}\ }\textbf {\bibinfo {volume}
  {7}},\ \bibinfo {pages} {17007} (\bibinfo {year} {2017})}\BibitemShut
  {NoStop}%
\bibitem [{\citenamefont {Janson}\ and\ \citenamefont
  {Kloeden}(2021)}]{Janson_NonAutonomous_Attractors}%
  \BibitemOpen
  \bibfield  {author} {\bibinfo {author} {\bibfnamefont {N.}~\bibnamefont
  {Janson}}\ and\ \bibinfo {author} {\bibfnamefont {P.}~\bibnamefont
  {Kloeden}},\ }\bibfield  {title} {\enquote {\bibinfo {title} {Robustness of a
  dynamical systems model with a plastic self-organising vector field to noisy
  input signals},}\ }\href@noop {} {\bibfield  {journal} {\bibinfo  {journal}
  {Eur. Phys. J. Plus}\ }\textbf {\bibinfo {volume} {136}},\ \bibinfo {pages}
  {720} (\bibinfo {year} {2021})}\BibitemShut {NoStop}%
\bibitem [{\citenamefont {Alzubaidi}\ \emph
  {et~al.}(2021{\natexlab{b}})\citenamefont {Alzubaidi}, \citenamefont {Zhang},
  \citenamefont {Humaidi} \emph {et~al.}}]{Alzubaidi_training_epoch_NN_JBD21}%
  \BibitemOpen
  \bibfield  {author} {\bibinfo {author} {\bibfnamefont {L.}~\bibnamefont
  {Alzubaidi}}, \bibinfo {author} {\bibfnamefont {J.}~\bibnamefont {Zhang}},
  \bibinfo {author} {\bibfnamefont {A.}~\bibnamefont {Humaidi}}, \emph
  {et~al.},\ }\bibfield  {title} {\enquote {\bibinfo {title} {Review of deep
  learning: concepts, {CNN} architectures, challenges, applications, future
  directions},}\ }\href
  {https://doi.org/https://doi.org/10.1186/s40537-021-00444-8} {\bibfield
  {journal} {\bibinfo  {journal} {J. Big Data}\ }\textbf {\bibinfo {volume}
  {8}},\ \bibinfo {pages} {53} (\bibinfo {year}
  {2021}{\natexlab{b}})}\BibitemShut {NoStop}%
\bibitem [{\citenamefont {Kloeden}\ and\ \citenamefont
  {Yang}(2020)}]{Kloeden_b2020}%
  \BibitemOpen
  \bibfield  {author} {\bibinfo {author} {\bibfnamefont {P.}~\bibnamefont
  {Kloeden}}\ and\ \bibinfo {author} {\bibfnamefont {M.}~\bibnamefont {Yang}},\
  }\href {https://www.worldscientific.com/doi/abs/10.1142/12053} {\emph
  {\bibinfo {title} {An Introduction to Nonautonomous Dynamical Systems and
  their Attractors}}}\ (\bibinfo  {publisher} {World Scientific},\ \bibinfo
  {year} {2020})\ \Eprint
  {https://arxiv.org/abs/https://www.worldscientific.com/doi/pdf/10.1142/12053}
  {https://www.worldscientific.com/doi/pdf/10.1142/12053} \BibitemShut
  {NoStop}%
\bibitem [{\citenamefont {Janson}, \citenamefont {Essex},\ and\ \citenamefont
  {Balanov}(2024)}]{Janson_explainable_plastic_DS_SSRN25}%
  \BibitemOpen
  \bibfield  {author} {\bibinfo {author} {\bibfnamefont {N.}~\bibnamefont
  {Janson}}, \bibinfo {author} {\bibfnamefont {A.}~\bibnamefont {Essex}},\ and\
  \bibinfo {author} {\bibfnamefont {A.}~\bibnamefont {Balanov}},\ }\href
  {http://dx.doi.org/10.2139/ssrn.5123206} {\enquote {\bibinfo {title}
  {Designing explainable cognitive systems and explaining neural networks with
  plastic dynamical systems},}\ }\bibinfo {howpublished}
  {\url{http://dx.doi.org/10.2139/ssrn.5123206}} (\bibinfo {year} {2024}),\
  \Eprint {https://arxiv.org/abs/ssrn.5123206} {SSRN:ssrn.5123206} \BibitemShut
  {NoStop}%
\bibitem [{\citenamefont {Saxena}, \citenamefont {Shobe},\ and\ \citenamefont
  {McNaughton}(2022)}]{Saxena_NN_deep_fighting_catastrophic_forgetting_PNAS22}%
  \BibitemOpen
  \bibfield  {author} {\bibinfo {author} {\bibfnamefont {R.}~\bibnamefont
  {Saxena}}, \bibinfo {author} {\bibfnamefont {J.~L.}\ \bibnamefont {Shobe}},\
  and\ \bibinfo {author} {\bibfnamefont {B.~L.}\ \bibnamefont {McNaughton}},\
  }\bibfield  {title} {\enquote {\bibinfo {title} {Learning in deep neural
  networks and brains with similarity-weighted interleaved learning},}\ }\href
  {https://doi.org/10.1073/pnas.2115229119} {\bibfield  {journal} {\bibinfo
  {journal} {Proc. Natl. Acad. Sci. U.S.A.}\ }\textbf {\bibinfo {volume}
  {119}},\ \bibinfo {pages} {e2115229119} (\bibinfo {year} {2022})},\ \Eprint
  {https://arxiv.org/abs/https://www.pnas.org/doi/pdf/10.1073/pnas.2115229119}
  {https://www.pnas.org/doi/pdf/10.1073/pnas.2115229119} \BibitemShut {NoStop}%
\bibitem [{Note1()}]{Note1}%
  \BibitemOpen
  \bibinfo {note} {Global bifurcations are irrelevant to this
  study.}\BibitemShut {Stop}%
\bibitem [{\citenamefont {Kuznetsov}(1998)}]{Kuznetsov_bif_theory_98}%
  \BibitemOpen
  \bibfield  {author} {\bibinfo {author} {\bibfnamefont {I.}~\bibnamefont
  {Kuznetsov}},\ }\href@noop {} {\emph {\bibinfo {title} {Elements of applied
  bifurcation theory}}},\ Vol.\ \bibinfo {volume} {112}\ (\bibinfo  {publisher}
  {Springer Verlag},\ \bibinfo {year} {1998})\BibitemShut {NoStop}%
\bibitem [{\citenamefont {Allgower}\ and\ \citenamefont
  {Georg}(1990)}]{Allgower_continuation_book90}%
  \BibitemOpen
  \bibfield  {author} {\bibinfo {author} {\bibfnamefont {E.}~\bibnamefont
  {Allgower}}\ and\ \bibinfo {author} {\bibfnamefont {K.}~\bibnamefont
  {Georg}},\ }\href {https://doi.org/doi.org/10.1007/978-3-642-61257-2_1}
  {\emph {\bibinfo {title} {Numerical Continuation Methods}}},\ Vol.~\bibinfo
  {volume} {13}\ (\bibinfo  {publisher} {Springer (Berlin, Heidelberg)},\
  \bibinfo {year} {1990})\BibitemShut {NoStop}%
\bibitem [{\citenamefont {Doedel}\ \emph {et~al.}(1997)\citenamefont {Doedel},
  \citenamefont {Champneys}, \citenamefont {Fairgrieve}, \citenamefont
  {Kuznetsov}, \citenamefont {Sandstede},\ and\ \citenamefont
  {Wang}}]{Doedel_AUTO_97}%
  \BibitemOpen
  \bibfield  {author} {\bibinfo {author} {\bibfnamefont {E.}~\bibnamefont
  {Doedel}}, \bibinfo {author} {\bibfnamefont {A.}~\bibnamefont {Champneys}},
  \bibinfo {author} {\bibfnamefont {T.}~\bibnamefont {Fairgrieve}}, \bibinfo
  {author} {\bibfnamefont {Y.}~\bibnamefont {Kuznetsov}}, \bibinfo {author}
  {\bibfnamefont {B.}~\bibnamefont {Sandstede}},\ and\ \bibinfo {author}
  {\bibfnamefont {X.}~\bibnamefont {Wang}},\ }\href@noop {} {\enquote {\bibinfo
  {title} {{AUTO97}: Continuation and bifurcation software for ordinary
  differential equations (with {H}om{C}ont)},}\ }\bibinfo {type} {Tech. Rep.}\
  (\bibinfo  {institution} {Technical Report, Concordia University},\ \bibinfo
  {year} {1997})\BibitemShut {NoStop}%
\bibitem [{\citenamefont {Ermentrout}(2002)}]{Ermentrout_XPPAUT_guide_book02}%
  \BibitemOpen
  \bibfield  {author} {\bibinfo {author} {\bibfnamefont {B.}~\bibnamefont
  {Ermentrout}},\ }\bibfield  {title} {\enquote {\bibinfo {title} {Simulating,
  analyzing, and animating dynamical systems: A guide to xppaut for researchers
  and students},}\ }\href {https://doi.org/10.1137/1.9780898718195} {\bibfield
  {journal} {\bibinfo  {journal} {SERBIULA (sistema Librum 2.0)}\ } (\bibinfo
  {year} {2002}),\ 10.1137/1.9780898718195}\BibitemShut {NoStop}%
\bibitem [{\citenamefont {Das}\ and\ \citenamefont
  {Schieve}(1995)}]{Das_Hopfield_NN_4-neuron_bifurcation_PhD95}%
  \BibitemOpen
  \bibfield  {author} {\bibinfo {author} {\bibfnamefont {P.}~\bibnamefont
  {Das}}\ and\ \bibinfo {author} {\bibfnamefont {W.}~\bibnamefont {Schieve}},\
  }\bibfield  {title} {\enquote {\bibinfo {title} {A bifurcation analysis of
  the four dimensional generalized hopfield neural network},}\ }\href@noop {}
  {\bibfield  {journal} {\bibinfo  {journal} {Phys. D: Nonlinear Phenom.}\
  }\textbf {\bibinfo {volume} {88}},\ \bibinfo {pages} {14--28} (\bibinfo
  {year} {1995})}\BibitemShut {NoStop}%
\bibitem [{\citenamefont
  {Beer}(1995)}]{Beer_small_NNs_bifurcations_also_some_large_NNs_AB95}%
  \BibitemOpen
  \bibfield  {author} {\bibinfo {author} {\bibfnamefont {R.}~\bibnamefont
  {Beer}},\ }\bibfield  {title} {\enquote {\bibinfo {title} {On the dynamics of
  small continuous-time recurrent neural networks},}\ }\href@noop {} {\bibfield
   {journal} {\bibinfo  {journal} {Adaptive Behavior}\ }\textbf {\bibinfo
  {volume} {3}},\ \bibinfo {pages} {469--509} (\bibinfo {year}
  {1995})}\BibitemShut {NoStop}%
\bibitem [{\citenamefont {Haschke}\ and\ \citenamefont
  {Steil}(2005)}]{Haschke_NN_discrete_2_3_neurons_non-learning_bifurcations_NC05}%
  \BibitemOpen
  \bibfield  {author} {\bibinfo {author} {\bibfnamefont {R.}~\bibnamefont
  {Haschke}}\ and\ \bibinfo {author} {\bibfnamefont {J.~J.}\ \bibnamefont
  {Steil}},\ }\bibfield  {title} {\enquote {\bibinfo {title} {Input space
  bifurcation manifolds of recurrent neural networks},}\ }\href
  {https://doi.org/https://doi.org/10.1016/j.neucom.2004.11.030} {\bibfield
  {journal} {\bibinfo  {journal} {Neurocomputing}\ }\textbf {\bibinfo {volume}
  {64}},\ \bibinfo {pages} {25--38} (\bibinfo {year} {2005})},\ \bibinfo {note}
  {trends in Neurocomputing: 12th European Symposium on Artificial Neural
  Networks 2004}\BibitemShut {NoStop}%
\bibitem [{\citenamefont {Huang}\ and\ \citenamefont
  {Huang}(2008)}]{Huang_Hopfield_NN_3_neuron_bif_chaos_AMC08}%
  \BibitemOpen
  \bibfield  {author} {\bibinfo {author} {\bibfnamefont {W.-Z.}\ \bibnamefont
  {Huang}}\ and\ \bibinfo {author} {\bibfnamefont {Y.}~\bibnamefont {Huang}},\
  }\bibfield  {title} {\enquote {\bibinfo {title} {Chaos of a new class of
  hopfield neural networks},}\ }\href@noop {} {\bibfield  {journal} {\bibinfo
  {journal} {Applied Mathematics and Computation}\ }\textbf {\bibinfo {volume}
  {206}},\ \bibinfo {pages} {1--11} (\bibinfo {year} {2008})}\BibitemShut
  {NoStop}%
\bibitem [{\citenamefont {Cervantes-Ojeda}, \citenamefont {Gómez-Fuentes},\
  and\ \citenamefont
  {Bernal-Jaquez}(2017)}]{Cervantes-Ojeda_NNs_N_2_bifurcation_diagrams_NC17}%
  \BibitemOpen
  \bibfield  {author} {\bibinfo {author} {\bibfnamefont {J.}~\bibnamefont
  {Cervantes-Ojeda}}, \bibinfo {author} {\bibfnamefont {M.}~\bibnamefont
  {Gómez-Fuentes}},\ and\ \bibinfo {author} {\bibfnamefont {R.}~\bibnamefont
  {Bernal-Jaquez}},\ }\bibfield  {title} {\enquote {\bibinfo {title} {Empirical
  analysis of bifurcations in the full weights space of a two-neuron
  {DTRNN}},}\ }\href
  {https://doi.org/https://doi.org/10.1016/j.neucom.2017.01.027} {\bibfield
  {journal} {\bibinfo  {journal} {Neurocomputing}\ }\textbf {\bibinfo {volume}
  {237}},\ \bibinfo {pages} {362--374} (\bibinfo {year} {2017})}\BibitemShut
  {NoStop}%
\bibitem [{\citenamefont {Njitacke}\ \emph {et~al.}(2019)\citenamefont
  {Njitacke}, \citenamefont {Kengne}, \citenamefont {Fozin}, \citenamefont
  {Leutcha},\ and\ \citenamefont
  {Fotsin}}]{Njitacke_Hopfield_4_neurons_bifurcations_IJDC19}%
  \BibitemOpen
  \bibfield  {author} {\bibinfo {author} {\bibfnamefont {Z.}~\bibnamefont
  {Njitacke}}, \bibinfo {author} {\bibfnamefont {J.}~\bibnamefont {Kengne}},
  \bibinfo {author} {\bibfnamefont {T.}~\bibnamefont {Fozin}}, \bibinfo
  {author} {\bibfnamefont {B.}~\bibnamefont {Leutcha}},\ and\ \bibinfo {author}
  {\bibfnamefont {H.}~\bibnamefont {Fotsin}},\ }\bibfield  {title} {\enquote
  {\bibinfo {title} {Dynamical analysis of a novel 4-neurons based hopfield
  neural network: emergences of antimonotonicity and coexistence of multiple
  stable states},}\ }\href@noop {} {\bibfield  {journal} {\bibinfo  {journal}
  {Int. J. Dyn. Contr.}\ }\textbf {\bibinfo {volume} {7}},\ \bibinfo {pages}
  {823--841} (\bibinfo {year} {2019})}\BibitemShut {NoStop}%
\bibitem [{\citenamefont {Hu}\ and\ \citenamefont
  {Wang}(2023)}]{Hu_Hopfield_NN_3-neuron_forced_multi-scroll_attractors_enctyption_MTA24}%
  \BibitemOpen
  \bibfield  {author} {\bibinfo {author} {\bibfnamefont {Z.}~\bibnamefont
  {Hu}}\ and\ \bibinfo {author} {\bibfnamefont {C.}~\bibnamefont {Wang}},\
  }\bibfield  {title} {\enquote {\bibinfo {title} {Hopfield neural network with
  multi-scroll attractors and application in image encryption},}\ }\href@noop
  {} {\bibfield  {journal} {\bibinfo  {journal} {Multimedia Tools and
  Applications}\ }\textbf {\bibinfo {volume} {83}},\ \bibinfo {pages} {97--117}
  (\bibinfo {year} {2023})}\BibitemShut {NoStop}%
\bibitem [{\citenamefont {Sharpe}\ and\ \citenamefont
  {Thorne}(1982)}]{Sharpe_Arc_Length_Parameterisation}%
  \BibitemOpen
  \bibfield  {author} {\bibinfo {author} {\bibfnamefont {R.}~\bibnamefont
  {Sharpe}}\ and\ \bibinfo {author} {\bibfnamefont {R.}~\bibnamefont
  {Thorne}},\ }\bibfield  {title} {\enquote {\bibinfo {title} {Numerical method
  for extracting an arc length parameterization from parametric curves},}\
  }\href {https://www.sciencedirect.com/science/article/pii/0010448582901713}
  {\bibfield  {journal} {\bibinfo  {journal} {Computer-Aided Design}\ }\textbf
  {\bibinfo {volume} {14}},\ \bibinfo {pages} {79--81} (\bibinfo {year}
  {1982})}\BibitemShut {NoStop}%
\bibitem [{\citenamefont
  {Hopfield}(1984{\natexlab{b}})}]{Hopfield_neurons_graded_response_PNAS84}%
  \BibitemOpen
  \bibfield  {author} {\bibinfo {author} {\bibfnamefont {J.}~\bibnamefont
  {Hopfield}},\ }\bibfield  {title} {\enquote {\bibinfo {title} {Neurons with
  graded response have collective computational properties like those of
  two-state neurons},}\ }\href@noop {} {\bibfield  {journal} {\bibinfo
  {journal} {PNAS}\ }\textbf {\bibinfo {volume} {81}},\ \bibinfo {pages}
  {3088--3092} (\bibinfo {year} {1984}{\natexlab{b}})}\BibitemShut {NoStop}%
\bibitem [{\citenamefont {Smith}\ and\ \citenamefont
  {Liles}(1982)}]{Smith_exposure-time_memory_formation_visual_PHFS82}%
  \BibitemOpen
  \bibfield  {author} {\bibinfo {author} {\bibfnamefont {L.}~\bibnamefont
  {Smith}}\ and\ \bibinfo {author} {\bibfnamefont {D.}~\bibnamefont {Liles}},\
  }\bibfield  {title} {\enquote {\bibinfo {title} {The effects of exposure time
  and retention time on location memory in visual information processing},}\
  }\href@noop {} {\bibfield  {journal} {\bibinfo  {journal} {Proc. of the Human
  Factors Soc. Annual Meeting}\ }\textbf {\bibinfo {volume} {26}},\ \bibinfo
  {pages} {812--815} (\bibinfo {year} {1982})}\BibitemShut {NoStop}%
\bibitem [{\citenamefont
  {Doya}(1992)}]{Doya_Bifurcations_in_the_learning_of_recurrent_neural_networks_IEEE92}%
  \BibitemOpen
  \bibfield  {author} {\bibinfo {author} {\bibfnamefont {K.}~\bibnamefont
  {Doya}},\ }\bibfield  {title} {\enquote {\bibinfo {title} {Bifurcations in
  the learning of recurrent neural networks},}\ }in\ \href
  {https://doi.org/10.1109/ISCAS.1992.230622} {\emph {\bibinfo {booktitle}
  {[Proceedings] 1992 IEEE International Symposium on Circuits and Systems}}},\
  Vol.~\bibinfo {volume} {6}\ (\bibinfo {year} {1992})\ pp.\ \bibinfo {pages}
  {2777--2780}\BibitemShut {NoStop}%
\bibitem [{\citenamefont {Raghavan}\ \emph {et~al.}(2024)\citenamefont
  {Raghavan}, \citenamefont {Tharwat}, \citenamefont {Hari}, \citenamefont
  {Satani}, \citenamefont {Liu},\ and\ \citenamefont
  {Thomson}}]{Raghavan_NN_learning_as_path_in_weights_space_NatMI24}%
  \BibitemOpen
  \bibfield  {author} {\bibinfo {author} {\bibfnamefont {G.}~\bibnamefont
  {Raghavan}}, \bibinfo {author} {\bibfnamefont {B.}~\bibnamefont {Tharwat}},
  \bibinfo {author} {\bibfnamefont {S.~N.}\ \bibnamefont {Hari}}, \bibinfo
  {author} {\bibfnamefont {D.}~\bibnamefont {Satani}}, \bibinfo {author}
  {\bibfnamefont {R.}~\bibnamefont {Liu}},\ and\ \bibinfo {author}
  {\bibfnamefont {M.}~\bibnamefont {Thomson}},\ }\bibfield  {title} {\enquote
  {\bibinfo {title} {Engineering flexible machine learning systems by
  traversing functionally invariant paths},}\ }\href@noop {} {\bibfield
  {journal} {\bibinfo  {journal} {Nat. Mac. Intell.}\ }\textbf {\bibinfo
  {volume} {6}},\ \bibinfo {pages} {1179--1196} (\bibinfo {year}
  {2024})}\BibitemShut {NoStop}%
\bibitem [{\citenamefont {Makarenkov}\ and\ \citenamefont
  {Lamb}(2012)}]{Makarenkov_bifurcations_non-smooth_PhD12}%
  \BibitemOpen
  \bibfield  {author} {\bibinfo {author} {\bibfnamefont {O.}~\bibnamefont
  {Makarenkov}}\ and\ \bibinfo {author} {\bibfnamefont {J.~S.}\ \bibnamefont
  {Lamb}},\ }\bibfield  {title} {\enquote {\bibinfo {title} {Dynamics and
  bifurcations of nonsmooth systems: A survey},}\ }\href
  {https://doi.org/https://doi.org/10.1016/j.physd.2012.08.002} {\bibfield
  {journal} {\bibinfo  {journal} {Physica D: Nonlinear Phenomena}\ }\textbf
  {\bibinfo {volume} {241}},\ \bibinfo {pages} {1826--1844} (\bibinfo {year}
  {2012})},\ \bibinfo {note} {dynamics and Bifurcations of Nonsmooth
  Systems}\BibitemShut {NoStop}%
\bibitem [{\citenamefont {Ermentrout}(2007)}]{xppaut07}%
  \BibitemOpen
  \bibfield  {author} {\bibinfo {author} {\bibfnamefont {B.}~\bibnamefont
  {Ermentrout}},\ }\bibfield  {title} {\enquote {\bibinfo {title} {{XPPAUT}},}\
  }\href@noop {} {\bibfield  {journal} {\bibinfo  {journal} {Scholarpedia}\
  }\textbf {\bibinfo {volume} {2}},\ \bibinfo {pages} {1399} (\bibinfo {year}
  {2007})}\BibitemShut {NoStop}%
\bibitem [{\citenamefont {Mira}\ and\ \citenamefont
  {{\'A}lvarez}(2003)}]{Atencia:2003aa}%
  \BibitemOpen
  \bibinfo {editor} {\bibfnamefont {J.}~\bibnamefont {Mira}}\ and\ \bibinfo
  {editor} {\bibfnamefont {J.}~\bibnamefont {{\'A}lvarez}},\ eds.,\ \href@noop
  {} {\emph {\bibinfo {title} {Computational Methods in Neural Modeling}}}\
  (\bibinfo  {publisher} {Springer Berlin Heidelberg},\ \bibinfo {address}
  {Berlin, Heidelberg},\ \bibinfo {year} {2003})\BibitemShut {NoStop}%
\bibitem [{\citenamefont {Krauth}, \citenamefont {Mezard},\ and\ \citenamefont
  {Nadal}(1988)}]{Krauth_basins_attraction_discrete_perceptron_NN_CS88}%
  \BibitemOpen
  \bibfield  {author} {\bibinfo {author} {\bibfnamefont {W.}~\bibnamefont
  {Krauth}}, \bibinfo {author} {\bibfnamefont {M.}~\bibnamefont {Mezard}},\
  and\ \bibinfo {author} {\bibfnamefont {J.-P.}\ \bibnamefont {Nadal}},\
  }\bibfield  {title} {\enquote {\bibinfo {title} {Basins of attraction in a
  perceptron-like neural network},}\ }\href@noop {} {\bibfield  {journal}
  {\bibinfo  {journal} {Complex Systems}\ }\textbf {\bibinfo {volume} {2}},\
  \bibinfo {pages} {387--408} (\bibinfo {year} {1988})}\BibitemShut {NoStop}%
\bibitem [{\citenamefont {Storkey}\ and\ \citenamefont
  {Valabregue}(1999)}]{Storkey_Hopfield_NN_dsicrete_basins_stat_NN99}%
  \BibitemOpen
  \bibfield  {author} {\bibinfo {author} {\bibfnamefont {A.}~\bibnamefont
  {Storkey}}\ and\ \bibinfo {author} {\bibfnamefont {R.}~\bibnamefont
  {Valabregue}},\ }\bibfield  {title} {\enquote {\bibinfo {title} {The basins
  of attraction of a new hopfield learning rule},}\ }\href
  {https://doi.org/https://doi.org/10.1016/S0893-6080(99)00038-6} {\bibfield
  {journal} {\bibinfo  {journal} {Neural Networks}\ }\textbf {\bibinfo {volume}
  {12}},\ \bibinfo {pages} {869--876} (\bibinfo {year} {1999})}\BibitemShut
  {NoStop}%
\bibitem [{\citenamefont {Davey}\ and\ \citenamefont
  {Hunt}(1999)}]{Davey_Hopfield_NN_discrete_basin_size_99}%
  \BibitemOpen
  \bibfield  {author} {\bibinfo {author} {\bibfnamefont {N.}~\bibnamefont
  {Davey}}\ and\ \bibinfo {author} {\bibfnamefont {S.~P.}\ \bibnamefont
  {Hunt}},\ }\bibfield  {title} {\enquote {\bibinfo {title} {The capacity and
  attractor basins of associative memory models},}\ }in\ \href@noop {} {\emph
  {\bibinfo {booktitle} {Foundations and Tools for Neural Modeling}}},\
  \bibinfo {editor} {edited by\ \bibinfo {editor} {\bibfnamefont
  {J.}~\bibnamefont {Mira}}\ and\ \bibinfo {editor} {\bibfnamefont {J.~V.}\
  \bibnamefont {S{\'a}nchez-Andr{\'e}s}}}\ (\bibinfo  {publisher} {Springer
  Berlin Heidelberg},\ \bibinfo {address} {Berlin, Heidelberg},\ \bibinfo
  {year} {1999})\ pp.\ \bibinfo {pages} {330--339}\BibitemShut {NoStop}%
\bibitem [{\citenamefont {Zhang}\ and\ \citenamefont
  {Zhang}(2008)}]{Zhang_Hopfield_NN_discrete_basin_radius_08}%
  \BibitemOpen
  \bibfield  {author} {\bibinfo {author} {\bibfnamefont {F.}~\bibnamefont
  {Zhang}}\ and\ \bibinfo {author} {\bibfnamefont {X.}~\bibnamefont {Zhang}},\
  }\bibfield  {title} {\enquote {\bibinfo {title} {The average radius of
  attraction basin of hopfield neural networks},}\ }in\ \href@noop {} {\emph
  {\bibinfo {booktitle} {Advances in Neural Networks - ISNN 2008}}},\ \bibinfo
  {editor} {edited by\ \bibinfo {editor} {\bibfnamefont {F.}~\bibnamefont
  {Sun}}, \bibinfo {editor} {\bibfnamefont {J.}~\bibnamefont {Zhang}}, \bibinfo
  {editor} {\bibfnamefont {Y.}~\bibnamefont {Tan}}, \bibinfo {editor}
  {\bibfnamefont {J.}~\bibnamefont {Cao}},\ and\ \bibinfo {editor}
  {\bibfnamefont {W.}~\bibnamefont {Yu}}}\ (\bibinfo  {publisher} {Springer
  Berlin Heidelberg},\ \bibinfo {address} {Berlin, Heidelberg},\ \bibinfo
  {year} {2008})\ pp.\ \bibinfo {pages} {253--258}\BibitemShut {NoStop}%
\bibitem [{\citenamefont {Lin}, \citenamefont {Yeap},\ and\ \citenamefont
  {Kiringa}(2023)}]{Lin_Hopfield_NN_basin_radius_stat_conf23}%
  \BibitemOpen
  \bibfield  {author} {\bibinfo {author} {\bibfnamefont {C.}~\bibnamefont
  {Lin}}, \bibinfo {author} {\bibfnamefont {T.~H.}\ \bibnamefont {Yeap}},\ and\
  \bibinfo {author} {\bibfnamefont {I.}~\bibnamefont {Kiringa}},\ }\bibfield
  {title} {\enquote {\bibinfo {title} {On the basin of attraction and capacity
  of restricted hopfield network as an auto-associative memory},}\ }in\ \href
  {https://doi.org/10.1109/CyberC58899.2023.00033} {\emph {\bibinfo {booktitle}
  {2023 International Conference on Cyber-Enabled Distributed Computing and
  Knowledge Discovery (CyberC)}}}\ (\bibinfo {year} {2023})\ pp.\ \bibinfo
  {pages} {146--154}\BibitemShut {NoStop}%
\bibitem [{\citenamefont {Sampath}\ and\ \citenamefont
  {Srivastava}(2020)}]{Sampath_NN_discrete_basins_PLOS20}%
  \BibitemOpen
  \bibfield  {author} {\bibinfo {author} {\bibfnamefont {S.}~\bibnamefont
  {Sampath}}\ and\ \bibinfo {author} {\bibfnamefont {V.}~\bibnamefont
  {Srivastava}},\ }\bibfield  {title} {\enquote {\bibinfo {title} {On stability
  and associative recall of memories in attractor neural networks},}\ }\href
  {https://api.semanticscholar.org/CorpusID:221788530} {\bibfield  {journal}
  {\bibinfo  {journal} {PLoS ONE}\ }\textbf {\bibinfo {volume} {15}} (\bibinfo
  {year} {2020})}\BibitemShut {NoStop}%
\bibitem [{\citenamefont {Parmananda}\ and\ \citenamefont
  {Eiswirth}(1996)}]{Parmananda96}%
  \BibitemOpen
  \bibfield  {author} {\bibinfo {author} {\bibfnamefont {P.}~\bibnamefont
  {Parmananda}}\ and\ \bibinfo {author} {\bibfnamefont {M.}~\bibnamefont
  {Eiswirth}},\ }\bibfield  {title} {\enquote {\bibinfo {title} {Stabilizing
  unstable fixed points using derivative control},}\ }\href@noop {} {\bibfield
  {journal} {\bibinfo  {journal} {J. Phys. Chem.}\ }\textbf {\bibinfo {volume}
  {100}},\ \bibinfo {pages} {16568--16570} (\bibinfo {year}
  {1996})}\BibitemShut {NoStop}%
\bibitem [{\citenamefont {Claussen}\ \emph {et~al.}(1998)\citenamefont
  {Claussen}, \citenamefont {Mausbach}, \citenamefont {Piel},\ and\
  \citenamefont {Schuster}}]{Claussen98}%
  \BibitemOpen
  \bibfield  {author} {\bibinfo {author} {\bibfnamefont {J.}~\bibnamefont
  {Claussen}}, \bibinfo {author} {\bibfnamefont {T.}~\bibnamefont {Mausbach}},
  \bibinfo {author} {\bibfnamefont {A.}~\bibnamefont {Piel}},\ and\ \bibinfo
  {author} {\bibfnamefont {H.}~\bibnamefont {Schuster}},\ }\bibfield  {title}
  {\enquote {\bibinfo {title} {Memory difference control of unknown unstable
  fixed points: Drifting parameter conditions and delayed measurement},}\
  }\href@noop {} {\bibfield  {journal} {\bibinfo  {journal} {Phys. Rev. E}\
  }\textbf {\bibinfo {volume} {58}},\ \bibinfo {pages} {7256--7260} (\bibinfo
  {year} {1998})}\BibitemShut {NoStop}%
\bibitem [{\citenamefont {Epureanu}\ and\ \citenamefont
  {Dowell}(2000)}]{Epureanu00}%
  \BibitemOpen
  \bibfield  {author} {\bibinfo {author} {\bibfnamefont {B.}~\bibnamefont
  {Epureanu}}\ and\ \bibinfo {author} {\bibfnamefont {E.}~\bibnamefont
  {Dowell}},\ }\bibfield  {title} {\enquote {\bibinfo {title} {Optimal
  multi-dimensional {OGY} controller},}\ }\href@noop {} {\bibfield  {journal}
  {\bibinfo  {journal} {Physica D}\ }\textbf {\bibinfo {volume} {139}},\
  \bibinfo {pages} {87--96} (\bibinfo {year} {2000})}\BibitemShut {NoStop}%
\bibitem [{\citenamefont {Ahlborn}\ and\ \citenamefont
  {Parlitz}(2004)}]{Ahlborn04}%
  \BibitemOpen
  \bibfield  {author} {\bibinfo {author} {\bibfnamefont {A.}~\bibnamefont
  {Ahlborn}}\ and\ \bibinfo {author} {\bibfnamefont {U.}~\bibnamefont
  {Parlitz}},\ }\bibfield  {title} {\enquote {\bibinfo {title} {Stabilizing
  unstable steady states using multiple delay feedback control},}\ }\href@noop
  {} {\bibfield  {journal} {\bibinfo  {journal} {Phys. Rev. Lett.}\ }\textbf
  {\bibinfo {volume} {93}} (\bibinfo {year} {2004})}\BibitemShut {NoStop}%
\bibitem [{\citenamefont {Abed}\ and\ \citenamefont {Fu}(1986)}]{Abed86}%
  \BibitemOpen
  \bibfield  {author} {\bibinfo {author} {\bibfnamefont {E.}~\bibnamefont
  {Abed}}\ and\ \bibinfo {author} {\bibfnamefont {J.}~\bibnamefont {Fu}},\
  }\bibfield  {title} {\enquote {\bibinfo {title} {Local feedback stabilization
  and bifurcations control, {I}. {H}opf bifurcation},}\ }\href@noop {}
  {\bibfield  {journal} {\bibinfo  {journal} {Syst. Control Lett.}\ }\textbf
  {\bibinfo {volume} {7}},\ \bibinfo {pages} {11--17} (\bibinfo {year}
  {1986})}\BibitemShut {NoStop}%
\bibitem [{\citenamefont {Abed}\ and\ \citenamefont {Fu}(1987)}]{Abed87}%
  \BibitemOpen
  \bibfield  {author} {\bibinfo {author} {\bibfnamefont {E.}~\bibnamefont
  {Abed}}\ and\ \bibinfo {author} {\bibfnamefont {J.}~\bibnamefont {Fu}},\
  }\bibfield  {title} {\enquote {\bibinfo {title} {Local feedback stabilization
  and bifurcations control, {II}. {S}tationary bifurcation},}\ }\href@noop {}
  {\bibfield  {journal} {\bibinfo  {journal} {Syst. Control Lett.}\ }\textbf
  {\bibinfo {volume} {8}},\ \bibinfo {pages} {467--473} (\bibinfo {year}
  {1987})}\BibitemShut {NoStop}%
\bibitem [{\citenamefont {Chen}, \citenamefont {Moiola},\ and\ \citenamefont
  {Wang}(2000)}]{Chen00}%
  \BibitemOpen
  \bibfield  {author} {\bibinfo {author} {\bibfnamefont {G.}~\bibnamefont
  {Chen}}, \bibinfo {author} {\bibfnamefont {J.}~\bibnamefont {Moiola}},\ and\
  \bibinfo {author} {\bibfnamefont {H.}~\bibnamefont {Wang}},\ }\bibfield
  {title} {\enquote {\bibinfo {title} {Bifurcation control: Theories, methods,
  and applications},}\ }\href {https://doi.org/10.1142/S0218127400000360}
  {\bibfield  {journal} {\bibinfo  {journal} {Int. J. Bifurcat. Chaos}\
  }\textbf {\bibinfo {volume} {10}},\ \bibinfo {pages} {511--548} (\bibinfo
  {year} {2000})}\BibitemShut {NoStop}%
\bibitem [{\citenamefont {Ramsauer}\ \emph {et~al.}(2020)\citenamefont
  {Ramsauer}, \citenamefont {Schafl}, \citenamefont {Lehner}, \citenamefont
  {Seidl}, \citenamefont {Widrich}, \citenamefont {Gruber}, \citenamefont
  {Holzleitner}, \citenamefont {Pavlovi'c}, \citenamefont {Sandve},
  \citenamefont {Greiff}, \citenamefont {Kreil}, \citenamefont {Kopp},
  \citenamefont {Klambauer}, \citenamefont {Brandstetter},\ and\ \citenamefont
  {Hochreiter}}]{Ramsauer_Hopfield_NN_continuous_memory_capacity_arxiv21}%
  \BibitemOpen
  \bibfield  {author} {\bibinfo {author} {\bibfnamefont {H.}~\bibnamefont
  {Ramsauer}}, \bibinfo {author} {\bibfnamefont {B.}~\bibnamefont {Schafl}},
  \bibinfo {author} {\bibfnamefont {J.}~\bibnamefont {Lehner}}, \bibinfo
  {author} {\bibfnamefont {P.}~\bibnamefont {Seidl}}, \bibinfo {author}
  {\bibfnamefont {M.}~\bibnamefont {Widrich}}, \bibinfo {author} {\bibfnamefont
  {L.}~\bibnamefont {Gruber}}, \bibinfo {author} {\bibfnamefont
  {M.}~\bibnamefont {Holzleitner}}, \bibinfo {author} {\bibfnamefont
  {M.}~\bibnamefont {Pavlovi'c}}, \bibinfo {author} {\bibfnamefont {G.~K.~F.}\
  \bibnamefont {Sandve}}, \bibinfo {author} {\bibfnamefont {V.}~\bibnamefont
  {Greiff}}, \bibinfo {author} {\bibfnamefont {D.~P.}\ \bibnamefont {Kreil}},
  \bibinfo {author} {\bibfnamefont {M.}~\bibnamefont {Kopp}}, \bibinfo {author}
  {\bibfnamefont {G.}~\bibnamefont {Klambauer}}, \bibinfo {author}
  {\bibfnamefont {J.}~\bibnamefont {Brandstetter}},\ and\ \bibinfo {author}
  {\bibfnamefont {S.}~\bibnamefont {Hochreiter}},\ }\bibfield  {title}
  {\enquote {\bibinfo {title} {Hopfield networks is all you need},}\ }\href
  {https://api.semanticscholar.org/CorpusID:220968978} {\bibfield  {journal}
  {\bibinfo  {journal} {ArXiv}\ }\textbf {\bibinfo {volume} {abs/2008.02217}}
  (\bibinfo {year} {2020})}\BibitemShut {NoStop}%
\end{thebibliography}

%

\end{document}